\documentclass[final,11pt,times]{elsarticle}
\usepackage{amsmath}
\usepackage{amssymb}
\usepackage{upgreek}
\usepackage{bm}
\usepackage[colorlinks=true,linkcolor=red,citecolor=blue,urlcolor=blue]{hyperref}
\usepackage{cleveref}
\usepackage{booktabs}
\usepackage{subfigure}
\usepackage{float}
\usepackage{multirow}
\usepackage{fullpage}
\usepackage[ruled, linesnumbered]{algorithm2e}
\usepackage{amsthm}
\usepackage{arydshln}

\allowdisplaybreaks[4]
\newtheorem{theorem}{Theorem}

\newtheorem{remark}{Remark}

\journal{Journal of Computational Physics}

\begin{document}

\begin{frontmatter}
\title{A Multiple Transferable Neural Network Method with Domain Decomposition for Elliptic Interface Problems}

\author[a]{Tianzheng Lu}
\ead{lutianzheng@buaa.edu.cn}
\author[b]{Lili Ju}
\ead{ju@math.sc.edu}
\author[a]{Liyong Zhu\corref{cor1}}
\ead{liyongzhu@buaa.edu.cn}

\address[a]{School of Mathematical Sciences, Beihang University, Beijing 100191, China}
\address[b]{Department of Mathematics, University of South Carolina, Columbia, SC 29208, USA}
\cortext[cor1]{Corresponding author.}

\begin{abstract}
The transferable neural network (TransNet) is a two-layer shallow neural network with pre-determined and uniformly distributed neurons in the hidden layer, and the least-squares solvers can be particularly used to compute the parameters of its output layer when applied to the solution of partial differential equations. In this paper, we integrate the TransNet technique with the nonoverlapping domain decomposition and the interface conditions to develop a novel multiple transferable neural network (Multi-TransNet) method for solving elliptic interface problems, which typically contain discontinuities in both solutions and their derivatives across interfaces. We first propose an empirical formula for the TransNet to characterize the relationship between the radius of the domain-covering ball, the number of hidden-layer neurons, and the optimal neuron shape. In the Multi-TransNet method, we assign each subdomain one distinct TransNet with an adaptively determined number of hidden-layer neurons to maintain the globally uniform neuron distribution across the entire computational domain, and then unite all the subdomain TransNets together by incorporating the interface condition terms into the loss function. The empirical formula is also extended to the Multi-TransNet and further employed to estimate appropriate neuron shapes for the subdomain TransNets, greatly reducing the parameter tuning cost. Additionally, we propose a normalization approach to adaptively select the weighting parameters for the terms in the loss function. Ablation studies and extensive experiments with comparison tests on different types of elliptic interface problems with low to high contrast diffusion coefficients in two and three dimensions are carried out to numerically demonstrate the superior accuracy, efficiency, and robustness of the proposed Multi-TransNet method.
\end{abstract}


\begin{keyword}
Two-layer  neural networks  \sep elliptic interface problems \sep transferability \sep domain decomposition \sep neuron shape
\end{keyword}

\end{frontmatter}

\section{Introduction}\label{sec:intro}
\begin{figure}[!ht]
	\centering
	\begin{minipage}[t]{.3\textwidth}
		\centering
		\includegraphics[width=\linewidth]{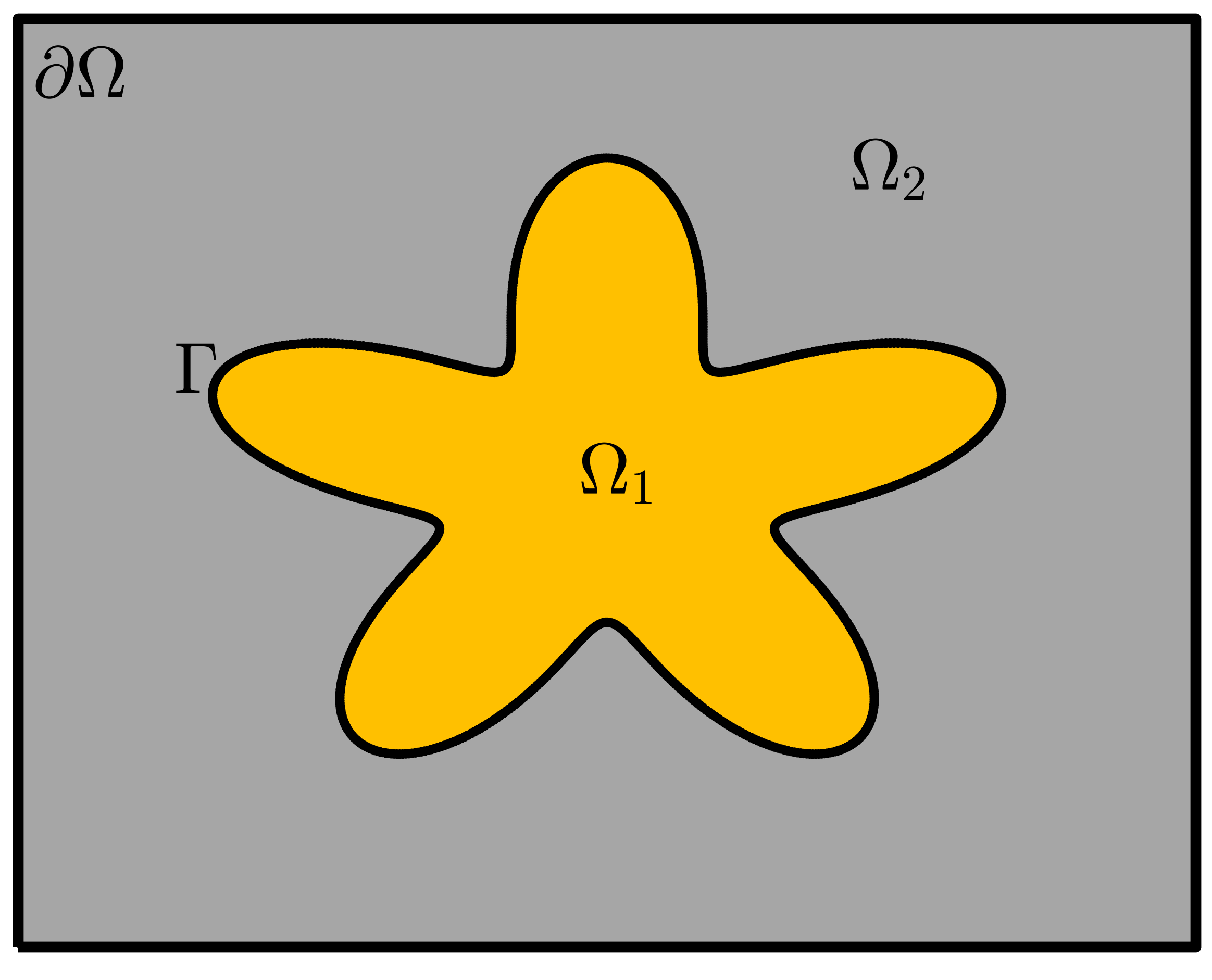}	
	\end{minipage}
	\hspace{2cm}
    \begin{minipage}[t]{.32\textwidth}
		\centering
		\includegraphics[width=\linewidth]{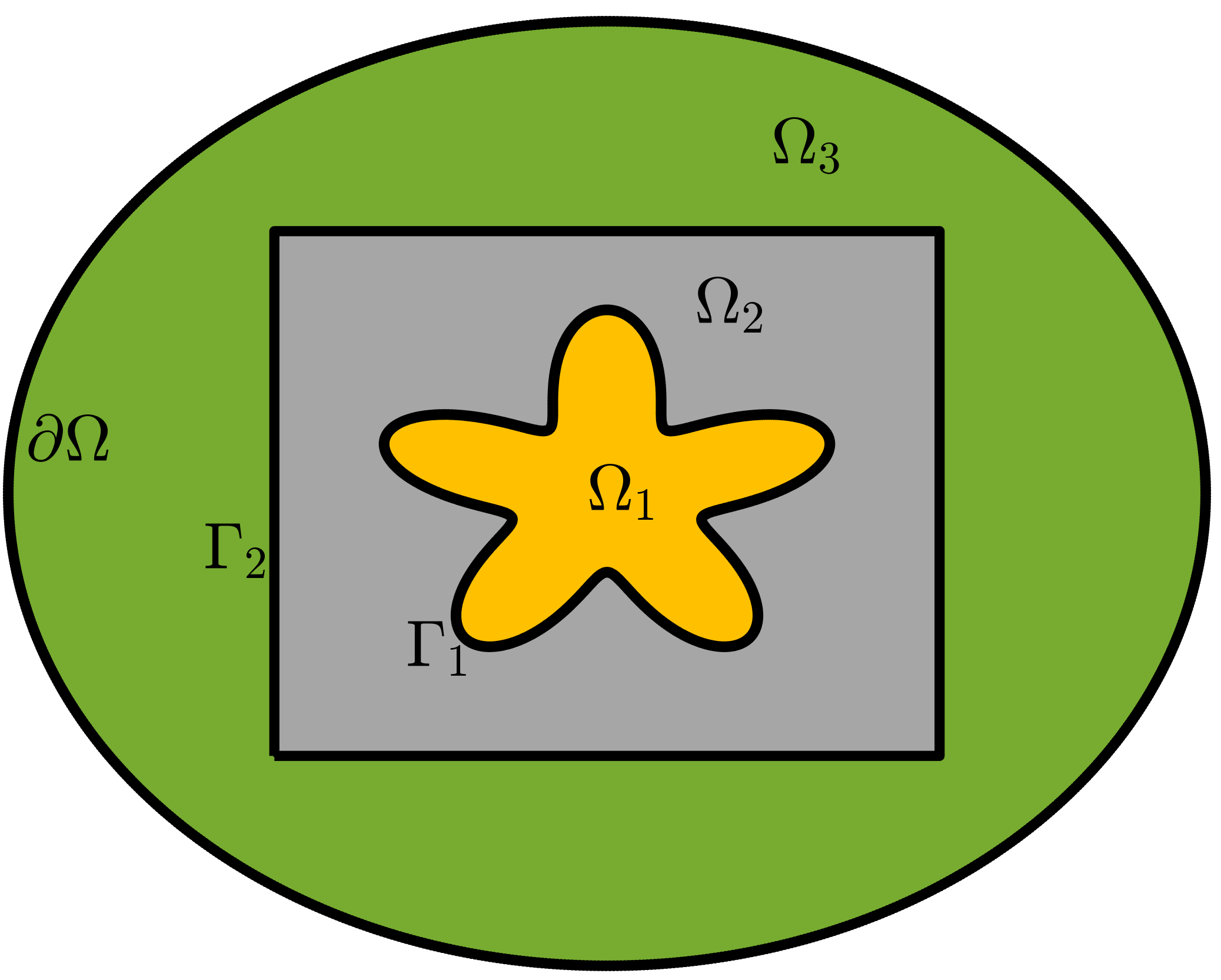}	
	\end{minipage}
    \vspace{-2mm}
    \caption{Some sample domains with interfaces. Left: $K=2$ with $\Gamma = \overline\Omega_1\cap\overline\Omega_2$; Right: $K=3$ with
    $\Gamma_1 = \overline\Omega_1\cap\overline\Omega_2$ and $\Gamma_2 = \overline\Omega_2\cap\overline\Omega_3$.}
    \label{fig:IP-domain}
\end{figure}

Interface problems arise in numerous physical applications, such as fluid mechanics \cite{hou1997hybrid, gibou2019sharp}, composite materials \cite{greengard1994numerical, bochkov2023numerical}, biological sciences \cite{peskin1977numerical, xia2014multiscale, egan2018fast} and electromagnetics \cite{chen2009adaptive, costabel1999singularities}, among others.
In this paper, we consider the following elliptic interface problem located in an open bounded, connected  domain $\Omega\in \mathbb{R}^d$:
\begin{equation}\label{equ:IP}
\left\{\begin{aligned}
	\mathcal{L}\left(u\right)  = f,&\quad \boldsymbol{x} \in \cup_{k=1}^{K}\Omega_k, \\
	[u]  =h_1,&\quad \boldsymbol{x} \in \Gamma_{i,j},\ \ i\ne j,\\
	[\mathcal{J}\left(u\right) \cdot \boldsymbol{n}_{i,j}]  =h_2,&\quad \boldsymbol{x} \in \Gamma_{i,j},\ \ i\ne j,\\\
	\mathcal{B}\left(u\right)  =g,&\quad \boldsymbol{x} \in \partial \Omega,
\end{aligned}\right.
\end{equation}
where $\Omega$ is partitioned into $K$ open connected subdomains $\{\Omega_1, \Omega_2,\dots,\Omega_K\}$ by the closed interfaces $\{\Gamma_{i,j}=\overline\Omega_i\cap\overline\Omega_j\ne\emptyset\;|\; i\ne j\}$, as illustrated in \autoref{fig:IP-domain}. The operators $\mathcal{L}$ and $\mathcal{J}$ take some differential forms in subdomains and interfaces, while $\mathcal{B}$ denotes the boundary operator defined on $\partial \Omega$, which are all assumed to be linear in this work.
The notation $[v]$ indicates the jump of $v$ across the interface $\Gamma_{i,j}$, i.e., $[v]=\left.v\right|_{\Omega_i} - \left.v\right|_{\Omega_j}$, and the vector $\boldsymbol{n}_{i,j}$ represents the unit outward normal vector on $\Gamma_{i,j}$ from $\Omega_i$ to $\Omega_j$. 

The low regularity of solutions across interfaces, coupled with the complex geometry of these interfaces, often leads to accuracy loss when applying standard numerical methods to solve the elliptic interface problem \eqref{equ:IP}. To address this, body-fitted meshes are introduced to ensure optimal or near-optimal convergence rates \cite{babuvska1970finite, chen1998finite}. However, generating robust grids is usually time-consuming and inefficient, particularly in the case of large deformations. Consequently, researchers have turned to modifying standard numerical methods on structured grids, a line of work that can be dated back to the immersed boundary method \cite{peskin1977numerical, peskin2002immersed}.
Since then, a variety of interface-capturing methods have emerged, which can generally be categorized into explicit and implicit approaches according to how they handle interfaces. Explicit methods, such as the front-tracking method \cite{unverdi1992front}, involve explicitly tracking interfaces. While these methods offer higher accuracy in capturing interface details, they may be less efficient in dealing with changes in interface topology. In contrast, implicit approaches are more flexible in accommodating such topological changes and are often preferred in practice for complex scenarios. Notable examples include the volume of fluid method \cite{hirt1981volume}, which tracks the volume fraction of phases within computational cells, inherently ensuring mass conservation across the domain; the level-set method \cite{osher1988fronts}, which employs a signed distance function to implicitly represent the interface as the zero level-set of a scalar function, allowing for more natural handling of complex interface geometries; the immersed interface method \cite{leveque1994immersed, li2006immersed}, which modifies standard finite difference stencils near interfaces by adding correction terms to account for jumps in the solution or its derivatives, maintaining optimal accuracy up to interfaces; and the ghost fluid method \cite{fedkiw1999non, liu2000boundary}, which defines ghost cells adjacent to interfaces to enforce interface conditions, enabling standard numerical schemes to be applied across discontinuities. Additionally, a variant of the ghost fluid method \cite{egan2020xgfm} ensures flux convergence by modifying the interface condition for fluxes in the ghost fluid method. We shall not provide an exhaustive discussion of all related traditional numerical methods here (see \cite{gibou2013high, gibou2018review} for more comprehensive reviews).

In recent years, with remarkable success of neural networks across various fields \cite{lecun2015deep}, neural network-based numerical methods for scientific computing, especially solving partial differential equations (PDEs) have emerged. Methods such as the Deep Ritz Method (DRM) \cite{yu2018deep}, the Deep Galerkin Method (DGM) \cite{sirignano2018dgm}, and the Physics-Informed Neural Networks (PINNs) \cite{raissi2019physics} have gained significant attention due to their mesh-free nature. In fact, this mesh-free property is especially advantageous for addressing problems involving complex geometric domains, such as interface problems. The elliptic interface problem was tackled using a deep neural network approach based on least-square functionals of the first-order system in \cite{cai2020deep}. DRM was also employed to address elliptic interface problems with high-contrast discontinuity coefficients in \cite{wang2020mesh}. 
Recently, localized neural network methods, including those utilizing domain decomposition techniques \cite{jagtap2020extended, li2020deep, chen2022bridging}, have attracted increasing interest. A piecewise deep neural network method was applied to elliptic interface problems in \cite{he2022mesh}, where it was numerically demonstrated that the method can accurately solve problems with complex-shaped interfaces. In \cite{wu2022inn}, the Interfaced Neural Network (INN) method, which utilizes multiple neural networks, was proposed, and the multiple-gradient descent approach was extended to adaptively adjust the weighting parameters in the loss function, thereby improving the robustness of solutions for elliptic interface problems. In \cite{mistani2023jax}, to achieve controllable accuracy and convergence, the authors integrated an advanced level-set based finite volume numerical scheme \cite{bochkov2020solving} into two deep neural networks. This hybrid method was applied in parallel to solve three-dimensional elliptic interface problems involving convoluted geometries.

In practice, to achieve better approximation power, the neural networks mentioned above are typically designed with deep architectures. However, this often leads to challenging optimization problems. Given the limitations of current optimization techniques, it is often difficult to significantly reduce optimization errors, hindering the attainment of high accuracy. Based on the universal approximation theorem for two-layer (i.e., single-hidden-layer) neural networks \cite{huang1998upper}, a viable approach is to employ such a shallow network structure and leverage high-order optimizers to efficiently lower optimization errors in practice. In \cite{hu2022discontinuity, tseng2023cusp}, a high-order full-batch optimizer was introduced to train a shallow neural network with augmented input for elliptic interface problems featuring cusps or discontinuities, leading to a significant improvement in accuracy. However, high-order optimizers are generally unstable and computationally prohibitive for large-scale training. 
An alternative approach is to bypass the challenging optimization problem, altogether by randomly initializing and fixing the parameters of the hidden layers, leaving only the output layer trainable. This results in a significantly simplified optimization problem, reducing it to a least squares problem that depends only the parameters of the output layer. The linearity or nonlinearity of the problem is fully determined by the nature of PDEs. Furthermore, this approach can be efficiently solved using well-established linear or nonlinear least-squares techniques, eliminating the need of gradient-based optimizers. 
The local Extreme Learning Machine (locELM) \cite{dong2021local} and the Random Feature Method (RFM) \cite{chen2022bridging} have been proposed for solving PDEs. Recently, both methods have been applied to elliptic interface problems in \cite{li2023local} and \cite{chi2024random}, demonstrating high accuracy and efficiency, often on a par with or even surpassing traditional numerical methods. However, the fixed parameters of hidden layers in both locELM and RFM lack interpretability, which typically leads to blind parameter selection even impractical manual adjustments. The recently proposed Transferable Neural Network (TransNet) \cite{zhang2024transferable} addresses this issue by developing a two-layer neural network with pre-determined, uniformly distributed, shape-shared neurons in the hidden layer, offering both intuitive interpretability and excellent transferability for solving PDEs.

In this paper, we develop a novel {\em Multiple Transferable Neural Network} (Multi-TransNet) method for solving elliptic interface problems. This approach leverages nonoverlapping domain decomposition, tailoring multiple distinct transferable neural networks to different subdomains based on their specific features and integrating these subdomain TransNets through interface conditions incorporated into the loss function.
The main contributions of our work are summarized as follows.
\begin{itemize}
    \item We extend the property of uniform distribution of hidden-layer neurons in TransNet to the global case in Multi-TransNet, enabling the adaptive assignment of the number of neurons for each subdomain TransNet, thereby enhancing the overall transferability.\vspace{-0.1cm}
    \item We propose an empirical formula for TransNet to characterize the relationship between some of its key parameters, which is then extended to Multi-TransNet and employed to predict appropriate neuron shapes for subdomain TransNets, significantly reducing the parameter tuning cost.\vspace{-0.1cm}
   \item We  propose a normalization approach to adaptively select the weighting parameters for the terms in the loss function, improving the accuracy and robustness  of the Multi-TransNet method.\vspace{-0.1cm}
    \item We conduct extensive ablation studies to confirm the effectiveness of the proposed strategies, and perform abundant experiments with comparative tests  on various types of elliptic interface problems in two and three dimensions to demonstrate the superior accuracy, efficiency, and robustness of the proposed Multi-TransNet method.
\end{itemize}

The rest of the paper is organized as follows. In \autoref{sec:revisit-TransNet}, we briefly revisit the TransNet method for solving general PDEs. In \autoref{sec:TransNet-auto-tune-gamma}, we propose an efficient empirical formula-based strategy  for automatically select an appropriate  shape parameter of TransNet. Subsequently, the Multi-TransNet method for elliptic interface problems is proposed and comprehensively discussed in \autoref{sec:Multi-TransNet}. Extensive ablation studies and comparison tests on various types of two- and three-dimensional elliptic interface problems are presented  in \autoref{sec:num-experim} to numerically  demonstrate the outstanding performance of the proposed Multi-TransNet method, followed by some concluding remarks in \autoref{sec:conclusions}.

\section{Review of the transferable neural network method}\label{sec:revisit-TransNet}
The TransNet  developed in \cite{zhang2024transferable}  is a two-layer neural network method  for solving PDEs of the following form
\begin{equation}\label{equ:pde}
\left\{
\begin{aligned}
 & \mathcal{L}(u) = f,\;\;  \quad \bm x \in \Omega,\\
 & \mathcal{B}(u) =  g,\;\;  \quad \bm x \in \partial \Omega.
\end{aligned}
\right.
\end{equation}
The key ingredient is the so-called neural feature space, denoted by $\mathcal{P}_{\mathrm{NN}}$, expanded by a group of neural basis functions (i.e., the neurons of the hidden layer):
\begin{equation}\label{equ:neural-feature-space}
	\mathcal{P}_{\mathrm{NN}}=\operatorname{span}\left\{1, \sigma\left(\boldsymbol{w}_1^{\top} \boldsymbol{x}+b_1\right), \ldots, \sigma\left(\boldsymbol{w}_M^{\top} \boldsymbol{x}+b_M\right)\right\},
\end{equation}
where $\sigma(\cdot)$ is the activation function, and $\boldsymbol{w}_m\in \mathbb{R}^d, b_m\in \mathbb{R}$ is the weight and bias of the $m$-th hidden-layer neuron, respectively. The TransNet solution $u_\mathrm{NN}\in \mathcal{P}_\mathrm{NN}$ is represented by a linear combination of the neural basis functions
\begin{equation}\label{equ:uNN-NotTranslate}
	u_{\mathrm{NN}}(\boldsymbol{x})=\sum_{m=1}^M \alpha_m \sigma\left(\boldsymbol{w}_m^{\top} \boldsymbol{x}+b_m\right)+\alpha_0,
\end{equation}
where $\alpha_0, \alpha_1, \alpha_2, \dots, \alpha_M$ are the parameters of the output layer in TransNet (one may directly take $\alpha_0=0$ in practice). 
\emph{The main idea of TransNet is that the neurons of its hidden layer will be pre-determined in a special sense of uniform distribution from the pure approximation point of view before applied to the solution of PDEs}. Below we briefly review the process of the TransNet method for solving the PDE problem \eqref{equ:pde}.

\subsection{Geometrization of the neural feature space}
Inspired by activation patterns of ReLU networks, the hidden-layer neurons of TransNet are geometrized based on a  {re-parameterization} of the neural feature space $\mathcal{P}_\mathrm{NN}$ without using any PDE information, and then are  naturally connected with the computational domain $\Omega$. Specifically, a hidden-layer neuron represented by a piecewise-linear function in a ReLU network can be viewed as a \emph{partition hyperplane}, i.e.,
\begin{equation}\label{equ:origi-parti-hyperplane}
	\boldsymbol{w}_m^\top \boldsymbol{x} + b_m=0,
\end{equation}
in a geometric space. It can be rewritten into
\begin{equation}\label{equ:reparam-parti-hyperplane}
    \gamma_m \boldsymbol{a}_m^\top \left(\boldsymbol{x} + r_m \boldsymbol{a}_m\right)=0,
\end{equation}
similar to the point-normal form of a hyperplane equation in $\mathbb{R}^d$, in which both the unit normal vector $\boldsymbol{a}_m \in \mathbb{R}^d$ of the partition hyperplane \eqref{equ:origi-parti-hyperplane} and the scalar $r_m \in \mathbb{R}$ representing the distance between the origin and the partition hyperplane \eqref{equ:origi-parti-hyperplane} are collectively referred to as the \emph{location parameter}, while the scalar $\gamma_m\in \mathbb{R}_+$ is named as the \emph{shape parameter} owing to its association with the geometric shape characteristic of  $\sigma\left(\gamma_m\left(\boldsymbol{a}_m^{\top} \boldsymbol{x}+r_m\right)\right)$ (see Figure 1 in \cite{zhang2024transferable} for their visualization). Evidently, the relationship between the original parameters $\boldsymbol{w}_m, b_m$ and the current ones $\boldsymbol{a}_m, r_m, \gamma_m$ is given by
\begin{equation}\label{equ:relationship}
		\boldsymbol{w}_m = \gamma_m \boldsymbol{a}_m, \quad
		b_m = \gamma_m r_m.
\end{equation}
Therefore, the TransNet solution \eqref{equ:uNN-NotTranslate} is rewritten as
\begin{equation}\label{equ:uNN-translate}
	u_{\mathrm{NN}}(\boldsymbol{x})
        =\sum_{m=1}^M \alpha_m \sigma\left(\gamma_m\left(\boldsymbol{a}_m^{\top} \boldsymbol{x}+r_m\right)\right)+\alpha_0.
\end{equation} 
The next step of TransNet is to obtain the location parameter $(\boldsymbol{a}_m, r_m)$ and the shape parameter $\gamma_m$ of the hidden-layer neuron, i.e., the \emph{pre-training}.

\subsection{Generating the hidden-layer neuron location --- uniform neuron distribution}\label{ssec:uniform-distri-arbi-ball}
In a ReLU network, a hidden-layer neuron can be regarded as a partition hyperplane \eqref{equ:origi-parti-hyperplane} to some extent, and a region intersected by multiple such partition hyperplanes can define a linear piece. Obviously, adding the number of linear pieces is the most straightforward way to enhance the approximation power of the network for solving PDEs. Furthermore, intuitively, the location of linear pieces should have a significant influence on the generalization of the network for solving PDEs defined in different domains. For these reasons, a concept of  \emph{uniform neuron distribution} is designed and related construction algorithm is rigorously proven in \cite{zhang2024transferable}. 
\begin{theorem}[Uniform neuron distribution in the unit ball $B_{1}(\boldsymbol{0})$]\label{thm:uniform-neuron-distribution}
	Given a set of  $M$ partition hyperplanes of $\mathbb{R}^d$ defined by
	\begin{equation*}
		\boldsymbol{a}_m^\top \boldsymbol{x} + r_m = 0,\quad m=1,2,\dots M.
	\end{equation*}
	If $\left\{\boldsymbol{a}_m\right\}_{m=1}^M$ are i.i.d. and uniformly distributed on the d-dimensional unit sphere, and $\left\{r_m\right\}_{m=1}^M$ are i.i.d. and uniformly distributed in $[0,1]$, then for a fixed $\uptau \in(0,1)$,
	\begin{equation*}
	\mathbb{E}\left[D_M^\uptau(\boldsymbol{x})\right]=\uptau,\quad \forall\, \|\boldsymbol{x}\|_2 \leqslant 1-\uptau,
	\end{equation*}
	where $D_M^\uptau(\boldsymbol{x})$ is the density function of the neurons defined by
	\begin{equation*}
	D_M^\uptau(\boldsymbol{x})=\frac{1}{M} \sum_{m=1}^M \chi_{\left\{d_m(\boldsymbol{x})<\uptau\right\}}(\boldsymbol{x}), 
	\end{equation*}
	with  $\chi_{\left\{d_m(\boldsymbol{x})<\uptau\right\}}(\boldsymbol{x})$ denoting the indicator function  whether the distance between $\boldsymbol{x}$ and the $m$-th partition hyperplane is smaller than $\uptau$.
\end{theorem}

Note that in practice, $\left\{\boldsymbol{a}_m\right\}_{m=1}^M$ satisfying the requirement in \autoref{thm:uniform-neuron-distribution} can be obtained by sampling from the $d$-dimensional standard Gaussian distribution in the Cartesian coordinate system and then normalizing the samples to unit vectors. {\autoref{thm:uniform-neuron-distribution} limits the domain of computation to a unit ball,  but one can naturally generalize it to the case of a ball $B_R(\boldsymbol{x}_c)$ centered at the point $\boldsymbol{x}_c\in \mathbb{R}^d$ with a radius $R\in \mathbb{R}_+$ with  an affine (scaling and translation) transformation}, and thus obtain the following  result.

\begin{theorem}[Uniform neuron distribution in the ball $B_{R}(\boldsymbol{x}_c)$]\label{thm:general-uniform-neuron-distribution}
	Given a set of  $M$ partition hyperplanes of $\mathbb{R}^d$ defined by
	\begin{equation*}
		\boldsymbol{a}_m^\top (\boldsymbol{x} - \boldsymbol{x}_c) + r_m = 0, \quad m=1,2,\dots M.
	\end{equation*}
If $\left\{\boldsymbol{a}_m\right\}_{m=1}^M$ are i.i.d. and uniformly distributed on the d-dimensional unit sphere, and $\left\{r_m\right\}_{m=1}^M$ are i.i.d. and uniformly distributed in $[0,R]$, then for a fixed $\uptau \in(0,R)$,
	\begin{equation}\label{res2}
	\mathbb{E}\left[D_M^\uptau(\boldsymbol{x})\right]=\frac{\uptau}{R},\quad \forall\, \|\boldsymbol{x} - \boldsymbol{x}_c\|_2 \leqslant R-\uptau.
	\end{equation}
\end{theorem}

\begin{figure}[!ht]
    \centering
	\includegraphics[width=.7\textwidth]{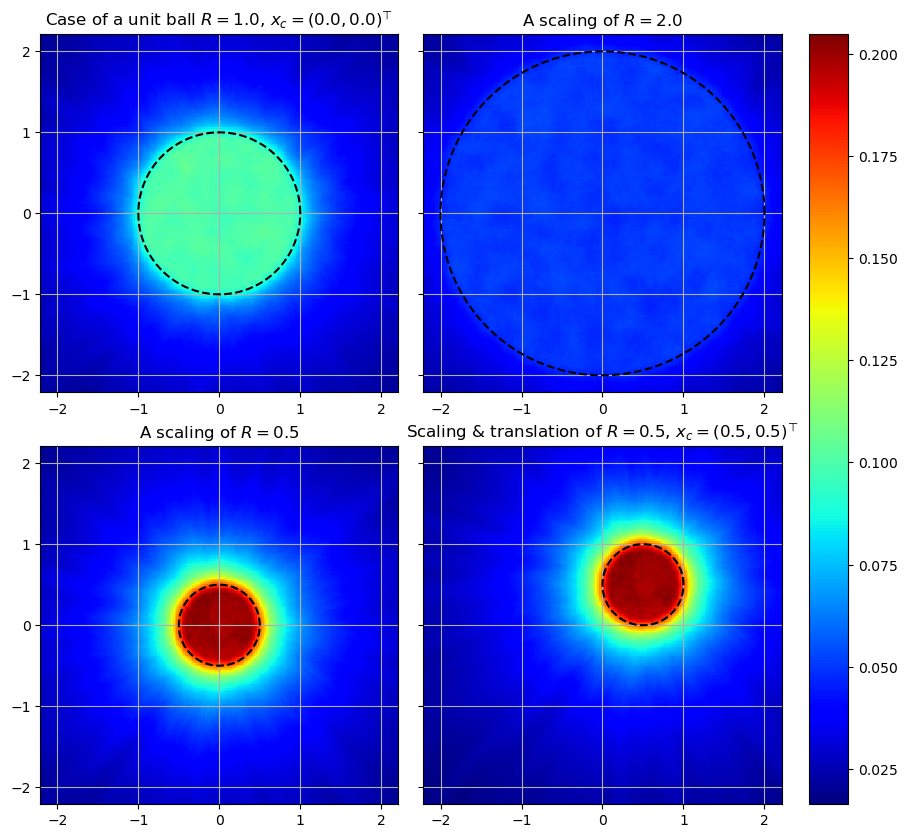}
	\vspace{-2mm}
	\caption{Illustration of the density distribution of the partition hyperplanes of hidden-layer neurons from TransNet with different $\boldsymbol{x}_c$ and $R$. Here $M=30,000$ and $\uptau=0.1$.
	Top-left to top-right: with a scaling of factor of 2.0, top-left to bottom-left:  with a scaling of  factor of 0.5,
	and top-left to bottom-right: with a scaling of  factor of 0.5 and a translation of the center $(0.5,0.5)$.}
	\label{fig:general-uniform-neuron-distribution}
\end{figure}

For the sake of geometric intuition, we illustrate the uniform neuron distribution with some examples in \autoref{fig:general-uniform-neuron-distribution}. It is clearly observed that the scaling of the hidden-layer neurons is achieved from the case of $R=1.0$ and $\boldsymbol{x}_c=(0,0)^\top$  (top-left) to the case of $R=2$ (top-right) and the case of $R=0.5$ (bottom-left), and the corresponding density of hidden-layer neurons doubles or halves. The hidden-layer neurons of bottom-left further undergo a translation from the center $(0, 0)$ to another center $(0.5, 0.5)$ (bottom-right), and the density remains unchanged. 

\begin{remark}[The selection of parameters $\boldsymbol{x}_c$ and $R$ for a general domain $\Omega$]\label{rem:translatation-and-scaling}
For a given domain $\Omega$,  the center of the ball $\boldsymbol{x}_c$ should be selected to be approximately the center of the domain $\Omega$ and 
the radius of the ball $R$ should be chosen to make the ball $B_{R}(\boldsymbol{x}_c)$ slightly over-cover the domain $\Omega$, to enhance the efficacy and density of the hidden-layer neurons within the domain.
\end{remark}

\subsection{Tuning the hidden-layer neuron shape and solving the PDE}\label{ssec:TransNet-GRFs}
In order to strive for simplicity and generalization, all hidden-layer neurons of TransNet are assumed to share the same shape parameter in \cite{zhang2024transferable}, i.e., $\gamma_m = \gamma$ for $\ m=1,2,\dots,M$. The Gaussian Random Fields (GRFs) are introduced to generate auxiliary functions for pre-tuning $\gamma$ based on the function approximation ability, and the grid search is used for  finding an optimal value of  $\gamma$ which minimizes the mean approximation errors over the whole set of auxiliary functions. The entire tuning process is \emph{offline} and only need to solve some linear least squares problems. The neural feature space \eqref{equ:neural-feature-space} is then completely constructed without using any PDE information, and the TransNet solution \eqref{equ:uNN-translate} is reduced to 
\begin{equation}\label{equ:uNN-final}
	u_{\mathrm{NN}}(\boldsymbol{x})
        =\sum_{m=1}^M \alpha_m \sigma\left(\gamma\left(\boldsymbol{a}_m^{\top} \left(\boldsymbol{x}-\boldsymbol{x}_c\right)+r_m\right)\right)+\alpha_0.
\end{equation} 

The final step of TransNet for solving the problem \eqref{equ:pde}  is to obtain the parameters $\boldsymbol{\alpha} = \{\alpha_0, \alpha_1, \dots, \alpha_M\}$ of the output layer.
This can be done by substituting the TransNet solution \eqref{equ:uNN-final} into \eqref{equ:pde} and  minimizing the following  physics-informed loss function
\begin{equation}\label{equ:loss}
{\mathrm{Loss_{TN}}}\,({\boldsymbol{\alpha}}) = \lambda_L \|\mathcal{L}(u_{\mathrm{NN}}(\boldsymbol{x}))\ - f(\boldsymbol{x})\|_2^2 + \lambda_B \|\mathcal{B}(u_{\mathrm{NN}}(\bm x)) - g(\bm x)\|_2^2
\end{equation}
over a set of training/collocation points, where $\lambda_L$ and $\lambda_B$ are some positive weighting parameters. This is essentially solving a linear (or nonlinear if the operators $\mathcal{L}$ and/or $\mathcal{B}$ in \eqref{equ:pde} are nonlinear) least squares problem, 
which can be efficiently done with preferred QR factorization related techniques.

\section{An efficient  empirical formula-based prediction strategy for the shape parameter of TransNet}\label{sec:TransNet-auto-tune-gamma}

As stated in \cite{zhang2024transferable}, the tuning process for the shape parameter  reviewed  in Section \ref{ssec:TransNet-GRFs} does not use any information from PDE problems and  is mainly targeted at providing  some reasonable values  for the shape parameter  $\gamma$ in practice and enhancing the transferability of the TransNet across various PDEs with different domains and boundary conditions. Additionally, the use of GRFs for generating auxiliary functions also requires the selection of a new parameter, the correlation length of GRF, an appropriate choice of which is usually related to the variations of the PDE solutions to be solved. 
Furthermore, the optimal value of $\gamma$ for a TransNet also varies along with the number of hidden-layer neurons $M$ even for the same PDE problem. Hence,  such tuning approach could fail to work and  inevitably sacrifice much accuracy of the TransNet method in some situations.

For a specific PDE problem \eqref{equ:pde},  assume that the ball $B_{R}(\boldsymbol{x}_c)$ has been determined for the domain $\Omega$ and the number of hidden-layer neurons $M$ is given.  The location parameters $\left\{(\boldsymbol{a}_m, r_m)\right\}_{m=1}^M$ then can be easily generated according to the process described in \autoref{thm:general-uniform-neuron-distribution}. 
Notice that the physics-informed loss function \eqref{equ:loss} is  a natural error indicator for accuracy of the TransNet solution \eqref{equ:uNN-translate}, so we can define the following posterior error indicator function for TransNet with respect to a given shape parameter $\gamma$: 
$$\eta(\gamma) = \min_{\boldsymbol{\alpha}}\, {\mathrm{Loss_{TN}}}\,({\boldsymbol{\alpha}}).$$
However, the derivative of $\eta(\gamma)$ is hard to derive, and consequently a  simple but effective way is to  utilize a  \emph{gradient-free} line search technique to find optimal shape parameter, i.e., solve the minimization problem 
\begin{equation}\label{tnopt}
\min_{\gamma}\,\eta(\gamma).
\end{equation}
Here we take the popular \emph{golden-section search} algorithm and the details  are described in \autoref{alg:TransNet-golden-shape}. Note that each evaluation of the objective function $\eta(\gamma)$  costs one TransNet solving. 
We will call this approach \emph{the training loss-based optimization strategy} for searching optimal shape parameter of the TransNet.

It is worth noting that when the number of hidden-layer neurons $M$ is small, the proposed training loss-based  optimization strategy is efficient  due to the low computational cost for a small-scale TransNet. However, when $M$ gets larger and larger, the cost of each TransNet solving increase significantly and this optimization  strategy could become computationally expensive and even not acceptable in practice, thus a more efficient selection strategy is still desired to resolve this issue.

\begin{algorithm}[!ht]
	\caption{The golden-section search algorithm  for the shape parameter of TransNet}
	\label{alg:TransNet-golden-shape}
	\KwIn{
A TransNet with the number of hidden-layer neurons $M$ and the location parameters $\left\{(\boldsymbol{a}_m, r_m)\right\}_{m=1}^M$ generated for the ball $B_{R}(\boldsymbol{x}_c)$, the search interval $[a,b]$, and the  number of iterations $Itr$
	}
	\KwOut{Optimal shape parameter $\gamma^{\text{opt}}$}
	Set $\omega=(\sqrt{5}-1)/2$, \hspace{2mm} $\gamma^{(1)} = b - \omega(b-a)$, \hspace{2mm} $\gamma^{(2)} = a + \omega(b-a)$ \;
    Solve \eqref{equ:pde} using the TransNet with $\gamma^{(1)}$ and $\gamma^{(2)}$ as the shape parameter and obtain the loss function value ${res}^{(1)}$  and ${res}^{(2)}$, respectively \;
	\For{$i = 1:Itr$}{
        \eIf{${res}^{(1)} \leqslant {res}^{(2)}$}{
            $b := \gamma^{(2)}$;\hspace{2mm} $\gamma^{(2)} := \gamma^{(1)}$;\hspace{2mm} 
            $\gamma^{(1)} := b - \omega(b-a)$;\hspace{2mm} ${res}^{(2)} := {res}^{(1)}$ \;
            Solve \eqref{equ:pde} using the TransNet with the shape parameter $\gamma^{(1)}$ to obtain a new loss value ${res}^{(1)}$ \;
        }
        {
            $a := \gamma^{(1)}$;\hspace{2mm} $\gamma^{(1)} := \gamma^{(2)}$;\hspace{2mm}
            $\gamma^{(2)} := a + \omega(b-a)$;\hspace{2mm} ${res}^{(1)} := {res}^{(2)}$ \;
            Solve \eqref{equ:pde}  using the TransNet with the shape parameter $\gamma^{(2)}$ to obtain a new loss value ${res}^{(2)}$ \;
        }
    }
    \eIf{${res}^{(1)} \leqslant {res}^{(2)}$}{
        $\gamma^\text{opt} := \gamma^{(1)}$ 
    }
    {
        $\gamma^\text{opt} := \gamma^{(2)}$ 
    }
\end{algorithm}

We observe from the results of \cite{zhang2024transferable}  that the optimal shape parameter $\gamma^\text{opt}$ is positively related to both the number of hidden-layer neurons $M$  and the variation speed of the solution. On the other hand, \eqref{equ:relationship} and \autoref{thm:general-uniform-neuron-distribution} tell us that the product of $\gamma$ and $R$ determines the value range of $b_m$, so there also exists a certain inverse relationship between $\gamma$ and $R$. Based on the above analysis, we propose an empirical formula \eqref{equ:emp-form} for characterizing the relation between these parameters of TransNet as follows:
\begin{equation}\label{equ:emp-form}
	\gamma \approx C \frac{M^{1/d}}{R},
\end{equation}
where $C$ is called the {\em empirical constant} which depends on the settings of the target PDE problem, and $d$ is the problem dimension.

By combining  the training loss-based optimization strategy   and the empirical formula \eqref{equ:emp-form},  we propose \emph{an empirical formula-based prediction strategy}  for estimating appropriate shape parameter of the TransNet. The strategy consists of a preprocessing step and a prediction step as follows: \vspace{-0.1cm}
\begin{itemize}
\item {\bf Preprocessing}: Select a small  number of hidden-layer neurons $M_0$ and use the training loss-based  optimization strategy (e.g., \autoref{alg:TransNet-golden-shape})  to find an optimal shape parameter $\gamma_0^\text{opt}$ for the current TransNet, and consequently, set the empirical constant $C := \frac{\gamma_0^\text{opt}R}{(M_0)^{1/d}}$
in the  empirical formula \eqref{equ:emp-form}. 
\vspace{-0.2cm}
\item {\bf Prediction}: Given any new (usually larger) value  of $M$ to be used for the TransNet, directly predict an appropriate  shape parameter $\gamma_*$ according to the  empirical formula \eqref{equ:emp-form}, that is, $\gamma_* := C \frac{M^{1/d}}{R}$.
\end{itemize}\vspace{-0.1cm}
Note that the prediction step  is completely explicit (does not involve any training or optimization) and the computational cost  for the preprocessing step is low since $M_0$ is small, and thus this strategy is very efficient in practice.

\section{A multiple transferable neural network method for elliptic interface problems}\label{sec:Multi-TransNet}
\begin{figure}[!ht]
    \centering
    \includegraphics[width=0.88\linewidth]{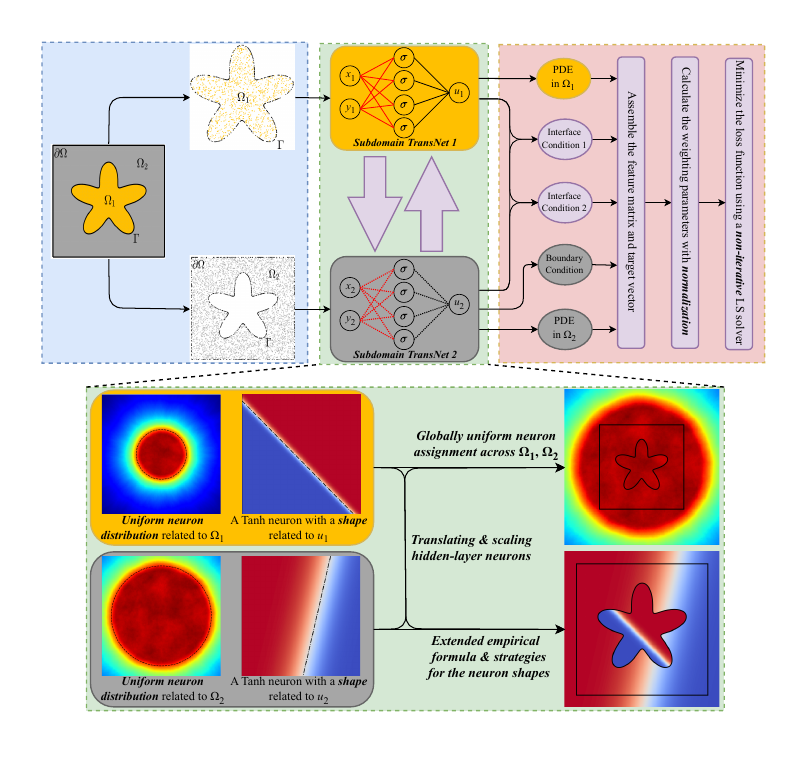}
    \caption{\fontsize{8.5pt}{8.5pt}\selectfont{Illustration of the solution process of the proposed Multi-TransNet method for the elliptic interface problem \eqref{equ:IP} with $K=2$.}}
    \label{fig:Multi-TransNet}
\end{figure}

It is well known that the activation function  plays a pivotal role in the expressive power of the neural network and it is generally continuous and even differentiable. However, the solution to interface problems often exhibits nonsmoothness or even discontinuities, which incurs poor performance of neural network methods. Specifically, for the elliptic interface problem \eqref{equ:IP} with $K$ subdomains, the solution in the whole domain is divided into $K$ parts located in different subdomains as illustrated in \autoref{fig:IP-domain}, in each of which the solution is usually smooth, but has cusps or jumps at the interfaces. In order to conquer such a problem, a natural idea is to respectively employ a neural network to approximate the smooth solution over each subdomain and then to unite them through the interface conditions on interfaces as done by many existing numerical methods. Therefore, we integrate multiple \emph{tailored} TransNets related to subdomains using the \emph{nonoverlapping domain decomposition} approach to develop a novel multiple transferable neural network method, abbreviated as Multi-TransNet, for solving the elliptic interface problem \eqref{equ:IP}. \autoref{fig:Multi-TransNet} illustrates the solution process of  the proposed  Multi-TransNet method  when $K=2$.

\subsection{Nonoverlapping domain decomposition and integrated subdomain TransNets}\label{ssec:NDD}
Without loss of generality, here we take the common case of $K=2$ subdomains ($\Omega_1$ and $\Omega_2$) and one interface $\Gamma$ (see Figure \ref{fig:IP-domain}-left)  for an illustration  of the proposed  Multi-TransNet method for solving the elliptic interface problem \eqref{equ:IP}. To be detailed, we employ two TransNets $u^\mathrm{NN}_1$ and $u^\mathrm{NN}_2$ with respective $M_1$ and $M_2$ hidden-layer neurons to approximate the solution in $\Omega_1$ and $\Omega_2$, respectively, and then the Multi-TransNet solution with totally $M=M_1+M_2$ hidden-layer neurons can be represented by
\begin{equation}\label{equ:Multi-TransNet-solution}
	u^\mathrm{NN} = u^\mathrm{NN}_1\,\chi_{\Omega_1}+u^\mathrm{NN}_2\,\chi_{\Omega_2},
	\end{equation}
where $\left\{\chi_{\Omega_k}\right\}_{k=1}^2$ are the indicator functions and 
\begin{equation}
	u^\mathrm{NN}_k = \boldsymbol{\alpha}_k^\top \boldsymbol{\phi}_k, \quad k=1,2,
\end{equation}
with $\boldsymbol{\alpha}_k=\left( \alpha_0^{(k)}, \alpha_1^{(k)},\dots, \alpha_{M_k}^{(k)}\right)^\top,\  \boldsymbol{\phi}_k = \left( \phi_0^{(k)}, \phi_1^{(k)},\dots, \phi_{M_k}^{(k)}\right)^\top$, and  the neural basis functions 
\begin{equation*}
	\phi_0^{(k)} = 1,\quad \phi_m^{(k)} (\boldsymbol{x}) = \sigma\left( \gamma_k\left(\left(\boldsymbol{x}-\boldsymbol{x}^{(k)}_{c}\right)^\top \boldsymbol{a}_m^{(k)} + r_m^{(k)} \right)\right),\quad m=1,2,\dots,M_k.
\end{equation*}
The corresponding loss function of our Multi-TransNet for solving \eqref{equ:IP} is then designed to be
\begin{equation}\label{equ:lossmul}
\begin{aligned}
{\mathrm{Loss_{MT}}}\,(\boldsymbol{\alpha}_1,\boldsymbol{\alpha}_2) ={}& 
\left\|\lambda_1 \left(\mathcal{L}\left(u^{\mathrm{NN}}\right) - f\right)\right\|_{\Omega_1,2}^2 + 
\left\|\lambda_2 \left(\mathcal{L}\left(u^{\mathrm{NN}}\right) - f\right)\right\|_{\Omega_2,2}^2 + 
\left\|\lambda_g \left(\mathcal{B}\left(u^{\mathrm{NN}}\right) - g\right)\right\|_{\partial\Omega,2}^2 \\
{}& + 
\left\|\lambda_{h,1} \left(\left[u^\mathrm{NN}\right]-h_1\right)\right\|_{\Gamma,2}^2 + 
\left\|\lambda_{h,2} \left(\left[\mathcal{J}\left(u^\mathrm{NN}\right) \cdot \boldsymbol{n}\right]-h_2\right)\right\|_{\Gamma,2}^2,
\end{aligned}
\end{equation}
where  $\lambda_{1}$, $\lambda_{2}$, $\lambda_{g}$, $\lambda_{h,1}$, $\lambda_{h,2}$ are some positive weighting parameters.

In practical implementation, assume that $N_1$, $N_2$, $N_g$ and $N_\Gamma$ training/collocation points have been sampled in $\Omega_1$, $\Omega_2$, $\partial \Omega$ and $\Gamma$, respectively, which are arranged as follows:
\begin{equation}
    \boldsymbol{X}_1=
    \begin{pmatrix}
        \left(\boldsymbol{x}^{(1)}_{1}\right)^\top\\
        \left(\boldsymbol{x}^{(1)}_{2}\right)^\top\\
        \vdots \\
        \left(\boldsymbol{x}^{(1)}_{N_1}\right)^\top
    \end{pmatrix}, \quad
    \boldsymbol{X}_2=
    \begin{pmatrix}
        \left(\boldsymbol{x}^{(2)}_{1}\right)^\top\\
        \left(\boldsymbol{x}^{(2)}_{2}\right)^\top\\
        \vdots \\
        \left(\boldsymbol{x}^{(2)}_{N_2}\right)^\top
    \end{pmatrix}, \quad
    \boldsymbol{X}_g=
    \begin{pmatrix}
        \boldsymbol{x}_{g,1}^\top\\
        \boldsymbol{x}_{g,2}^\top\\
        \vdots \\
        \boldsymbol{x}_{g,N_g}^\top
    \end{pmatrix}, \quad
    \boldsymbol{X}_{\Gamma}=
    \begin{pmatrix}
        \boldsymbol{x}_{\Gamma,1}^\top\\
        \boldsymbol{x}_{\Gamma,2}^\top\\
        \vdots \\
        \boldsymbol{x}_{\Gamma,N_\Gamma}^\top
    \end{pmatrix}.
\end{equation}
Next, we substitute the Multi-TransNet solution \eqref{equ:Multi-TransNet-solution} into the elliptic interface problem \eqref{equ:IP} and evaluate them at the above training/collocation points, and then minimize the squared residual (i.e., the loss function \eqref{equ:lossmul} in the discrete sense), i.e., 
\begin{equation}\label{equ:Multi-TransNet-ls}
	\min_{\boldsymbol{\alpha}} \|\boldsymbol{F}\boldsymbol{\alpha}- \boldsymbol{T}\|_2^2,
\end{equation}
where the feature matrix $\boldsymbol{F}$, the target vector $\boldsymbol{T}$ and the parameters $\boldsymbol{\alpha}$ of the output layer are respectively assembled in the following manner:
\begin{equation}\label{equ:assemble}
	\boldsymbol{F} = \left(\begin{array}{c:c}
		\lambda_1 \mathcal{L}  \left(\boldsymbol{\phi}^\top_1\left(\boldsymbol{X}_1\right)\right) & \boldsymbol{O} \\
		\boldsymbol{O} & \lambda_2 \mathcal{L} \left(\boldsymbol{\phi}^\top_2\left(\boldsymbol{X}_2\right)\right) \\
		\boldsymbol{O} & \lambda_{g} \mathcal{B} \left( \boldsymbol{\phi}^\top_2\left(\boldsymbol{X}_g \right)\right) \\
		\lambda_{h,1} \boldsymbol{\phi}^\top_1\left(\boldsymbol{X}_\Gamma\right) & -\lambda_{h,1}\boldsymbol{\phi}^\top_2\left(\boldsymbol{X}_\Gamma\right) \\
		\lambda_{h,2} \mathcal{J}\left(\boldsymbol{\phi}^\top_1\left(\boldsymbol{X}_\Gamma\right)\right)\cdot \boldsymbol{n} & - 
        \lambda_{h,2} \mathcal{J}\left(\boldsymbol{\phi}^\top_2\left(\boldsymbol{X}_\Gamma\right)\right)\cdot \boldsymbol{n}
	\end{array}\right), \quad
	\boldsymbol{T} = \begin{pmatrix}
		\lambda_1 f \left(\boldsymbol{X}_1\right) \\
		\lambda_2 f \left(\boldsymbol{X}_2\right) \\
		\lambda_{g} g\left( \boldsymbol{X}_g \right) \\
		\lambda_{h,1} h_1\left(\boldsymbol{X}_\Gamma\right) \\
		\lambda_{h,2} h_2\left(\boldsymbol{X}_\Gamma\right) 
	\end{pmatrix}, \quad 
	\boldsymbol{\alpha}=
	\begin{pmatrix}
		\boldsymbol{\alpha}_1 \\
		\boldsymbol{\alpha}_2 
	\end{pmatrix}.
\end{equation}
It is clearly $\boldsymbol{F}\in \mathbb{R}^{\left(N_1+N_2+N_g+2N_\Gamma\right)\times \left(M+2\right)}$, $\boldsymbol{T}\in \mathbb{R}^{N_1+N_2+N_g+2N_\Gamma}$ and $\boldsymbol{\alpha}\in \mathbb{R}^{M+2}$, and the minimization problem \eqref{equ:Multi-TransNet-ls}
is in fact  a linear least squares problem, thus $\boldsymbol{F}$ is also called  the least squares coefficient matrix.

\begin{remark}
To briefly demonstrate the advantage of the proposed Multi-TransNet method over the TransNet method, we consider a 1D elliptic interface problem as follows:
\begin{equation}\label{equ:1D-IP}
	\left\{\begin{aligned}
		-\left(\beta u^\prime\right)^\prime & = f,\quad\;\, x \in \Omega_1\cup \Omega_2=(0,1/4)\cup (1/4, 1), \\
		[u] & =h_1,\quad x \in \Gamma=\{1/4\}, \\
		\left[\beta u^\prime\right] & =h_2,\quad x \in \Gamma=\{1/4\}, \\
		u & =g,\quad\;\, x \in \partial \Omega=\{0, 1\},
	\end{aligned}\right.
\end{equation}
where $\beta=1$ for $x\in (0,1/4)$ and $\beta=10$ for $x\in (1/4, 1)$.
The exact solution is set to be
\begin{equation}\label{equ:1D-IP-ExaSol}
	u(x) = \begin{cases}
		\sin (2\pi x), &x \in [0,1/4], \\
		\cos (2\pi x), &x \in (1/4,1],
	\end{cases}
\end{equation}
and $f, h_1, h_2, g$ in \eqref{equ:1D-IP} can be obtained by substituting \eqref{equ:1D-IP-ExaSol} into \eqref{equ:1D-IP}. 
\autoref{fig:NDD} illustrates the comparison of the numerical solutions and gradients obtained by using one TransNet with 100 hidden-layer neurons and the proposed Multi-TransNet with 5 hidden-layer neurons for each of the two subdomain TransNets for the 1D elliptic interface problem \eqref{equ:1D-IP}. It is easily observed that the TransNet method fails while the Multi-TransNet method catches the discontinuous solution and gradient accurately.

\begin{figure}[!ht]
    \centering
    \begin{minipage}[t]{.48\textwidth}
		\centering
		\includegraphics[width=\linewidth]{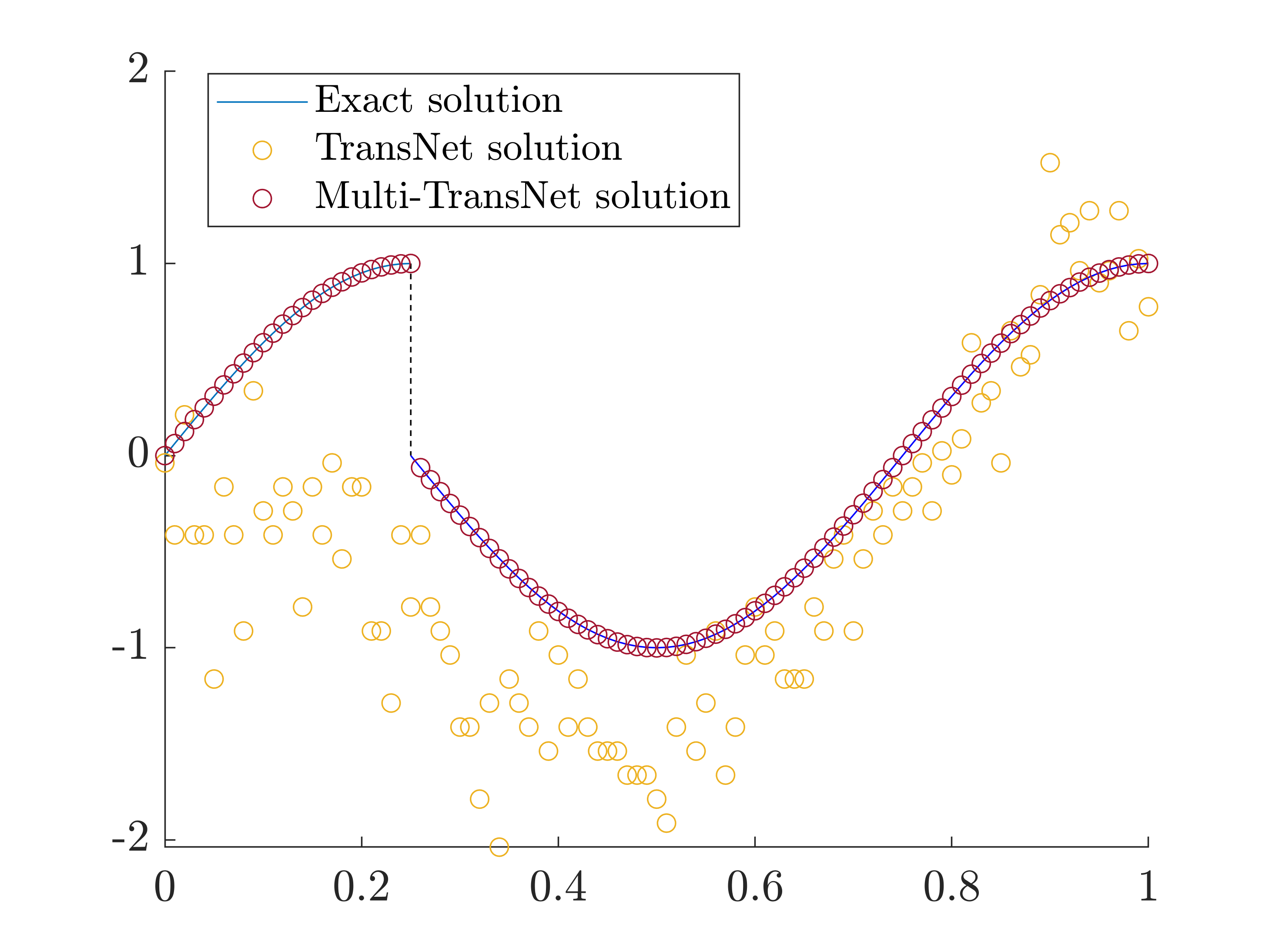}	
	\end{minipage}
	\begin{minipage}[t]{.48\textwidth}
		\centering
		\includegraphics[width=\linewidth]{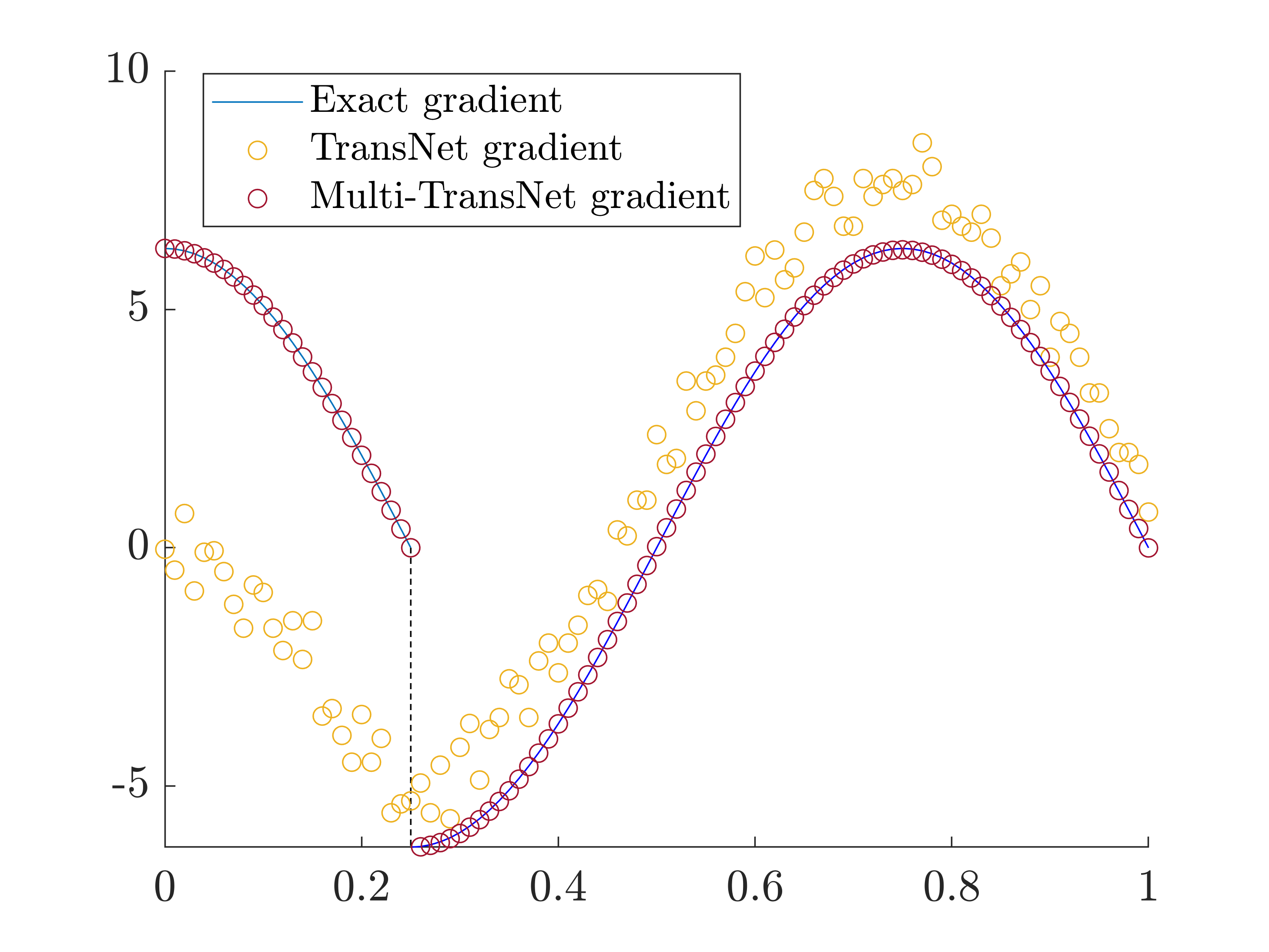}	
	\end{minipage}
    \vspace{-2mm}
    \caption{Comparisons of the numerical solutions (left) and gradients (right) produced by one TransNet with 100 hidden-layer neurons and the proposed Multi-TransNet with 5 hidden-layer neurons for each of the two subdomain TransNets for the 1D elliptic interface problem \eqref{equ:1D-IP}.}
    \label{fig:NDD}
\end{figure}
\end{remark}

The extension of the Multi-TransNet method to the case of $K>2$ subdomains $\{\Omega_k\}_{k=1}^K$ with multiple interfaces $\{\Gamma_{i,j}\}$  is straightforward. We assign each subdomain $\Omega_k$ a TransNet $u_k^{\mathrm{NN}}$ with $M_k$ hidden-layer neurons to form a total of $K$ subdomain TransNets $\{u_k^{\mathrm{NN}}\}_{k=1}^{K}$ such that the Multi-TransNet solution is given by
\begin{equation}\label{equ:Multi-TransNet-sol-general}
\begin{aligned}
u^\mathrm{NN} 
&\;=\sum_{k=1}^K u^\mathrm{NN}_{k}~\chi_{\Omega_k}
=\sum_{k=1}^K \boldsymbol{\alpha}_k^\top \boldsymbol{\phi}_k~\chi_{\Omega_k}\\
&\;=\sum_{k=1}^K \left(\alpha_0^{(k)}+\sum_{m=1}^{M_k}\alpha_m^{(k)}\sigma\left(\gamma_k\left(\left({\bm x}-{\bm x}^{(k)}_{c}\right)^\top {\bm a}_m^{(k)} + r_m^{(k)}\right)\right)\right)~\chi_{\Omega_k}.
\end{aligned}
\end{equation}
The way of constructing the loss function, assembling the feature matrix, the target vector and the parameters of the output layer is similar to  \eqref{equ:lossmul}-\eqref{equ:assemble}.

\subsection{Globally uniform neuron distribution across subdomains} \label{ssub-dglob}
Rather than arbitrarily assigning or equally distributing the number of hidden-layer neurons among all subdomains, we  consider the globally uniform neuron distribution of Multi-TransNet for the elliptic interface problem \eqref{equ:IP} with $K$ subdomains. 
At first, from the view of the neuron location, the hidden-layer neurons for each of the $K$ subdomain  TransNets should be respectively translated to the approximate center $\boldsymbol{x}^{(k)}_c$ of the corresponding subdomain $\Omega_k$, and also be respectively scaled to enable the ball $B_{R_k}(\boldsymbol{x}^{(k)}_c)$ to slightly over-cover the subdomain $\Omega_k$, as emphasized in \autoref{rem:translatation-and-scaling}. 
Since the uniform distribution of hidden-layer neurons can bring better transferability and accuracy for TransNet for solving general PDE problems as demonstrated in \cite{zhang2024transferable}, a natural question is how to extend the uniform distribution feature of hidden-layer neurons on each subdomain to the \emph{globally} uniform distribution across the whole domain for the Multi-TransNet.

\begin{figure}[!t]
	\centering
	\includegraphics[width=0.95\textwidth]{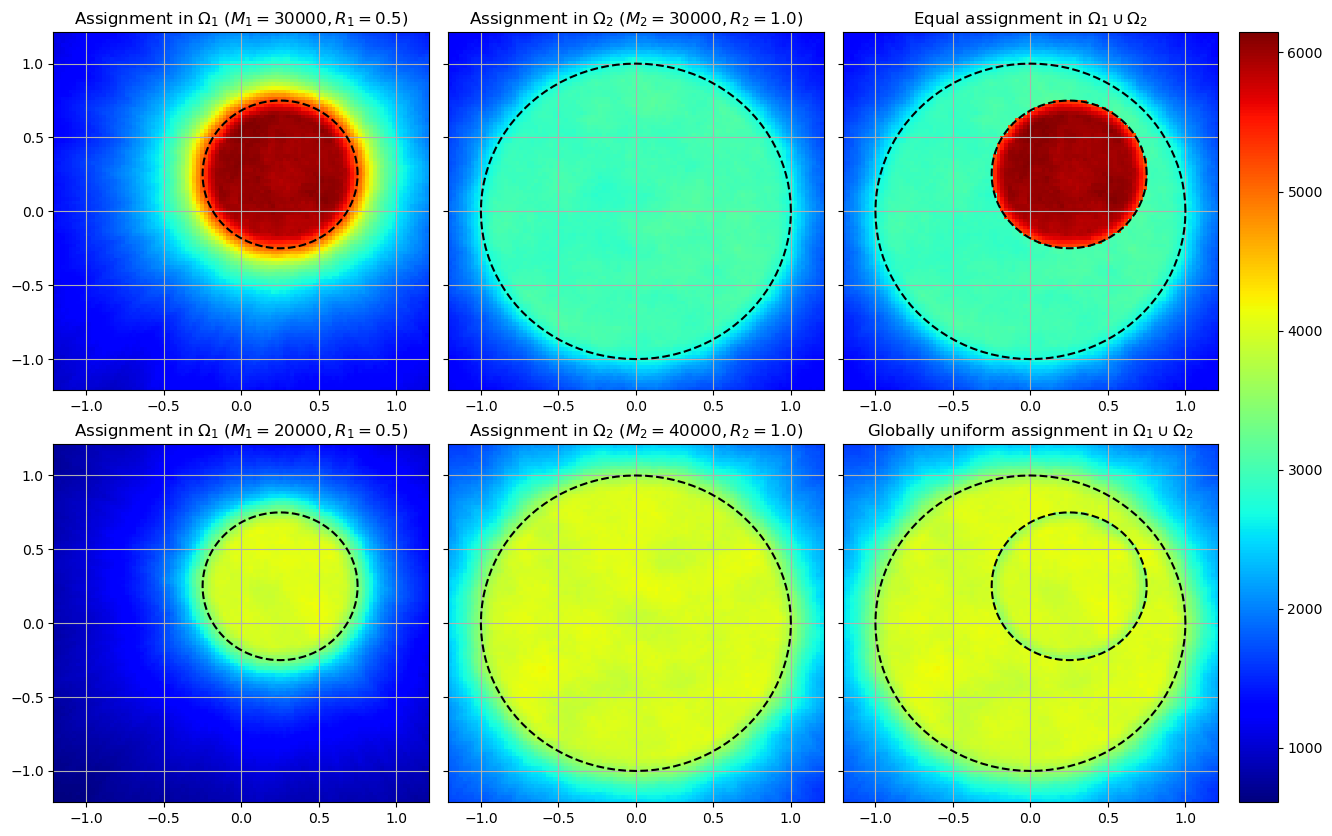}
	\vspace{-2mm}
	\caption{Comparison of the number distribution of hidden-layer neurons of the Multi-TransNet with two subdomain TransNets with the subdomains $\Omega_1 = B_{0.5}\left(0.25, 0.25\right)$ and $\Omega_2 = B_{1}\left(0, 0\right)\setminus B_{0.5}\left(0.25, 0.25\right)$. Here $M_1+M_2$ and $\uptau$ are fixed to 60000 and 0.1, respectively. Top row: $M_1=M_2=30000$; Bottom row: $M_1=20000$ and $M_2=40000$ thus ${M_1}/{R_1} = {M_2}/{R_2}$.}
	\label{fig:general-uniform-neuron-number-distribution-TransNets}
\end{figure}

Let us rewrite the main result \eqref{res2} in \autoref{thm:general-uniform-neuron-distribution} into the following form:
\begin{equation}
	\mathbb{E}\left[\sum_{m=1}^M \chi_{\left\{d_m(\boldsymbol{x})<\uptau\right\}}(\boldsymbol{x})\right]=\frac{M}{R}\uptau,\quad \forall\, \|\boldsymbol{x} - \boldsymbol{x}_c\|_2 \leqslant R-\uptau,
\end{equation}
which implies that the expected number of hidden-layer neurons whose corresponding partition hyperplanes intersect the neighborhood $B_{\uptau}(\boldsymbol{x})$ is proportional to ${M}/{R}$. It is natural  to assume that the same neighborhood size $\uptau$ is used for defining the acting range of hidden-layer neurons in each subdomain TransNet of the Multi-TransNet. As a result, in order to obtain the globally uniform distribution of hidden-layer neurons, we need to make the ratio ${M_k}/{R_k}$ as equal as possible (note $M_k$ is always an integer), i.e., 
\begin{equation}\label{equ:equal-ratio-M-R}
    \frac{M_1}{R_1} \approx \frac{M_2}{R_2} \approx \cdots \approx \frac{M_K}{R_K}.
\end{equation}
\autoref{fig:general-uniform-neuron-number-distribution-TransNets} illustrates comparison of the number distribution of hidden-layer neurons of the Multi-TransNet with two subdomain TransNets  under the setting of the equal assignment (i.e., $M_1=M_2$) and ${M_1}/{R_1} ={M_2}/{R_2}$, where $\Omega_1 = B_{0.5}\left(0.25, 0.25\right)$ and $\Omega_2 = B_{1}\left(0, 0\right)\setminus B_{0.5}\left(0.25, 0.25\right)$.
It is easy to see that the setting of ${M_1}/{R_1} = {M_2}/{R_2}$ indeed ensures the globally uniform distribution of hidden-layer neurons while the setting of the equal assignment only achieves the locally uniform distribution of those and enables more hidden-layer neurons to gather in the interior circular subdomain $\Omega_1$. Hence, for achieving globally uniform neuron distribution, we propose to use  the relation  \eqref{equ:equal-ratio-M-R} to assign the number of neurons for each subdomain TransNet,  which is equivalent to 
\begin{equation}\label{equ:equal-ratio-M-R-des}
    M_k \approx M\frac{R_k}{\sum_{k=1}^K R_k},
\end{equation}
where $M=\sum_{k=1}^K M_k$ is the total number of hidden-layer neurons of the Multi-TransNet. 

\subsection{Extended empirical formula and strategies for the hidden-layer neuron shapes of Multi-TransNet} \label{ssub-eef}
For simplicity of illustration, we first pay attention to the shape parameters of a Multi-TransNet with two subdomain TransNets. As introduced in \autoref{ssec:TransNet-GRFs}, the shape parameter determines the variation speed of the neural basis function  and affects the approximation performance of TransNet, so does Multi-TransNet. The solutions in the subdomains $\Omega_1$ and $\Omega_2$ could be significantly distinct, thus we need to provide  respective shape parameter for each of the two subdomain  TransNets (i.e., $\gamma_1$ and $\gamma_2$) and a corresponding efficient tuning strategy for optimal ones is again needed. Similar to the TransNet case,
we define a posterior error indicator function for the Multi-TransNet with respect to given shape parameters $(\gamma_1, \gamma_2)$ based on the loss function \eqref{equ:lossmul} as follows:
$$\eta(\gamma_1,\gamma_2) = \min_{(\boldsymbol{\alpha}_1,\boldsymbol{\alpha}_2)}\, {\mathrm{Loss_{MT}}}\,({\boldsymbol{\alpha}}).$$
Furthermore, we propose to \emph{bridge} the shape parameters of the two subdomain TransNets by extending the empirical formula \eqref{equ:emp-form} introduced for the TransNet in \autoref{sec:TransNet-auto-tune-gamma} into the Multi-TransNet case. Specifically, we assume that the two subdomain TransNets with the respective parameters, $(M_1, R_1)$ and $(M_2, R_2)$, share the same empirical constant $C$ in \eqref{equ:emp-form}, i.e.,
\begin{equation}\label{111}
		\gamma_1 \approx C \frac{M_1^{1/d}}{R_1}, \quad
		\gamma_2 \approx C \frac{M_2^{1/d}}{R_2},
\end{equation}
then consequently the following relation holds
\begin{equation}\label{equ:gamma-rela}
	{\gamma_2}\approx {\gamma_1} \frac{R_1}{R_2} \left(\frac{M_2}{M_1}\right)^{1/d}.
\end{equation}

Therefore, the  bivariate optimization problem for searching optimal shape parameters $(\gamma_1, \gamma_2)$ of the Multi-TransNet
can be converted into a univariate one with respect to $\gamma_1$ (or $\gamma_2$) only, i.e.,
\begin{equation}\label{mtopt}
\min_{\gamma_1}\,  \eta\left(\gamma_1,\gamma_1\textstyle\frac{R_1}{R_2} \left(\frac{M_2}{M_1}\right)^{1/d}\right).
\end{equation}
 Then the training loss-based optimization strategy 
also can be slightly modified and used to find optimal $\gamma_1$ for the problem \eqref{mtopt}, which  then can be applied to  \eqref{111} to compute the needed empirical constant $C$. Consequently, the efficient empirical formula-based  prediction strategy 
presented in \autoref{sec:TransNet-auto-tune-gamma} also  can be similarly generalized to the case of Multi-TransNet for  selecting  appropriate values of the shape parameters $\gamma_1$ and  $\gamma_2$.

In the case of $K$ subdomain TransNets with the parameters $\left\{M_k\right\}_{k=1}^K$ and $\left\{R_k\right\}_{k=1}^K$, for the set of shape parameters $\left\{\gamma_k\right\}_{k=1}^K$ of the corresponding Multi-TransNet, we have under the same principle:
\begin{equation}\label{222}
		\gamma_k \approx C \frac{M_k^{1/d}}{R_k},  \quad k=1,2,\cdots, K,
\end{equation}
and thus 
\begin{equation}
    {\gamma_k}\approx \gamma_1\frac{R_1}{R_k} \left(\frac{M_k}{M_1}\right)^{1/d}, \quad k=2,3,\dots, K.
\end{equation}
Thus the two selection strategies  discussed above  can be straightforwardly extended.

\begin{remark}
In particular, when $M_1\approx M_2 \approx  \cdots \approx M_k$,  we have
\begin{equation}
	\gamma_k \approx \gamma_1\frac{R_1}{R_k}, \quad k=2, 3, \cdots, K,
\end{equation}
and if the globally uniform distribution pattern  is adopted, i.e.,  ${M_1}/{R_1} \approx {M_2}/{R_2}\approx\cdots \approx {M_k}/{R_k}$, then it holds
\begin{equation}
\gamma_k \approx \gamma_1\left(\frac{R_1}{R_k}\right)^{1-\frac{1}{d}}, \quad k=2, 3, \cdots, K.
\end{equation}
\end{remark}

\subsection{Adaptively determining the loss weighting parameters through normalization} \label{ssub-hyp}
The choice of the  weighting parameters in the loss function of neural network methods may significantly affect the performance of neural network methods.
Inspired by the work of \cite{chi2024random}, we propose to determine the loss weighting  parameters based on the magnitudes of 
elements of the augmented matrix combining the least squares coefficient matrix with the target vector for the proposed Multi-TransNet method. To be specific, let us again  take the case of $K=2$ for explanation, and the weighting  parameters in the loss function \eqref{equ:assemble} will be calculated through the following  \emph{normalization}:
\begin{equation}\label{equ:penalty}
\begin{aligned}
    \lambda_k &= \frac{1}{\max\limits_{\substack{1\leqslant i \leqslant N_k \\ 
    1\leqslant j \leqslant M_k+2}}
    \left|\left(\hspace{-2mm}
    \begin{array}{c:c}
         \mathcal{L} \left(\boldsymbol{\phi}_k^\top \left(\boldsymbol{X}_k \right) \right)\hspace{1mm} & \hspace{1mm} f\left(\boldsymbol{X}_k\right)\hspace{-2mm} \\
    \end{array} \right)_{i,j}\right|}, \quad k=1,2,\\
    \lambda_g &= \frac{1}{\max\limits_{\substack{1\leqslant i \leqslant N_g \\ 
    1\leqslant j \leqslant M_2+2}}
    \left|\left(\hspace{-2mm} 
    \begin{array}{c:c}
         \mathcal{B} \left(\boldsymbol{\phi}_2^\top \left(\boldsymbol{X}_g \right) \right)\hspace{1mm} & \hspace{1mm} g\left(\boldsymbol{X}_g\right)\hspace{-2mm} \\
    \end{array} \right)_{i,j}\right|},\\
    \lambda_{h,1} &= \frac{1}{\max\limits_{\substack{1\leqslant i \leqslant N_\Gamma \\ 
    1\leqslant j \leqslant M+3}}
    \left|\left(\hspace{-2mm} 
    \begin{array}{c:c:c}
         \boldsymbol{\phi}_1^\top \left(\boldsymbol{X}_\Gamma \right) \hspace{1mm} & \hspace{1mm}
         -\boldsymbol{\phi}_2^\top \left(\boldsymbol{X}_\Gamma \right) \hspace{1mm} & \hspace{1mm}
         h_1\left(\boldsymbol{X}_\Gamma\right)\hspace{-2mm} \\
    \end{array} \right)_{i,j}\right|},\\
    \lambda_{h,2} &= \frac{1}{\max\limits_{\substack{1\leqslant i \leqslant N_\Gamma \\ 
    1\leqslant j \leqslant M+3}}
    \left|\left(\hspace{-2mm}
    \begin{array}{c:c:c}
         \mathcal{J} \left(\boldsymbol{\phi}_1^\top \left(\boldsymbol{X}_\Gamma \right) \right) \cdot \boldsymbol{n} \hspace{1mm} & \hspace{1mm}
         -\mathcal{J} \left(\boldsymbol{\phi}_2^\top \left(\boldsymbol{X}_\Gamma \right) \right) \cdot \boldsymbol{n} \hspace{1mm} & \hspace{1mm}
         h_2\left(\boldsymbol{X}_\Gamma\right)\hspace{-2mm} \\
    \end{array} \right)_{i,j}\right|}, 
\end{aligned}
\end{equation}
where the dashed vertical lines denote the matrix augmentation operation. The extension of such normalization approach to the Multi-TransNet with  multiple interfaces and $K$ subdomains is again straightforward.
We also note the the approach developed in \cite{chi2024random} applies the normalization to the  least squares coefficient  matrix only
for obtaining the loss weighting parameters, i.e., removing the parts involving $f$, $g$, $h_1$ and h$_2$ in the case of \eqref{equ:penalty}.
  
Finally, we summarize in \autoref{alg:Multi-TransNet} the main steps of the proposed Multi-TransNet method for solving the elliptic interface problem \eqref{equ:IP}, i.e., compute its Multi-TransNet solution \eqref{equ:Multi-TransNet-sol-general}.

\begin{algorithm}[!ht]
	\caption{The Multi-TransNet method for solving the elliptic interface problem \eqref{equ:IP}}
	\label{alg:Multi-TransNet}
	 \KwIn{
     The total number of hidden-layer neurons $M$.
	}
	\KwOut{The output-layer parameters $\{\boldsymbol{\alpha}_k\}_{k=1}^K$ of Multi-TransNet
    }
	
    Select appropriately the parameters $\{(\boldsymbol{x}^{(k)}_{c},R_k)\}_{k=1}^K$ for all subdomain TransNets according to the location and radius magnitude of the $K$ subdomains $\{\Omega_k\}_{k}^{K}$.
    
    Calculate the number of hidden-layer neurons for all subdomain  TransNets,  $\{M_k\}_{k=1}^K$, by using  \eqref{equ:equal-ratio-M-R-des} with $\sum_{k=1}^K M_k = M$.
	
	Generate the location parameters of hidden-layer neurons  for all subdomain  TransNets, $\left\{(\boldsymbol{a}_m^{(1)}, r_m^{(1)})\right\}_{m=1}^{M_1}$,  $\dots$, $\left\{(\boldsymbol{a}_m^{(K)}, r_m^{(K)})\right\}_{m=1}^{M_K}$
	according to the sampling process described in \autoref{thm:general-uniform-neuron-distribution}.
	
	Determine the shape parameters of hidden-layer neurons for all subdomain  TransNets, $\{\gamma_k\}_{k=1}^K$, by  using the empirical formula-based prediction strategy.
	
	Construct the neural basis functions for all subdomain  TransNets, $\{\boldsymbol{\phi}_k\}_{k=1}^{K}$.

	Sample the training/collocation points (randomly or uniformly) on subdomains, domain boundaries and interfaces.		
	
	Assemble the feature matrix and target vector (as done in \eqref{equ:assemble}) and adaptively determine the weighting parameters for all terms in the loss function with the normalization technique (as done in \eqref{equ:lossmul} and \eqref{equ:penalty}).
	    
      Find the output-layer parameters $\{\boldsymbol{\alpha}_k\}_{k=1}^K$ by minimizing the loss function \eqref{equ:Multi-TransNet-ls} (i.e., solving a linear least squares problem).
\end{algorithm}

\section{Numerical experiments}\label{sec:num-experim}
To demonstrate the superior accuracy, efficiency and robustness of the proposed Multi-TransNet method,  we perform abundant numerical tests in this section, including ablation studies and comparison with recent neural network methods and traditional numerical techniques for solving elliptic interface problems. 

The hyperbolic tangent $\mathit{tanh}$ is adopted as the activation as in \cite{zhang2024transferable} owing to its good smoothness. 
The training/collocation points are sampled uniformly in the computational domain (including interior domains, boundaries and/or interfaces), and the test points are generated via the Latin hypercube sampling method \cite{tang1993orthogonal}, which ensures that the test points are more evenly distributed across the range of each dimension. The number of training/collocation points is determined by the spacing size due to uniform sampling adopted (though the proposed Multi-TransNet method is mesh-free), and the number of test points is set to $2^d$ times as large as that of training/collocation points, where $d$ is the problem dimension.
When the golden-section search algorithm (\autoref{alg:TransNet-golden-shape}) is used,  we set the initial search interval  to be $[0, 5]$ and the number of search iteration to be 7.
For measurement of the accuracy of a numerical solution $\boldsymbol{u}$, we use the discrete relative L$_2$ error (referred to as RL$_2$ for short) and the discrete relative L$_\infty$ error (referred to as RL$_\infty$ for short) over all the test points, which are defined as
\[
\text{RL}_2 (\boldsymbol{u}) = \frac{\|\boldsymbol{u}-\boldsymbol{u}_\text{exact}\|_2}{\|\boldsymbol{u}_\text{exact}\|_2}, \quad
\text{RL}_\infty (\boldsymbol{u}) = \frac{\|\boldsymbol{u}-\boldsymbol{u}_\text{exact}\|_\infty}{\|\boldsymbol{u}_\text{exact}\|_\infty}.
\]

The API $\mathit{torch.linalg.lstsq}$ from PyTorch is called for solving the least squares problems. All the experiments are implemented on an Ubuntu 20.04.6 LTS server with a 3.00-GHz Intel Xeon Gold 6248R CPU and a NVIDIA GeForce RTX 4090 GPU. All experimental results are obtained by repeating 10 runs, removing two extrema and then taking the averages of the remaining to diminish the influence of randomness. 
The right-hand terms of PDEs, boundary conditions and/or interface conditions in each example can be derived from the given exact solution, and hence are not listed separately.

\subsection{Ablation studies}\label{sec:abla-studies}
In this subsection, we perform systematic ablation studies to verify  effectiveness and benefits of some important components of the TransNet and the proposed Multi-TransNet methods.

\subsubsection{TransNet}\label{ssec:abla-TransNet}
In this subsection, we test the TransNet method
by solving the classic 2D Poisson problem defined in the domain $\Omega=[0,2]^2$:
\begin{equation}\label{equ:TransNet-abla-ex}
	\left\{\begin{aligned}
		\Delta u & = f,\quad (x,y) \in \Omega, \\
		u & =g,\quad (x,y) \in \partial \Omega.
	\end{aligned}\right.
\end{equation}
The exact solution is taken as $u(x,y) = \sin x \sin y$. We uniformly sample $49^2$ points in the interior and $51\times4$ points on the boundary, with an overall spacing of approximately $0.04$, for training.

\begin{figure}[!ht]
	\centering
	\includegraphics[width=.88\textwidth]{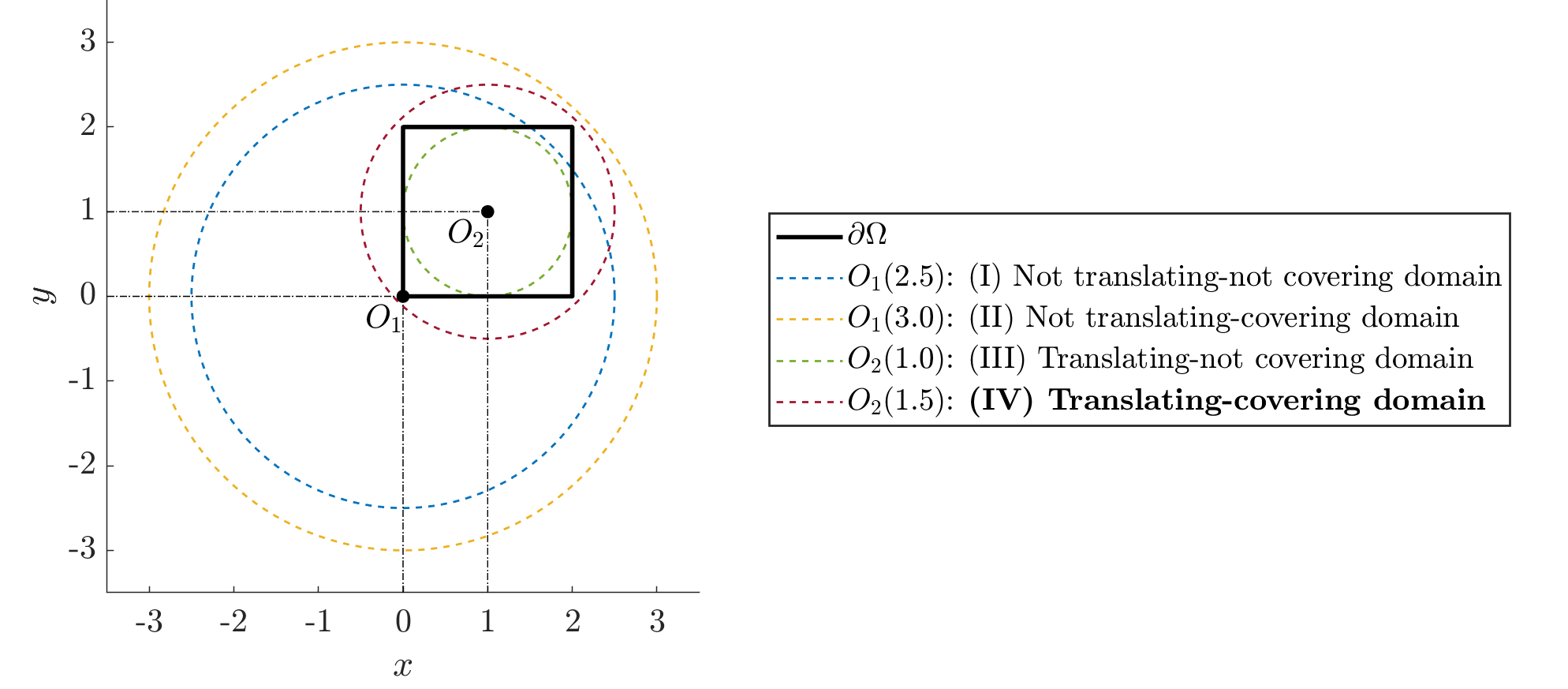}
	\vspace{-2mm}
	\caption{Illustration of different parameter settings for the problem \eqref{equ:TransNet-abla-ex} in \autoref{tab:TransNet-abla-2D-empi-const}.}
	\label{fig:TransNet-abla-2D-trans-scal}
\end{figure}

\begin{figure}[!ht]
	\centering
	\begin{minipage}[t]{.48\textwidth}
		\centering
		\includegraphics[width=\linewidth]{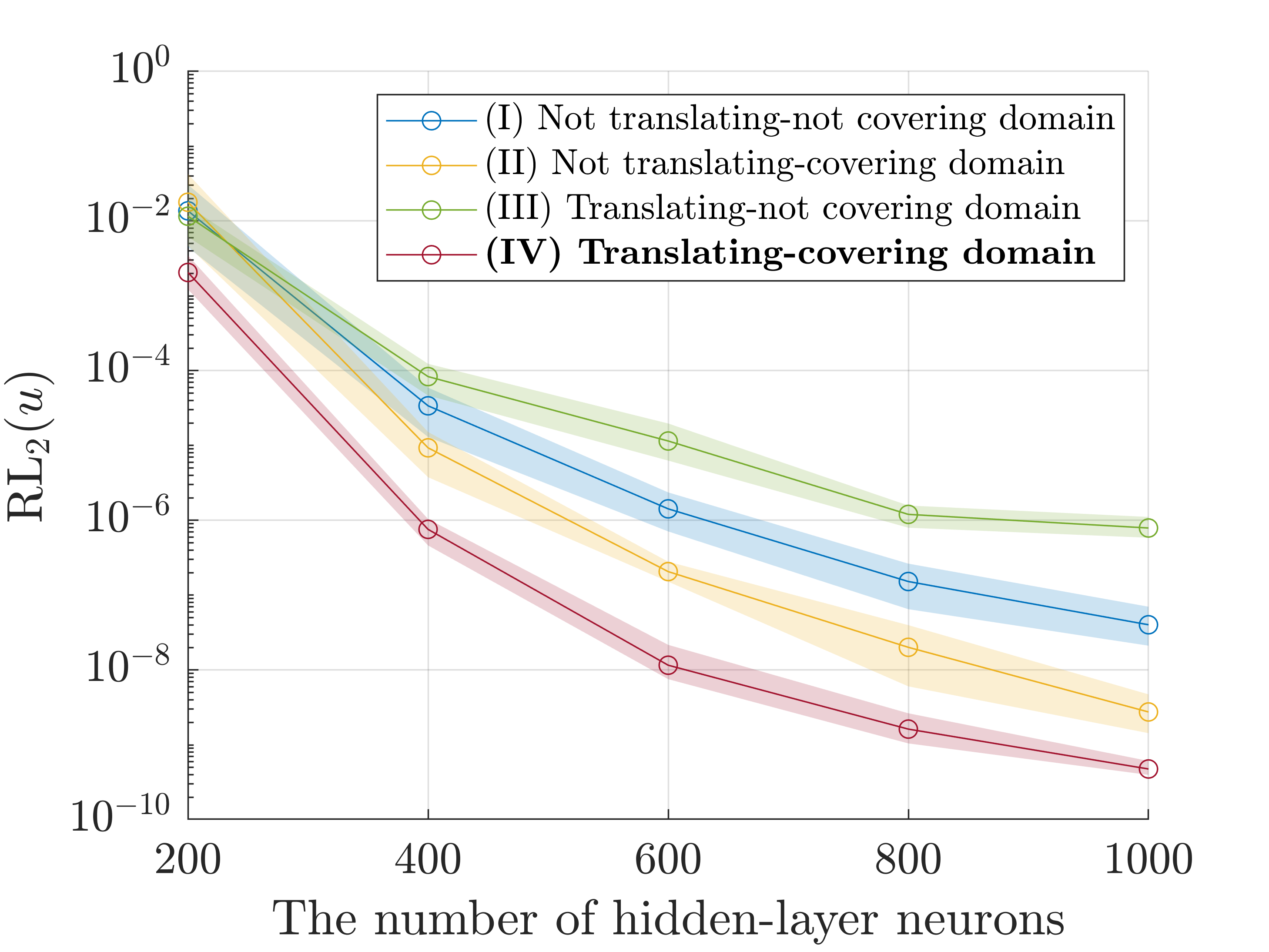}	
	\end{minipage}
	\begin{minipage}[t]{.48\textwidth}
		\centering
		\includegraphics[width=\linewidth]{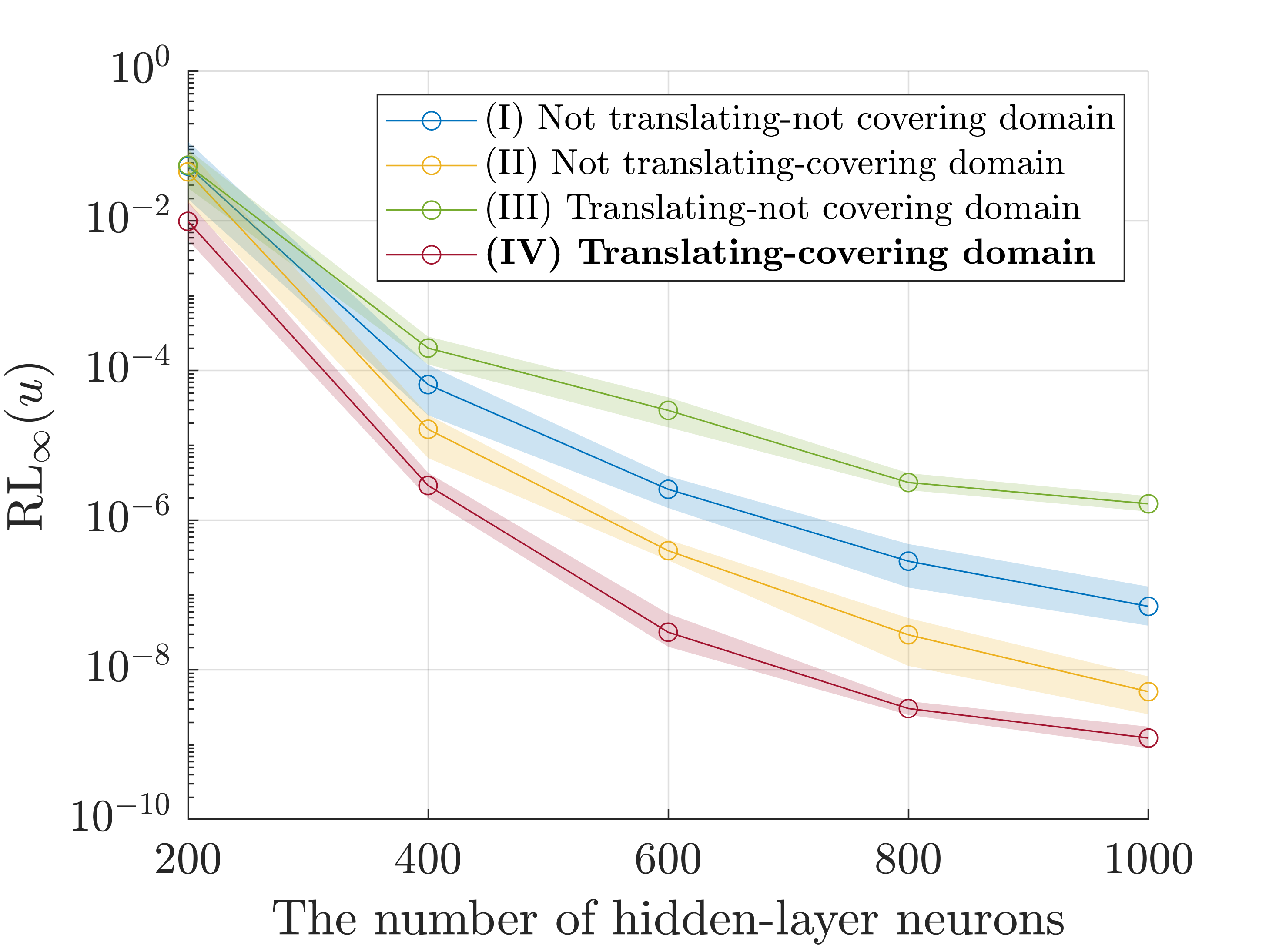}
	\end{minipage}
    \vspace{-2mm}
    \caption{Comparisons of the relative L$_2$ (left) and L$_\infty$ (right) errors of numerical solutions produced by the TransNets with different combinations of translating and scaling parameters (see \autoref{tab:TransNet-abla-2D-empi-const}) for the problem \eqref{equ:TransNet-abla-ex}.}
    \label{fig:TransNet-abla-2D-adap-rel-err}
\end{figure}

\paragraph*{Benefits of translating and scaling hidden-layer neurons}
We numerically investigate the effects of translating and scaling hidden-layer neurons of TransNet  discussed in \autoref{ssec:uniform-distri-arbi-ball}. For this propose, we design four combinations for $B_{R}(\boldsymbol{x}_c)$ according to whether to translate to the center of the domain $\Omega$ and whether to scale  to cover the domain: (I) Not translating-not covering domain, (II) Not translating-covering domain, (III) Translating-not covering domain, and (IV) Translating-covering domain. 
The specific settings for parameters $\boldsymbol{x}_c$ and $R$ are listed in \autoref{tab:TransNet-abla-2D-empi-const}, and their illustration is shown in \autoref{fig:TransNet-abla-2D-trans-scal}. 
For all  the four settings,  the number of hidden-layer neurons of the TransNet $M$ is gradually increased from 200 to 1000 by an increment of 200 each time, and the corresponding shape parameter $\gamma$ is obtained by using the training loss-based  optimization strategy for TransNet presented in \autoref{sec:TransNet-auto-tune-gamma}. All results of the relative numerical errors are shown in \autoref{fig:TransNet-abla-2D-adap-rel-err}. 
First, we observe that both the relative L$_2$ and L$_\infty$ errors of all four combinations rapidly and steadily decay as the number of hidden-layer neurons is growing. The setting of (IV) Translating-covering domain performs the best, in which the errors drop  from the level of $O(10^{-2})$ to 
 the level of $O(10^{-9})$ when $M$ increases from 200 to 1000. 
Subsequently,  with translating or not translating, the accuracy of the case of covering domain is remarkably better than that of the case of not covering domain, which is more and more noticeable along the growing of the number of hidden-layer neurons. These phenomena manifest that scaling  to ensure that the ball $B_{R}(\boldsymbol{x}_c)$  covers the domain is rather crucial to the performance of TransNet. 
Furthermore, with covering domain, the case of translating markedly outperforms the case of not translating. 
These results clearly demonstrate benefits of translating and scaling hidden-layer neurons  for TransNet.

\begin{table}[!ht]
	\centering
	\caption{Parameter settings corresponding to different combinations of translating and scaling for the problem \eqref{equ:TransNet-abla-ex}.}\label{tab:TransNet-abla-2D-empi-const}\footnotesize
    \vspace{2mm}
    \begin{tabular}{lcc}
		\toprule
		Translating \& scaling & $\boldsymbol{x}_{c}$ & $R$ \\
		\midrule
		(I) Not translating-not covering domain & (0, 0) & 2.5 \\
		(II) Not translating-covering domain & (0, 0) & 3.0   \\
		(III) Translating- not covering domain & (1, 1) & 1.0  \\
		\textbf{(IV) Translating-covering domain} & \textbf{(1, 1)} & \textbf{1.5}  \\
		\bottomrule
	\end{tabular}
\end{table}

\begin{figure}[!ht]
	\centering
	\begin{minipage}[t]{.48\textwidth}
		\centering
		\includegraphics[width=\linewidth]{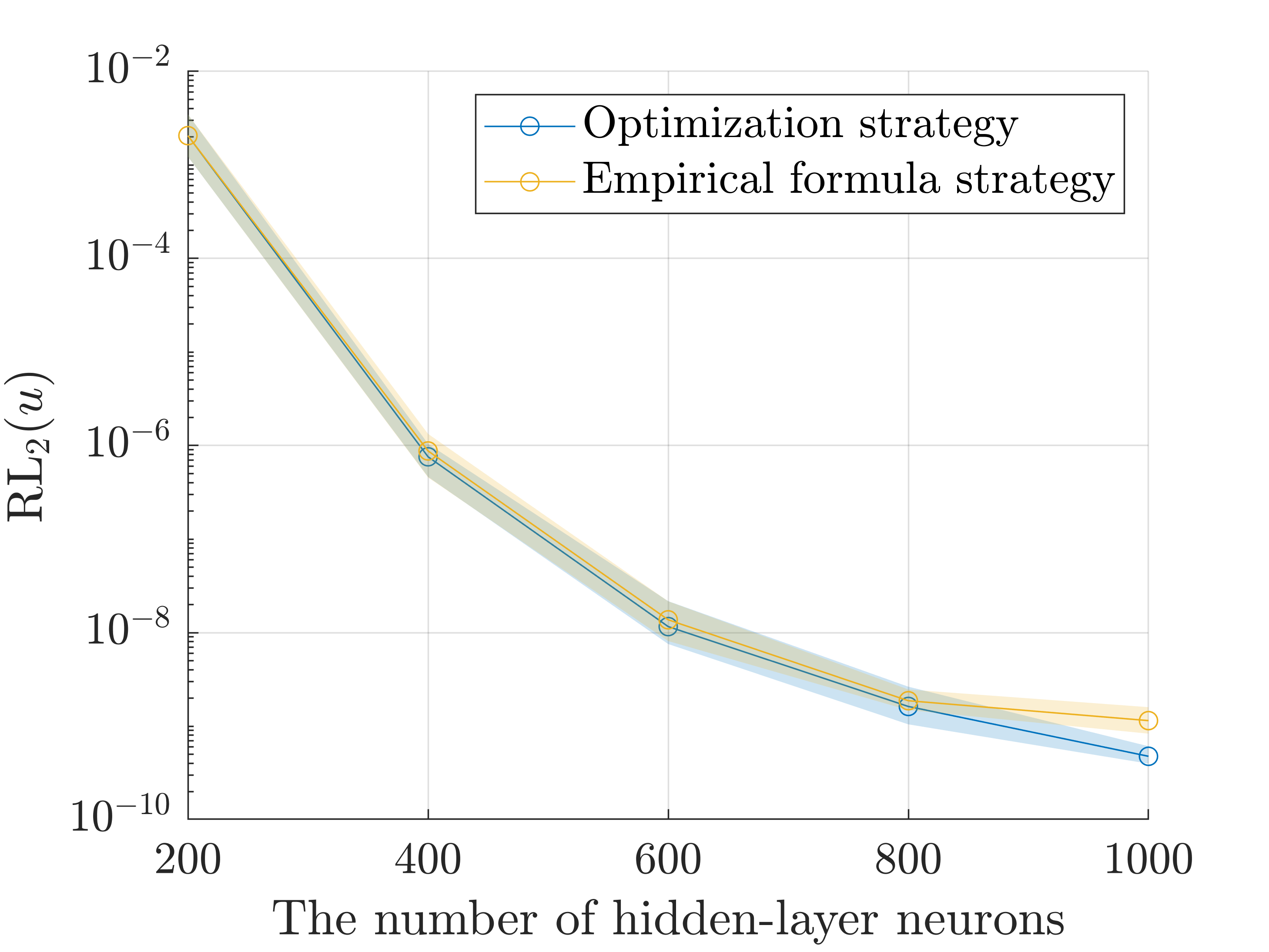}	
	\end{minipage}
	\begin{minipage}[t]{.48\textwidth}
		\centering
		\includegraphics[width=\linewidth]{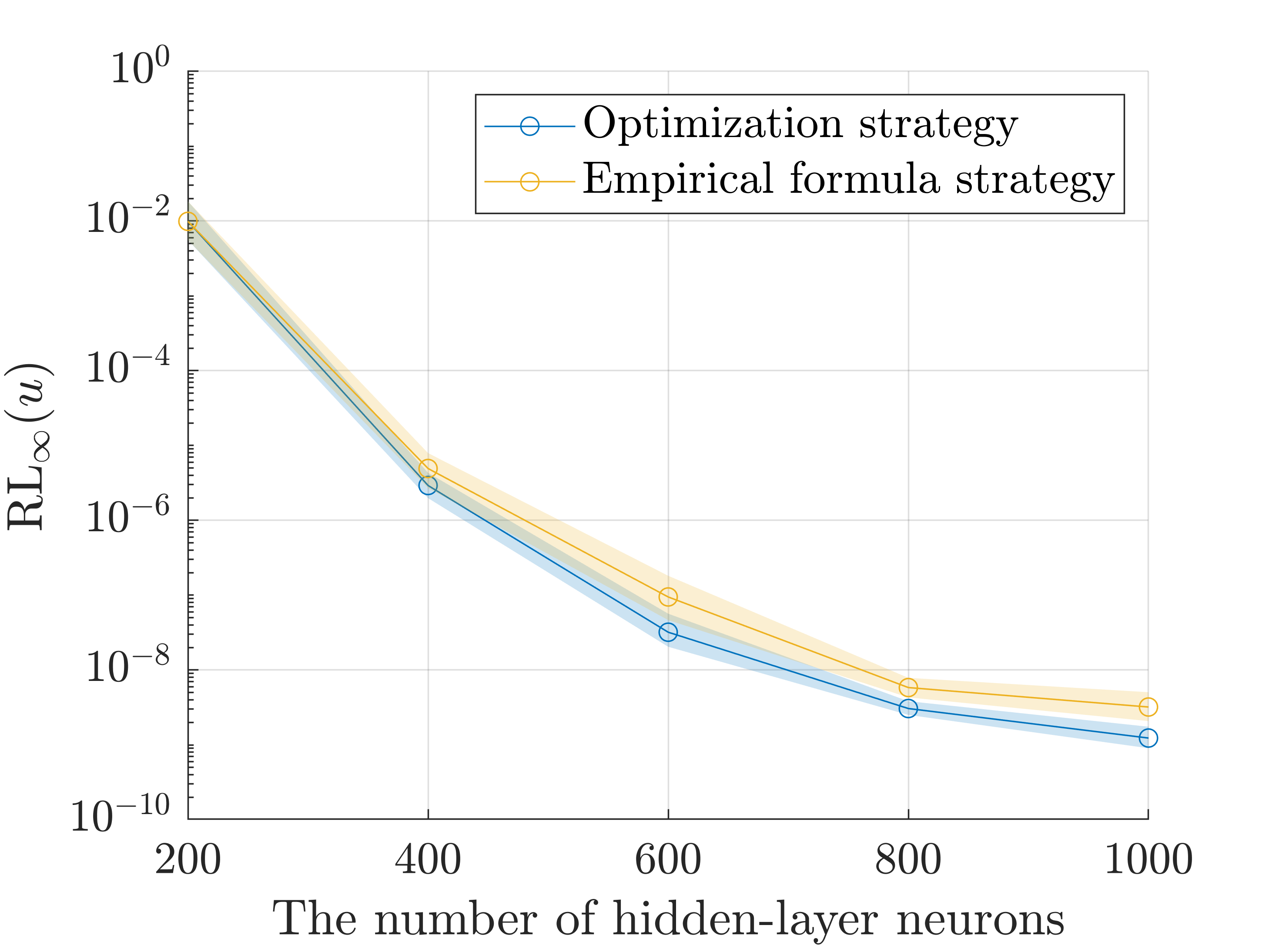}
	\end{minipage}
    \vspace{-2mm}
    \caption{Comparisons of the relative L$_2$ (left) and L$_\infty$ (right) errors of numerical solutions produced by the TransNets 
with the shape parameter $\gamma$ determined by the optimization strategy and by the empirical formula strategy for  the problem \eqref{equ:TransNet-abla-ex}.}\label{fig:TransNet-abla-2D-adap-vs-empi}
\end{figure}

\paragraph*{Effectiveness of the empirical formula-based prediction strategy}
We  compare the solution accuracy of the TransNets with the shape parameter $\gamma$ determined by the training loss-based optimization strategy and by the empirical formula-based prediction strategy (see \autoref{sec:TransNet-auto-tune-gamma}), respectively. The ball $B_{R}(\boldsymbol{x}_c)$ centered at ${\bm x}_c=(1,1)$ with radius $R=1.5$ (the setting of (IV) Translating-covering domain in \autoref{tab:TransNet-abla-2D-empi-const})  is used to cover the domain $\Omega$ and generate the hidden-layer neurons of the TransNet.
The number of hidden-layer neurons $M$ is again gradually increased from 200 to 1000 by an increment of 200 each time. 
For the training loss-based optimization strategy, the optimal shape parameter with respect to each number of hidden-layer neurons  is obtained by the 
golden-section search algorithm. 
For the empirical formula-based prediction strategy, the empirical constant $C$ in the  empirical formula \eqref{equ:emp-form} is first estimated using the optimization strategy for the TransNet with $200$ hidden-layer neurons in the preprocessing step (it is found $C\approx$ 8.4779e-2), and then the shape parameters of the TransNets with other numbers of hidden-layer neurons are automatically calculated  in a completely explicit manner. The comparison results on relative errors are shown in \autoref{fig:TransNet-abla-2D-adap-vs-empi}. It is seen that  the accuracy of the TransNets produced by the empirical formula strategy is similar to that of the TransNets by the optimization strategy and  marginally lower when the number of hidden-layer neurons reaches 1000.
Hence, the proposed empirical formula-based prediction strategy  is as effective as the training loss-based optimization strategy for TransNet, in addition to its
high efficiency.

\subsubsection{Multi-TransNet}\label{ssec:abla-Multi-TransNet}
In this subsection, we test the Multi-TransNet method
by  solving the following diffusion interface problem defined in the domain $\Omega=[0,2]^2$:
\begin{equation}\label{equ:abla-ex}
    \left\{\begin{aligned}
        -\nabla \cdot(\beta \nabla u) & = f, & & (x,y) \in \Omega_1 \cup \Omega_2, \\
        [u] & =h_1, & & (x,y) \in \Gamma, \\
        [\beta \nabla u \cdot \boldsymbol{n}] & =h_2, & & (x,y) \in \Gamma, \\
        u & =g, & & (x,y) \in \partial \Omega,
    \end{aligned}\right.
\end{equation}
where the interface is a circle $\Gamma=\{(x,y)\;|\;(x-1)^2 + (y-1)^2 = 0.5\}$. It divides $\Omega$ into two subdomains (i.e., $K=2$) $\Omega_1$ (inside) and $\Omega_2$ (outside). The exact solution is taken as $u(x,y) = \chi_{\Omega_1} \sin x \sin y  + \chi_{\Omega_2} \cos x \cos y $, and the diffusion coefficient is defined by the piecewise constant $\beta(x,y) = 1\chi_{\Omega_1} + 10\chi_{\Omega_2}$. The proposed Multi-TransNet method clearly will consist of two subdomain TransNets $u^\mathrm{NN}_1$ (associated to $\Omega_1$)
and $u^\mathrm{NN}_2$ (associated to $\Omega_2$). Inspired by the ablation studies for single TransNet in the previous subsection, we also apply  translating and scaling to the two subdomain TransNets 
so that the two balls  $B_{R_1}(\boldsymbol{x}^{(1)}_c)$ and $B_{R_2}(\boldsymbol{x}^{(2)}_c)$ can properly cover $\Omega_1$ and $\Omega_2$
respectively. To be detailed, the following  parameters are used: ${\bm x}^{(1)}_c=(1,1)$, $R_1=1$ and  ${\bm x}^{(2)}_c=(1,1)$, $R_2=1.5$.
We uniformly sample $49^2$ points in the interior, $51\times4$ points on the boundary, and $120$ points on the interface, with an overall spacing of approximately $0.04$, for training.

\begin{figure}[!ht]
	\centering
	\begin{minipage}[t]{.45\textwidth}
		\centering
		\includegraphics[width=\linewidth]{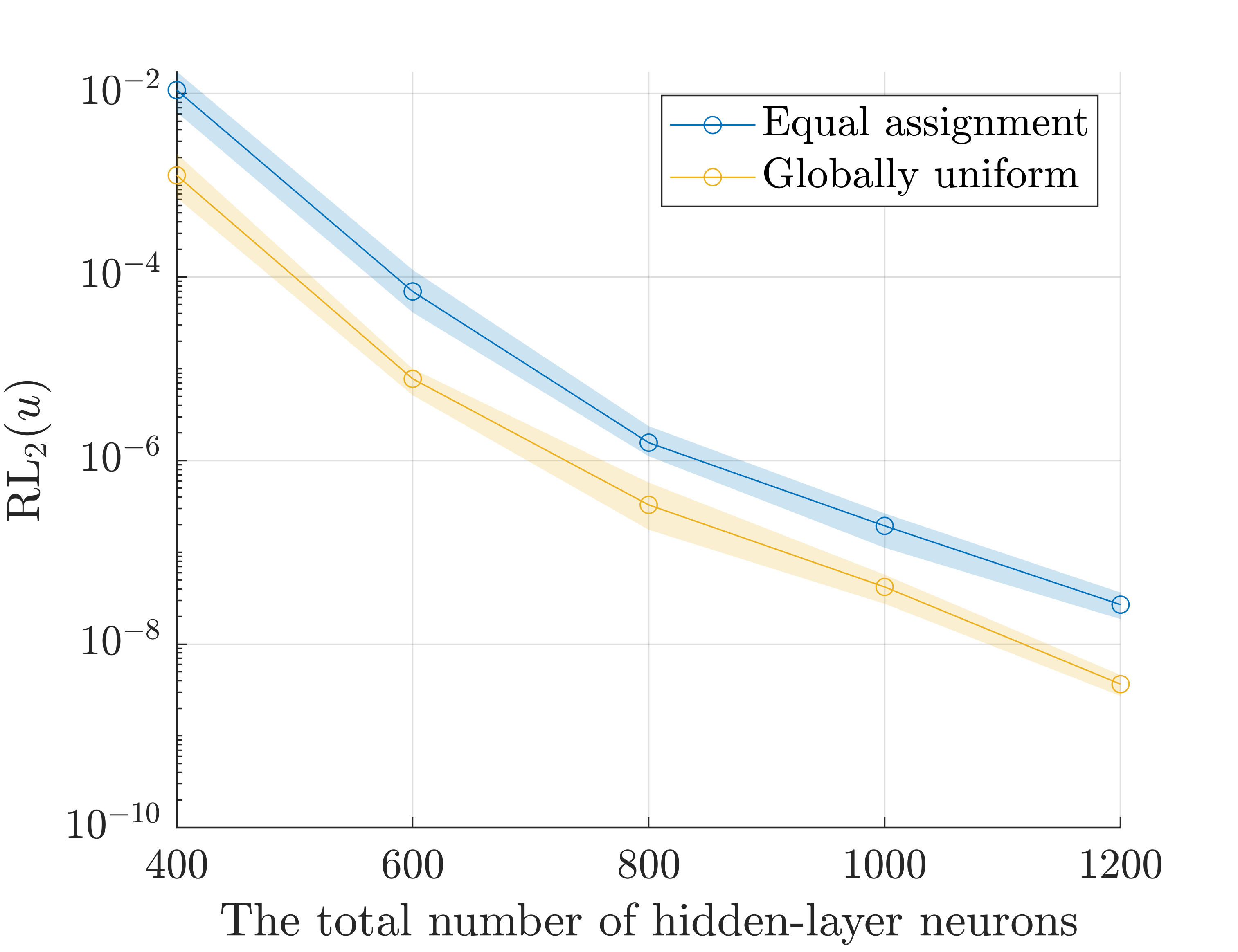}	
	\end{minipage}
	\begin{minipage}[t]{.45\textwidth}
		\centering
		\includegraphics[width=\linewidth]{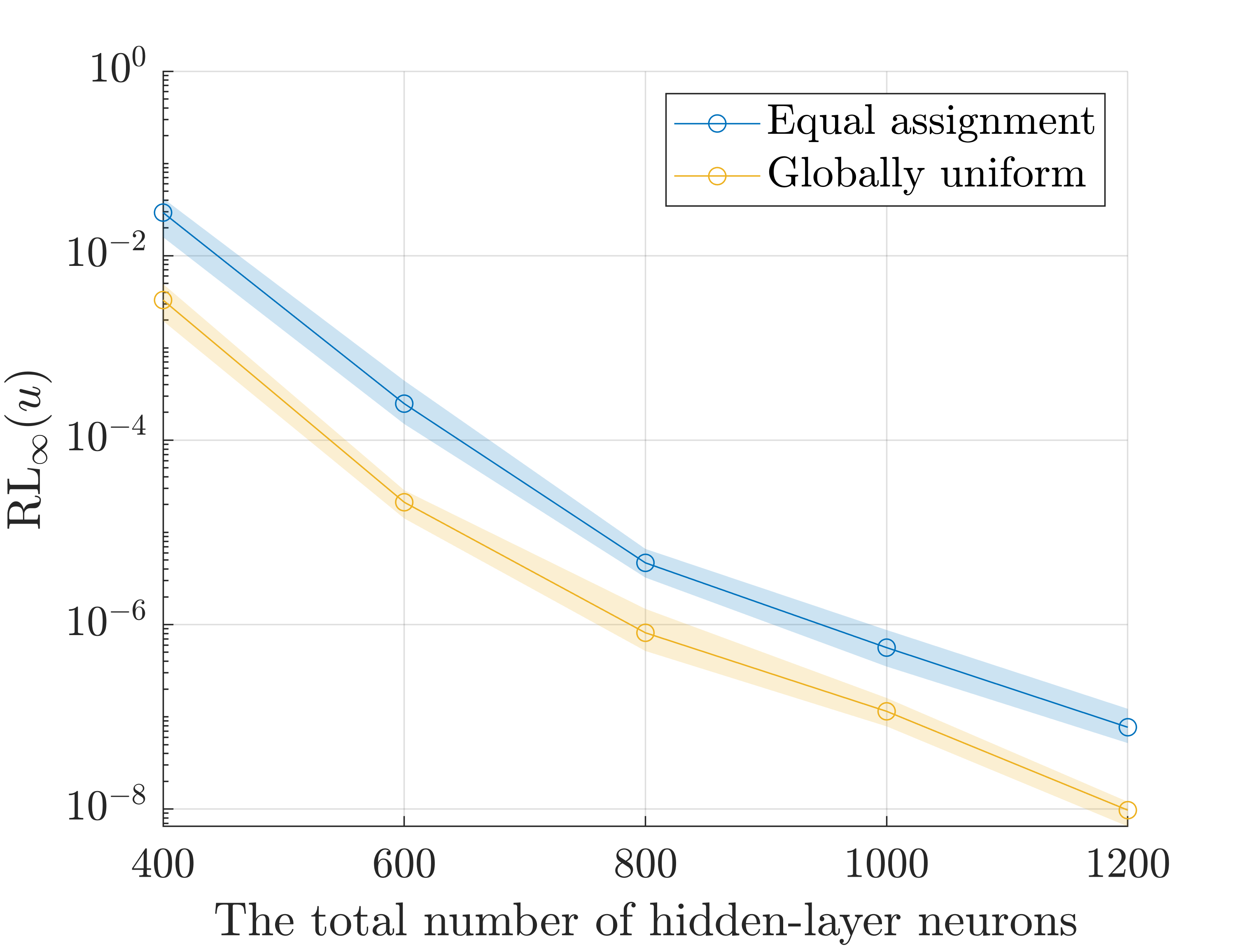}
	\end{minipage}
    \vspace{-2mm}
    \caption{Comparisons of the relative L$_2$ (left) and L$_\infty$ (right) numerical errors of solutions produced by the Multi-TransNet with the approach of equal assignment of the number of hidden-layer neurons among the subdomains and with the approach of globally uniform neuron distribution across the entire computational domain for the problem \eqref{equ:abla-ex}.}
    \label{fig:abla-2D-adap-glob-vs-loc}
\end{figure}

\paragraph*{Effect of globally uniform neuron distribution}
We  test and compare the performance of the Multi-TransNet method under the approach of  equal assignment of  the number of hidden-layer neurons among the subdomains (i.e., $M_1=M_2$) and  the approach of globally uniform neuron distribution across the entire computational domain (i.e., $M_1/R_1=M_2/R_2$).
The total number of hidden-layer neurons $M=M_1+M_2$ of the Multi-TransNet is gradually increased from 400 to 1200 by an increment of 200 each time, and corresponding shape parameters $(\gamma_1,\gamma_2)$ of the subdomain TransNets are obtained by using the  training loss-based optimization strategy 
for Multi-TransNet in \autoref{ssub-eef}.
The performance comparisons in terms of the relative errors are illustrated in \autoref{fig:abla-2D-adap-glob-vs-loc}.
It is clearly observed that the globally uniform neuron distribution approach  significantly prevails over the equal assignment approach (the errors of 
the former are approximately 10 times smaller than those of  the latter).
Although this is just one specific example, it implies that the globally uniform neuron distribution approach for Multi-TransNet is able to   bring better transferability and accuracy in general.

\begin{figure}[!ht]
	\centering
	\begin{minipage}[t]{.45\textwidth}
		\centering
		\includegraphics[width=\linewidth]{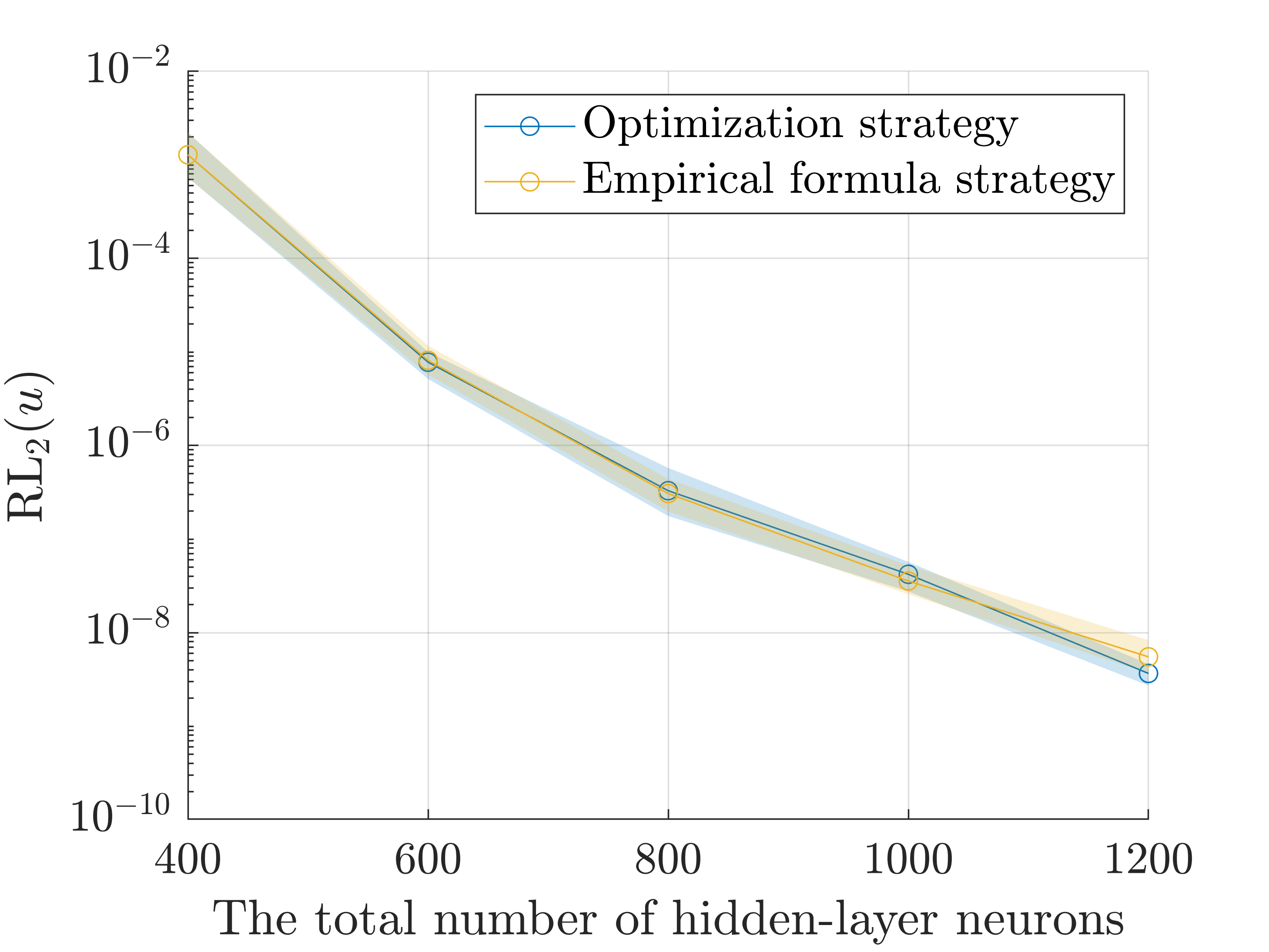}	
	\end{minipage}
	\begin{minipage}[t]{.45\textwidth}
		\centering
		\includegraphics[width=\linewidth]{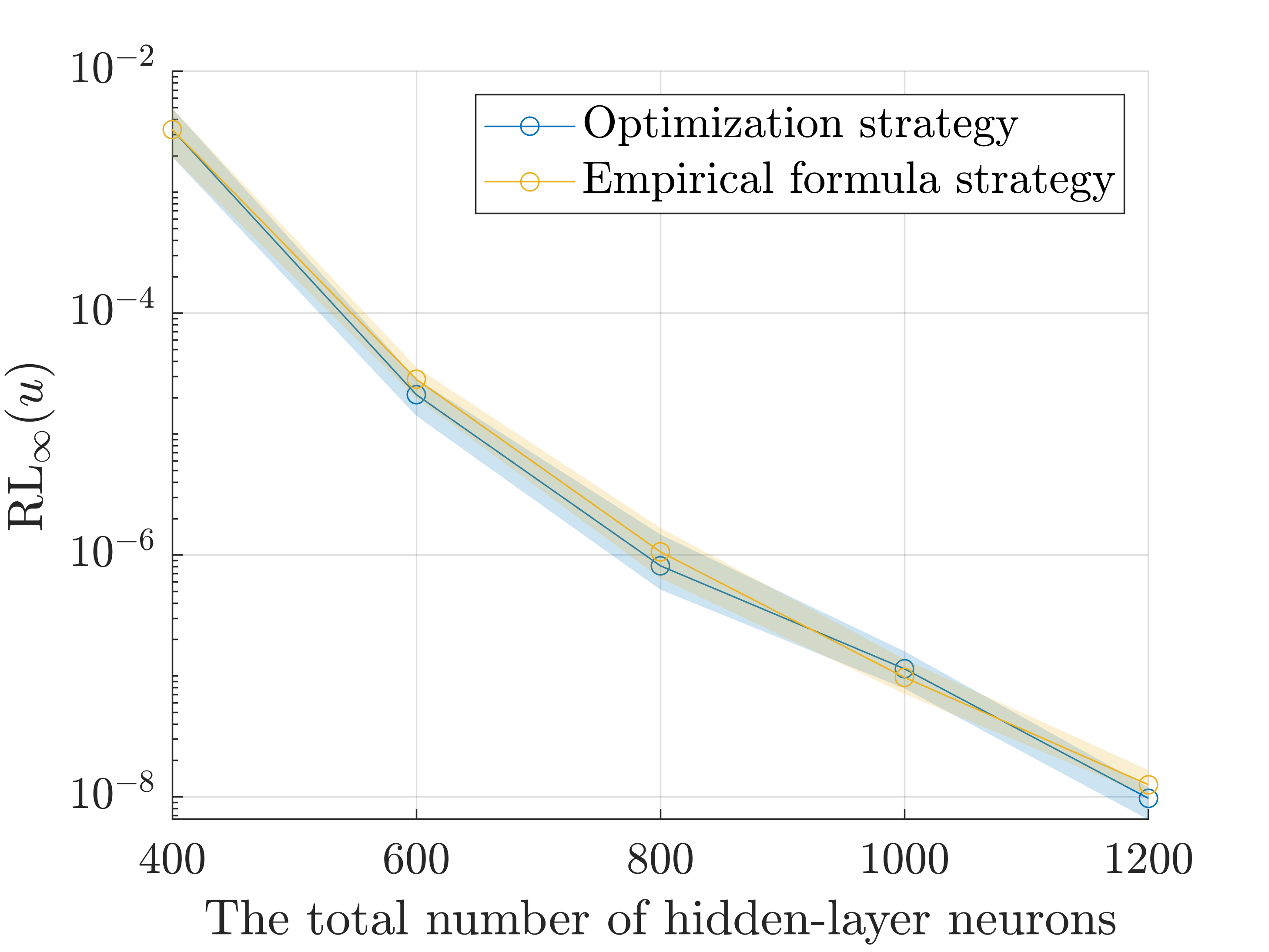}
	\end{minipage}
    \vspace{-2mm}
    \caption{Comparisons of the relative L$_2$ (left) and L$_\infty$ (right) errors of numerical solutions produced by the Multi-TransNet with the shape parameters $(\gamma_1,\gamma_2)$ determined by the optimization strategy and by the empirical formula strategy for the problem \eqref{equ:abla-ex}.}
    \label{fig:abla-2D-adap-vs-empi}
\end{figure}

\paragraph*{Effectiveness of the empirical formula-based prediction strategy}
We compare the solution accuracy of the
Multi-TransNet with the shape parameters $(\gamma_1,\gamma_2)$ determined by the training loss-based optimization strategy and the
empirical formula-based prediction strategy (see \autoref{ssub-eef}), respectively. The total number of hidden-layer neurons of the Multi-TransNet
$M$ is again gradually increased from 400 to 1200 by an increment of 200 each time, and  the globally uniform neuron distribution approach
is used for assignment of the numbers of hidden-layer neurons among the two subdomain TransNets.
For the training loss-based
optimization strategy, the optimal  shape parameters with respect to each total number of hidden-layer neurons
are obtained by the golden-section search algorithm. For the empirical formula-based prediction strategy,
an estimated value 7.3037e-2 for the empirical constant $C$ in the extended empirical formula \eqref{111} is first found by applying the optimization strategy to the
Multi-TransNet with totally 400 hidden-layer neurons in the preprocessing step, and then the shape parameters of
the Multi-TransNet with other total numbers of hidden-layer neurons are automatically calculated.
Their comparison results on  the relative  errors are shown in \autoref{fig:abla-2D-adap-vs-empi}.
It is observed that the accuracies produced by both strategies are nearly the same, which demonstrate that 
the proposed empirical formula-based prediction strategy for shape parameters   is effective and efficient for the Multi-TransNet method.

\begin{figure}[!htb]
	\centering
	\begin{minipage}[t]{.45\textwidth}
		\centering
		\includegraphics[width=\linewidth]{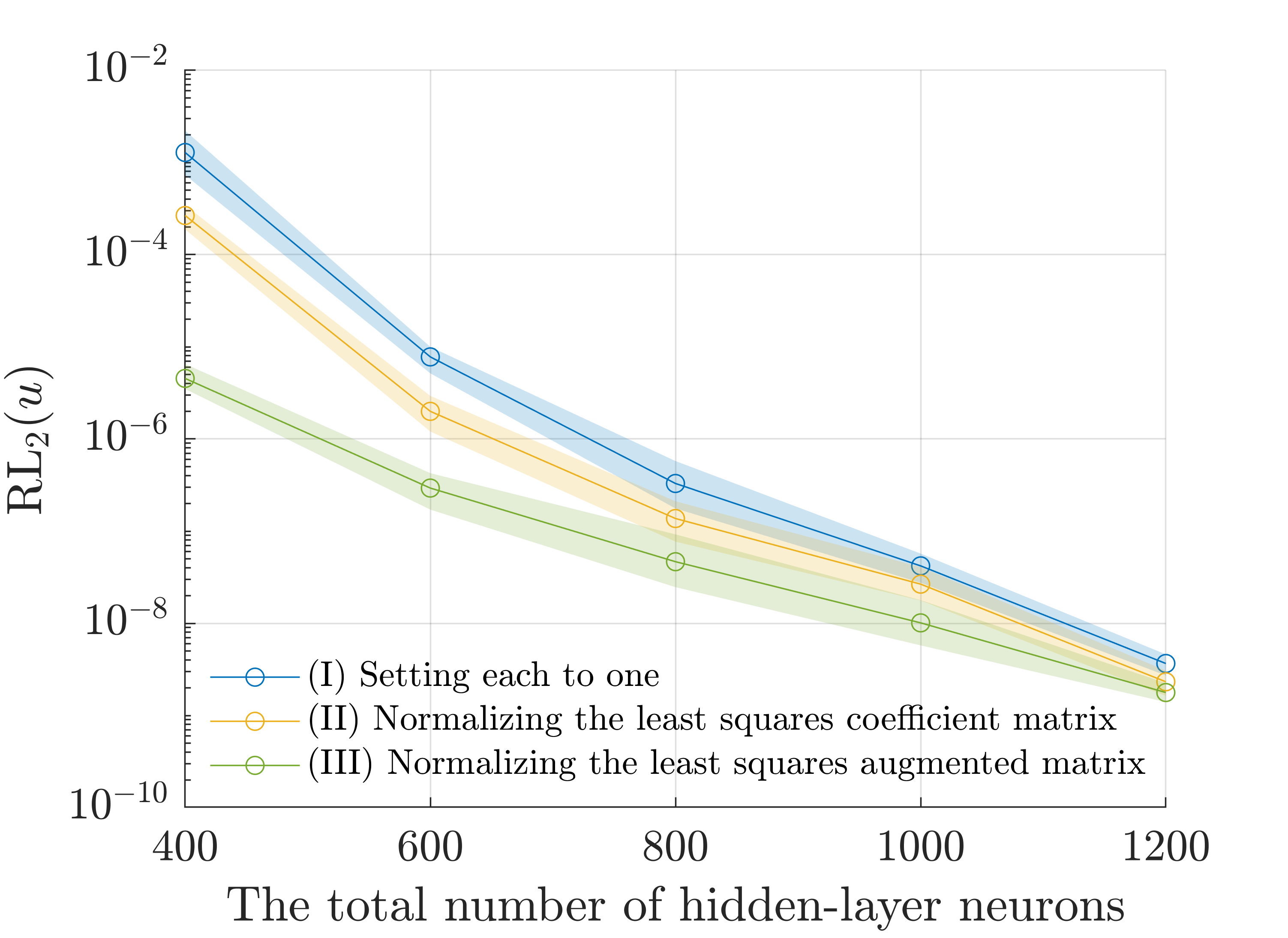}	
	\end{minipage}
	\begin{minipage}[t]{.45\textwidth}
		\centering
		\includegraphics[width=\linewidth]{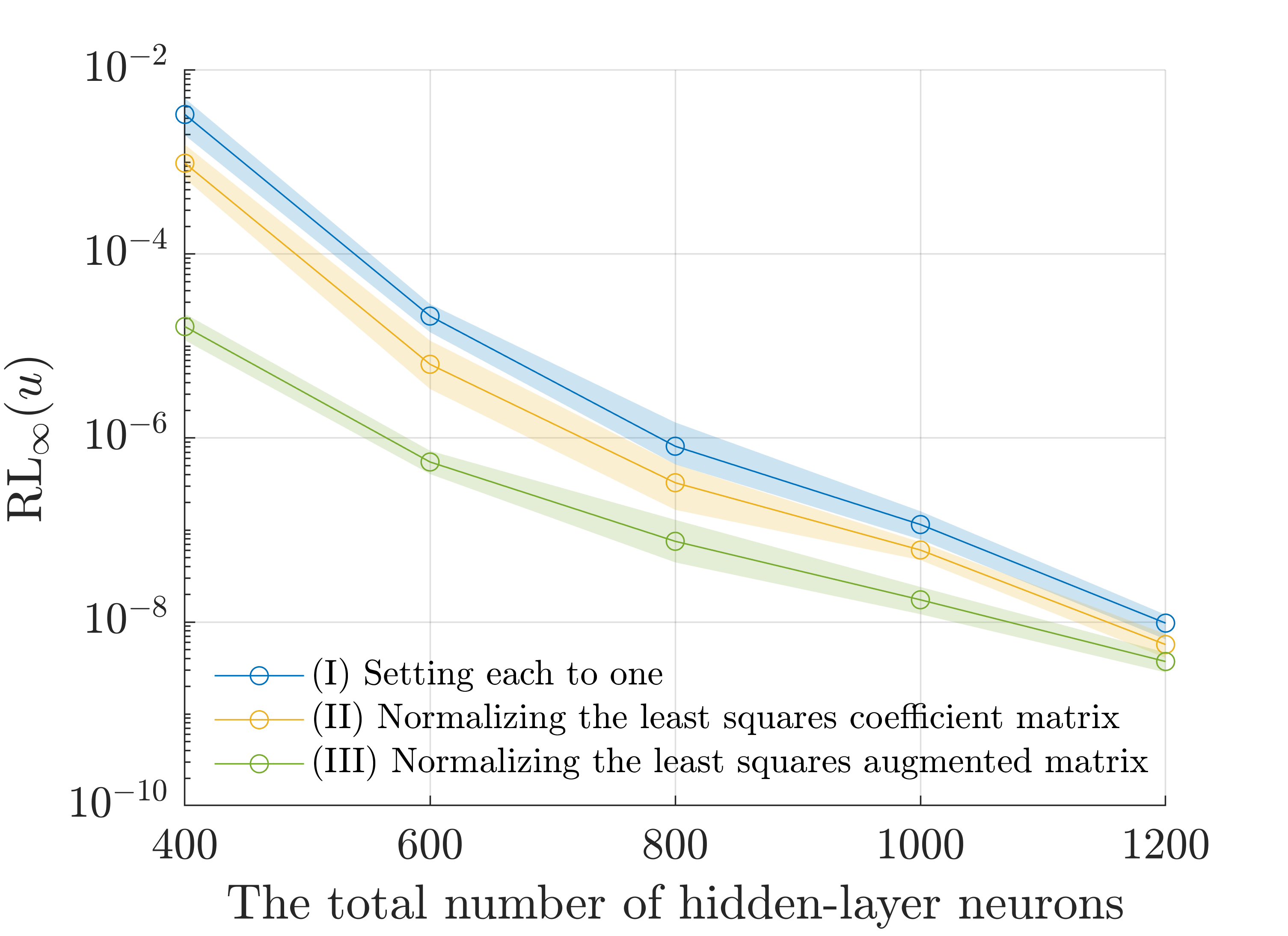}
	\end{minipage}
    \vspace{-2mm}
    \caption{Comparisons of the relative L$_2$ (left) and L$_\infty$ (right) errors of numerical solutions produced by the Multi-TransNet under three different choices of the loss weighting parameters for the problem \eqref{equ:abla-ex}.}
    \label{fig:abla-2D-penalty}
\end{figure}

\paragraph*{Impact of the weighting parameters in the loss function}
Regarding the impact of the weighting parameters of the loss function (i.e., $\lambda_1$, $\lambda_2$, $\lambda_g$, $\lambda_{h,1}$, $\lambda_{h,2}$ in \eqref{equ:lossmul})  to the Multi-TransNet solution,  
we compare three choices of determining the loss weighting parameters: (I) Setting each to one (i.e., equal weighting balances); (II) Normalizing the least squares coefficient matrix only, and (III) Normalizing the least squares augmented matrix, i.e., \eqref{equ:penalty}. 
Note that except the difference in the loss weighting parameters, all other parameter settings of Multi-TransNet are the same. For example,
 the globally uniform neuron distribution strategy is used for assignment of subdomain hidden-layer neurons and the shape parameters $(\gamma_1,\gamma_2)$ are determined by the training loss-based optimization strategy under the loss weighting  parameters choice (I).
 The comparison results are shown in \autoref{fig:abla-2D-penalty}.
 It is clearly observed that the relative errors produced under the choice (III) are always the best,  especially when the total number of hidden-layer neurons is relatively small. The performance of the choice (II) is slight better than that of the choice (I).

\subsection{Applications of the Multi-TransNet to typical elliptic interface problems}
In this subsection, we apply the proposed Multi-TransNet to a series of typical elliptic interface problems in two and three dimensions to demonstrate 
its outstanding performance  in terms of superior accuracy, efficiency and robustness, including a 2D Stokes interface problem with a circular interface, a  2D diffusion interface problem with multiple interfaces,   a 3D elasticity interface problem with an ellipsoidal interface, and a 3D diffusion interface problem with a convoluted immersed interface. The default configuration of the following Multi-TransNet is a combination of properly translating and scaling hidden-layer neurons based on subdomains, the globally uniform neuron distribution, the empirical formula-based prediction strategy for the shape parameters, and the normalization approach for weighting parameters in the loss function. 

\subsubsection{A 2D Stokes interface problem with a circular interface}\label{ssec:stokes-interf}
\begin{figure}[!t]
	\centering
	\includegraphics[width=\textwidth]{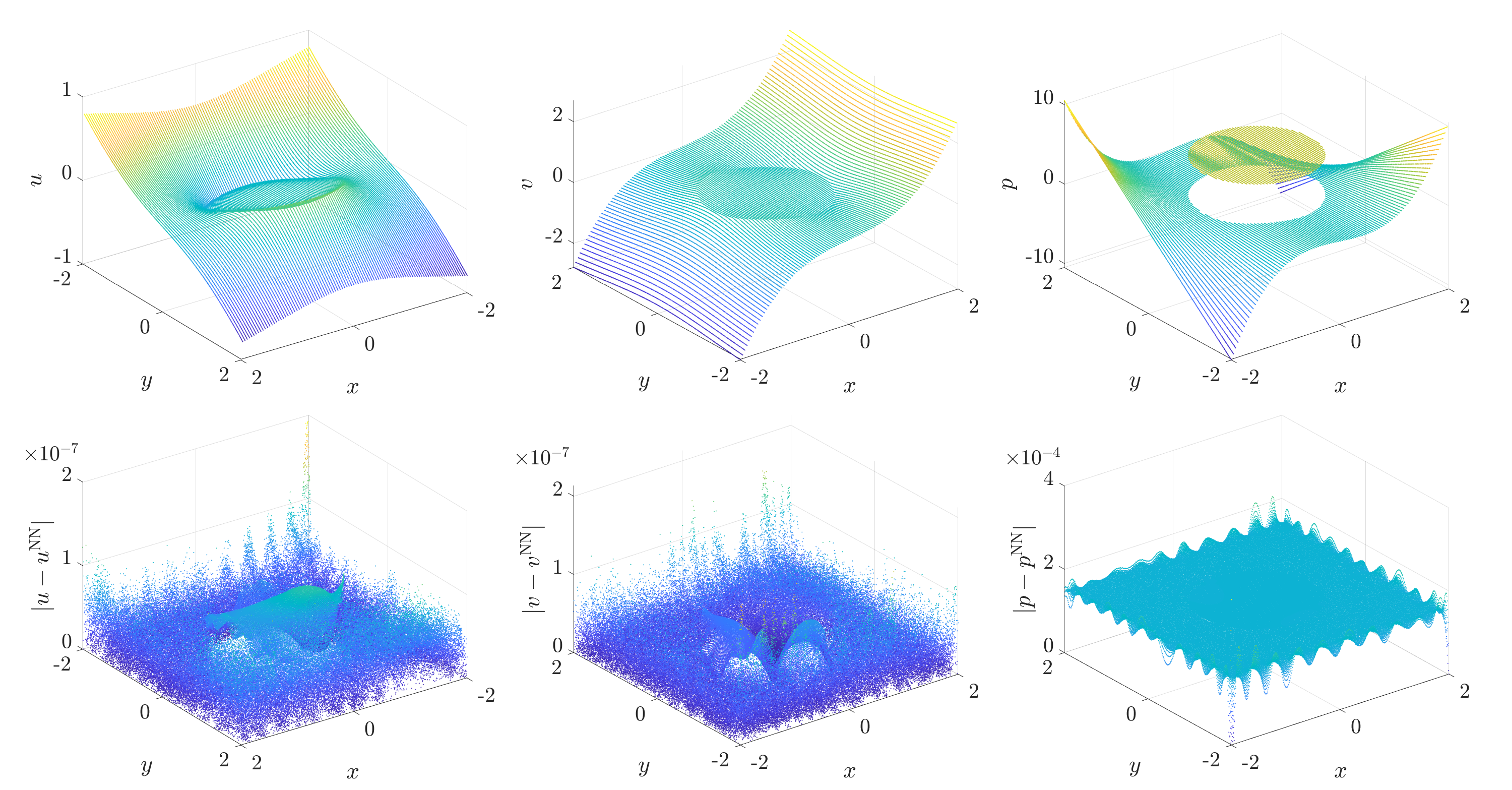}
	\vspace{-6mm}
	\caption{Top row: the exact solution of the 2D Stokes interface problem \eqref{equ:stokes-interf}; Bottom row: the pointwise absolute errors of the numerical  solution,  when the viscosity $(\mu_1, \mu_2)=(1,10)$, produced by the Multi-TransNet with totally $1000$ hidden-layer neurons.
	From left to right:  $u$, $v$  and $p$.}
	\label{fig:stokes-interf-exasol-abserr}
\end{figure}

In this example, we consider the two-phase Stokes interface problem \cite{dong2023kernel} in fluid mechanics field. The  domain $\Omega=[-2,2]^2$ is divided into two subdomains $\Omega_1$ (the inside one) and $\Omega_2$ (the outside one) by a circular interface $\Gamma=\{(x,y)\;|\;x^2+y^2=1\}$. The problem is described as follows:

\begin{equation}\label{equ:stokes-interf}
\left\{\begin{aligned}
    -\mu \Delta \boldsymbol{u}+\nabla p & =\boldsymbol{f}, & & (x,y) \in \Omega_1 \cup \Omega_2, \\
    \nabla \cdot \boldsymbol{u} & =0, & & (x,y) \in \Omega_1 \cup \Omega_2, \\
    \left[ \boldsymbol{u} \right] &=\boldsymbol{0}, & & (x,y) \in \Gamma, \\
    \left[ \boldsymbol{\sigma}(\boldsymbol{u}, p) \boldsymbol{n} \right] & =\boldsymbol{h}_2, & & (x,y) \in \Gamma, \\
    \boldsymbol{u} & =\boldsymbol{g}, & & (x,y) \in \partial \Omega,
\end{aligned}\right.
\end{equation}
where the viscosity $\mu(x,y) = \mu_k$ for $(x,y)\in \Omega_k$ ($k=1,2$), the velocity vector $\boldsymbol{u}=\left(u,v\right)^\top$ and the stress tensor $\boldsymbol{\sigma}(\boldsymbol{u},p) = -p\boldsymbol{I} + \mu \left(\nabla \boldsymbol{u} + \nabla \boldsymbol{u}^\top\right)$. Its exact solution is given by
\begin{equation}\label{equ:stokes-interf-exasol}
    \begin{aligned}
        u(x, y) & = 
        \begin{cases}
            \frac{y}{4}\left(x^2+y^2\right), &\qquad\qquad (x, y) \in \Omega_1, \\
            \frac{y}{\sqrt{x^2+y^2}} - \frac{3y}{4}, &\qquad\qquad (x, y) \in \Omega_2,
        \end{cases} \\
        v(x, y) & = 
        \begin{cases}
            -\frac{x y^2}{4}, &\, (x, y) \in \Omega_1, \\
            -\frac{x}{\sqrt{x^2+y^2}} + \frac{x}{4}\left(3+x^2\right), &\, (x, y) \in \Omega_2,
        \end{cases} \\
        p(x, y) & = 
        \begin{cases}
            5.0, &\qquad\quad\   (x, y) \in \Omega_1, \\
            \left(-\frac{3}{4} x^3+\frac{3}{8} x\right) y, &\qquad\quad\   (x, y) \in \Omega_2.
        \end{cases}
    \end{aligned}
\end{equation}
Note that  $u, v$ are continuous but nonsmooth across the interface while $p$ is discontinuous as shown in the top row of \autoref{fig:stokes-interf-exasol-abserr}. The {viscosity} $\mu_1=1$ is fixed and  different values  will be chosen for $\mu_2$ to reflect various ratio contrast (i.e., ${\mu_1}/{\mu_2}$) cases.
We  sample $N_f = 149^2$ (in the interior), $N_g=301\times 4$ (on the boundary) and $N_\Gamma=1000$ (on the interface) training points, and
 use two subdomain TransNets with a combination of translating and scaling parameters $(\boldsymbol{x}^{(1)}_{c}, R_1) = ((0.0, 0.0), 1.25)$, $(\boldsymbol{x}^{(2)}_{c}, R_2) = ((0.0, 0.0), 3.0)$ for the Multi-TransNet method to approximate the solution of the problem \eqref{equ:stokes-interf}.

\begin{figure}[!ht]
	\centering
	\includegraphics[width=\textwidth]{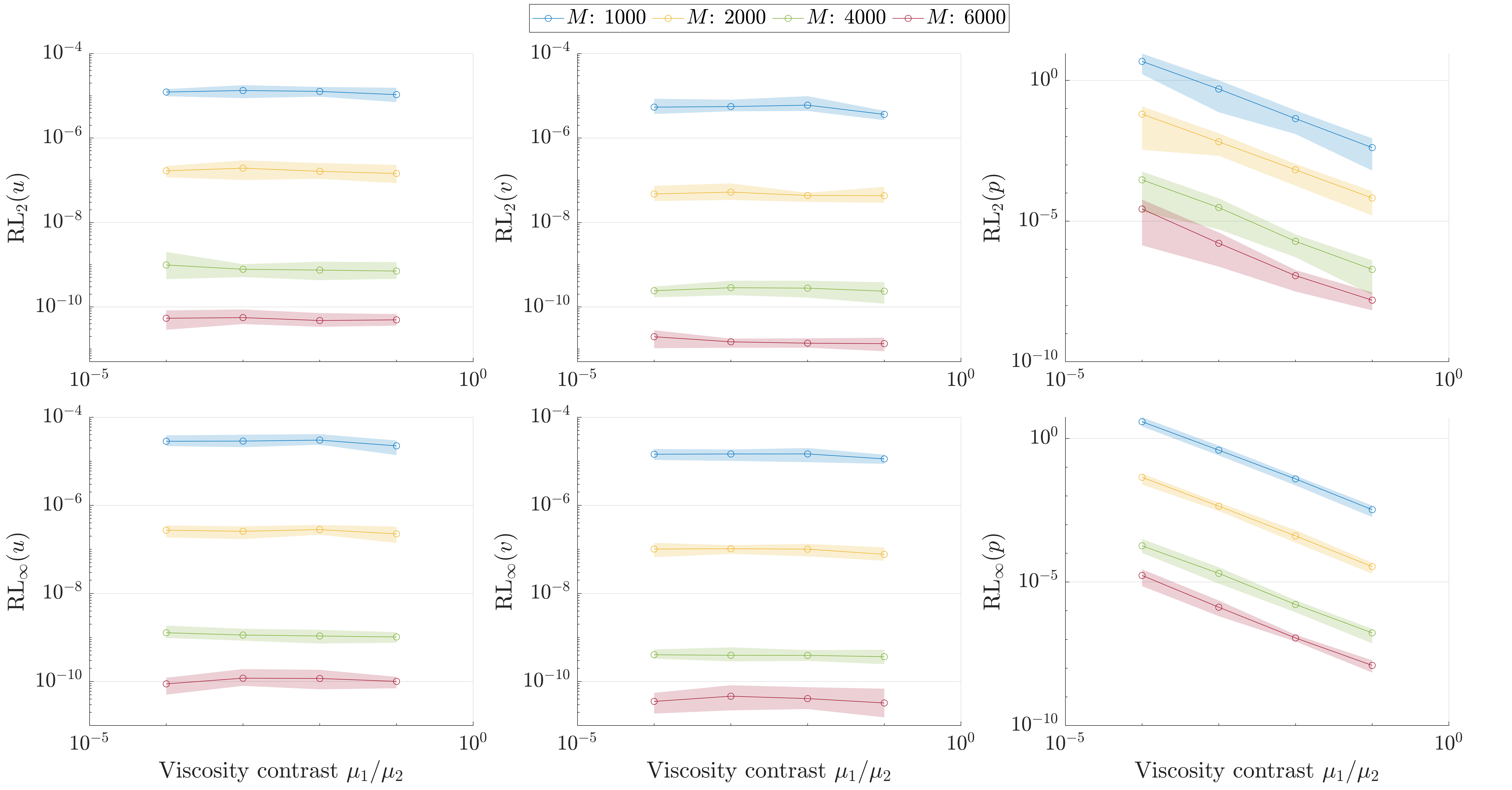}
	\vspace{-6mm}
	\caption{The relative L$_2$ (top row) and L$_\infty$ (bottom row) errors of the numerical solutions produced by the Multi-TransNet method with different total  numbers of hidden-layer neurons ($M = 1000, 2000, 4000, 6000$) for the  2D Stokes interface problem \eqref{equ:stokes-interf} with different viscosity contrasts ({$\mu_1/\mu_2 = 10^{-1}, 10^{-2}, 10^{-3}, 10^{-4}$}). From left to right:  $u$, $v$  and $p$.}
	\label{fig:stokes-accu-contrast-M}
\end{figure}

\begin{figure}[!ht]
	\centering
	\includegraphics[width=0.55\textwidth]{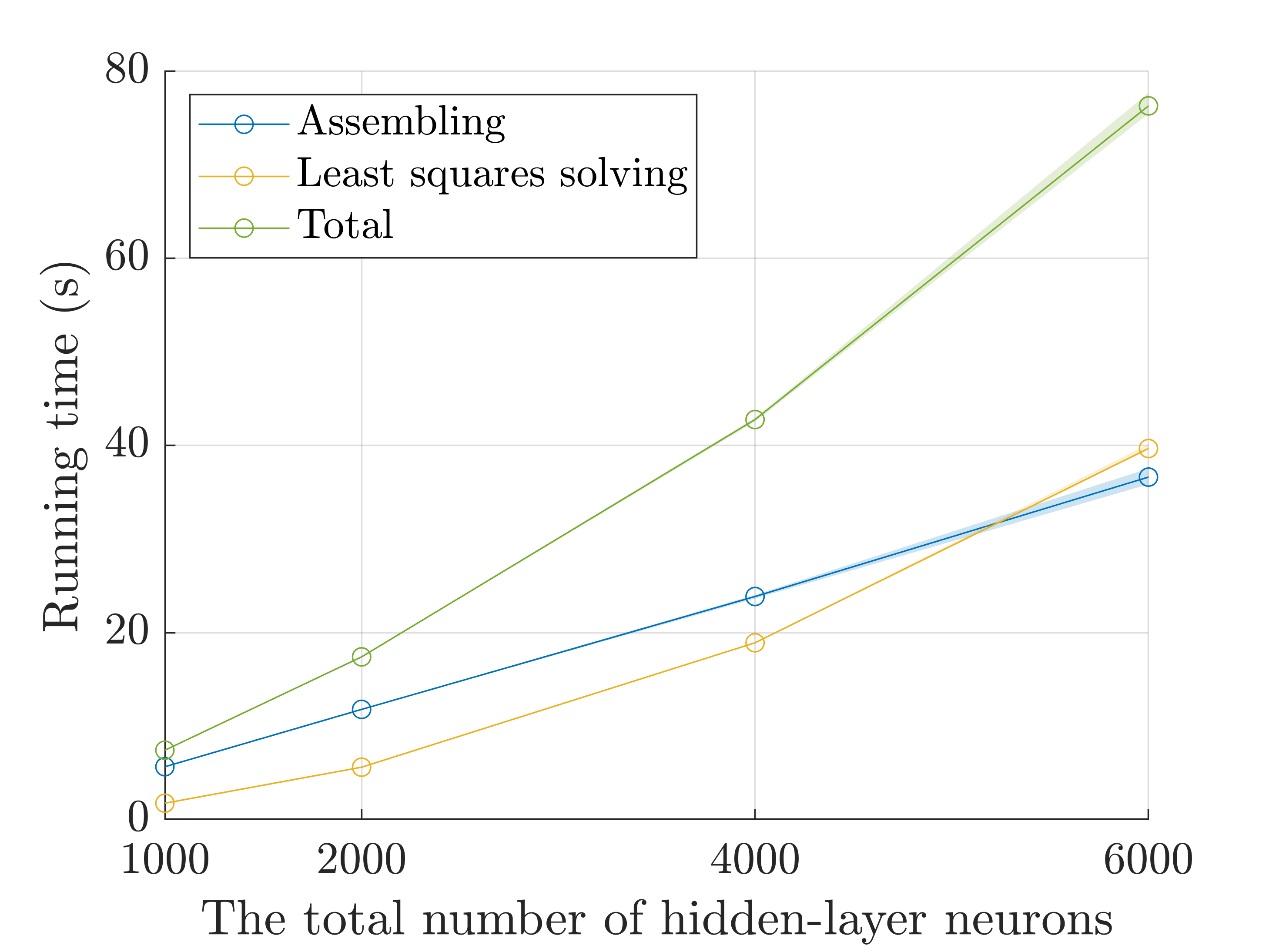}
	\caption{The average running times (measured in seconds) of the Multi-TransNet method with different total  numbers of hidden-layer neurons ($M = 1000, 2000, 4000, 6000$) for the 2D Stokes interface problem \eqref{equ:stokes-interf}.}
	\label{fig:stokes-rt}
\end{figure}

The bottom row of \autoref{fig:stokes-interf-exasol-abserr} illustrates the pointwise   absolute  errors of the numerical solution produced by the Multi-TransNet method with totally $M=M_1+M_2=1000$ hidden-layer neurons under the training loss-based optimization strategy for optimal  shape parameters $(\gamma_1,\gamma_2)$ in the low-contrast case {of the viscosity $\mu_2=10$}. Meanwhile, an estimated value 8.6107e-2 is obtained for the empirical constant $C$ in the extended empirical formula \eqref{111} for Multi-TransNet. 
Subsequently,  the empirical formula-based prediction strategy is employed to determine appropriate shape parameters for the Multi-TransNet with other numbers of hidden-layer neurons, which will also be used to  solve other viscosity contrast cases. The numerical results for $\mu_1/\mu_2=10^{-1}$, $10^{-2}$, $10^{-3}$ and $10^{-4}$ are illustrated in \autoref{fig:stokes-accu-contrast-M}, from which it is clearly observed that the relative L$_2$ and L$_\infty$ errors of the Multi-TransNet solutions quickly decline when the total number of hidden-layer neurons gradually increases from $M=1000, 2000, 4000$ to $M=6000$, and finally reach a level of $O(10^{-10})$ for the velocity components $u$ and $v$.
Furthermore, the relative  errors of both velocity components remain flat regardless of the viscosity contrast, while that of the pressure decreases with the lowering of the viscosity contrast (note similar behaviors were also observed on the pressure solutions by the traditional numerical method \cite{dong2023kernel} and the Random Feature Method (RFM) \cite{chi2024random} for this example). \autoref{fig:stokes-rt} reports the average running times of the Multi-TransNet method with respect to different total numbers of hidden-layer neurons. We observe that the system assembling time seems  linearly proportional to the total number of hidden-layer neurons, but  the least squares solving time  increases much faster along with the growing of the total number of hidden-layer neurons.

Finally, we compare the performance of the Multi-TransNet with the RFM in terms of accuracy and efficiency, and the results are listed in \autoref{tab:stokes-interf-vs-rfm}. We observe that both neural network methods exhibit good robustness for different viscosity contrasts. Nonetheless, it is easy to see that the Multi-TransNet strikingly outperforms the RFM regardless of the number of constraints (i.e., the row number of the resulting least squares system) and the number of DOFs (i.e., the column number of the resulting least squares system,  or the product of the total number of hidden-layer neurons and the number of unknown variables). In particular, the relative L$_2$ errors of the Multi-TransNet with even smaller number of DOFs are several magnitudes smaller than those of the RFM for the velocity components $u$ and $v$.

\begin{table}[!ht]
	\centering
	\caption{Comparison of the performance of the Multi-TransNet and the RFM \cite{chi2024random} for the 2D Stokes interface problem \eqref{equ:stokes-interf} with different viscosity contrasts  ({$\mu_1/\mu_2 = 10^{-1}, 10^{-2}, 10^{-3}$}). Note that \#DOFs = 3$M$ for the  Multi-TransNet method in this problem.}
    \vspace{0.2cm}
	\label{tab:stokes-interf-vs-rfm}
	{\footnotesize
		\begin{tabular}{cccrlll}
        \toprule
        ${\mu_1}/{\mu_2}$ & Method & \#Constraints & \#DOFs & RL$_2 (u)$ & RL$_2 (v)$ & RL$_2 (p)$ \\
        \midrule
        \multirow{4}{*}{$10^{-1}$} 
        & \multirow{2}{*}{RFM} & \multirow{2}{*}{84,801}          
          & 9,600  & 6.92e-6 & 1.71e-6 & 3.33e-4 \\
        &&& 38,400 & 1.18e-8 & 3.30e-9 & 4.84e-8 \\
        \cmidrule{2-7}
        & \multirow{2}{*}{Multi-TransNet} & \multirow{2}{*}{73,012} 
          & 6,000  & 1.39e-7 & 4.04e-8 & 6.45e-5 \\
        &&& 18,000 & 4.88e-11 & 1.33e-11 & 1.43e-8 \\
        \cmidrule{1-7}
        \multirow{4}{*}{$10^{-2}$} 
        & \multirow{2}{*}{RFM} & \multirow{2}{*}{84,801}
          & 9,600  & 4.93e-6 & 1.48e-6 & 2.63e-3 \\
        &&& 38,400 & 1.52e-8 & 3.21e-9 & 5.99e-7 \\
        \cmidrule{2-7}
        & \multirow{2}{*}{Multi-TransNet} & \multirow{2}{*}{73,012}
          & 6,000  & 1.58e-7 & 4.52e-8 & 6.89e-4 \\
        &&& 18,000 & 4.59e-11 & 1.37e-11 & 1.20e-7 \\
        \cmidrule{1-7}
        \multirow{4}{*}{$10^{-3}$} 
        & \multirow{2}{*}{RFM} & \multirow{2}{*}{84,801}
          & 9,600  & 4.49e-6 & 1.44e-6 & 3.31e-2 \\
        &&& 38,400 & 5.97e-9 & 1.64e-9 & 2.76e-6 \\
        \cmidrule{2-7}
        & \multirow{2}{*}{Multi-TransNet} & \multirow{2}{*}{73,012}
          & 6,000  & 1.90e-7 & 5.02e-8 & 6.16e-3 \\
        &&& 18,000 & 5.12e-11 & 1.51e-11 & 1.50e-6 \\
        \bottomrule
        \end{tabular}
	}
\end{table}

\subsubsection{A 2D diffusion interface problem with multiple interfaces}\label{ssec:multi-interf}
\begin{figure}[!ht]
	\centering
	\begin{minipage}[t]{.32\textwidth}
		\centering
		\includegraphics[width=\linewidth]{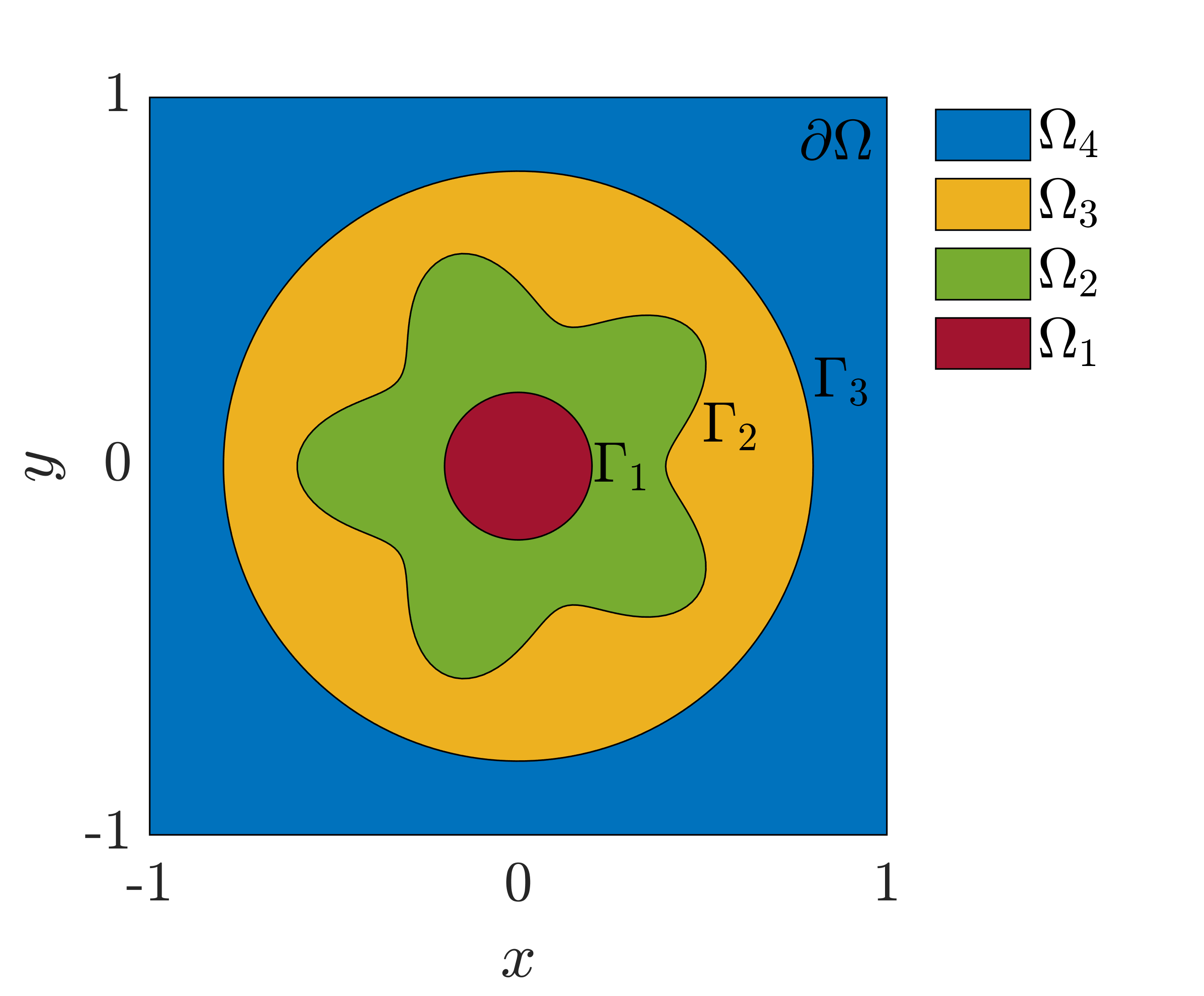}	
	\end{minipage}
	\begin{minipage}[t]{.33\textwidth}
		\centering
		\includegraphics[width=\linewidth]{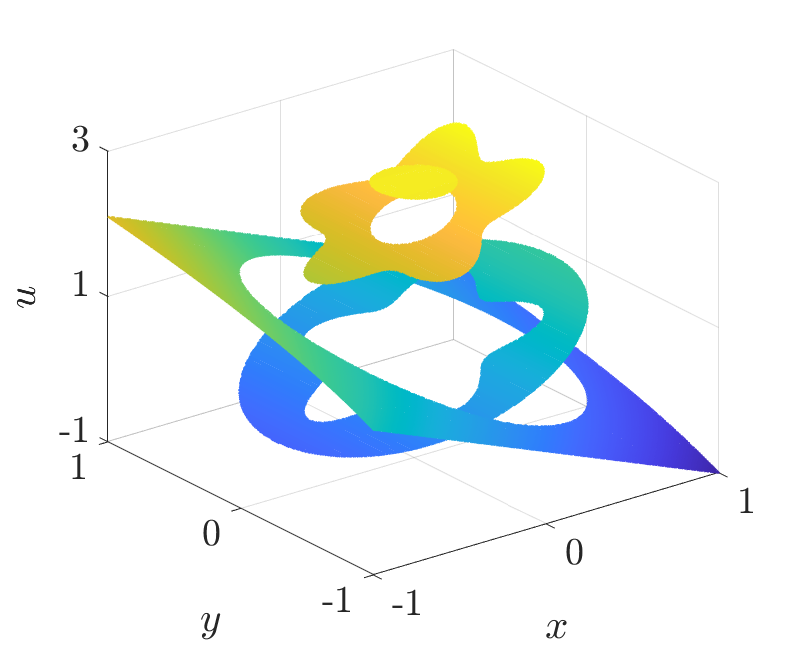}
	\end{minipage}
    \begin{minipage}[t]{.33\textwidth}
		\centering
		\includegraphics[width=\linewidth]{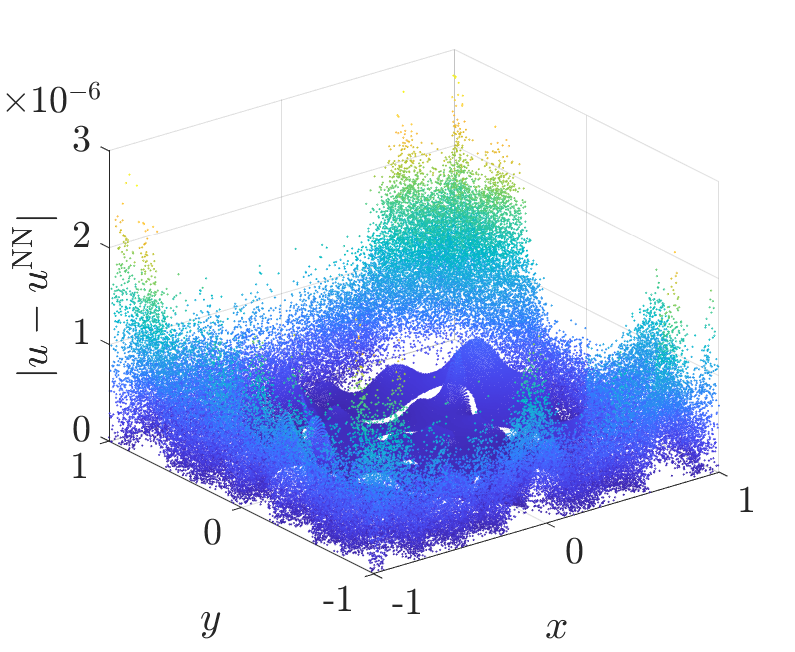}
	\end{minipage}
    \vspace{-2mm}
    \caption{Left: illustration of the partition of the domain into four subdomains by three interfaces; Middle: the exact solution of the 2D diffusion interface \eqref {equ:multi-interf}; Right: the pointwise absolute errors of numerical solution, when the diffusion coefficient {$(\beta_1, \beta_2, \beta_3, \beta_4)=(1,10,10^2,10^3)$}, produced by the Multi-TransNet method with totally 400 hidden-layer neurons.}\label{fig:multi-interf-domain-exasol-abserr}
\end{figure}

In this example, we consider the two-dimensional diffusion interface problem with multiple interfaces \cite{guittet2015solving}. The domain $\Omega=[-1,1]^2$ is divided into four subdomains $\Omega_1$, $\Omega_2$, $\Omega_3$ and $\Omega_4$ (from the inside to the outside) by the following three closed curves:
\begin{align*}
    \Gamma_1&=\left\{(x, y)\;|\; x^2+y^2=0.2^2\right\}, \\
    \Gamma_2&=\left\{(x, y)\;|\; x^2+y^2=\left(0.5-0.1 \cos (5 \theta)\right)^2\right\}, \\
    \Gamma_3&=\left\{(x, y)\;|\; x^2+y^2=0.8^2\right\},
\end{align*}
as plotted in the left of \autoref{fig:multi-interf-domain-exasol-abserr}, where $\theta\in [0,2\pi)$ is the angle between the positive x-axis and the segment connecting the point $(x,y)$ and the origin. The problem is formulated as

\begin{equation}\label{equ:multi-interf}
    \left\{\begin{aligned}
        -\nabla \cdot(\beta \nabla u) & = f, & & (x,y) \in \Omega_1 \cup \Omega_2 \cup \Omega_3 \cup \Omega_4, \\
        [u] & =h_1, & & (x,y) \in \Gamma_i,\ \ i=1,2,3, \\
        [\beta \nabla u \cdot \boldsymbol{n}_i] & =h_2, & & (x,y) \in \Gamma_i,\ \ i=1,2,3, \\
        u & =g, & & (x,y) \in \partial \Omega,
    \end{aligned}\right.
\end{equation}
where $\beta(x,y)=\sum_{k=1}^{4}\beta_k \chi_{\Omega_k}$. The corresponding exact solution is given by
\begin{equation}\label{equ:multi-interf-exasol}
    u(x, y)=
    \begin{cases}
        \cos y+1.8, &(x, y) \in \Omega_1, \\
        e^x+1.3,      &(x, y) \in \Omega_2, \\
        \sin x+0.5, &(x, y) \in \Omega_3, \\
        -x+\ln  (y+2), &(x, y) \in \Omega_4,
    \end{cases}
\end{equation}
as shown in the middle of \autoref{fig:multi-interf-domain-exasol-abserr}. {The diffusion coefficient $\beta_1=1$ is fixed and  different values  will be selected for $\beta_2, \beta_3, \beta_4$ to reflect various ratio contrast (i.e., ${\beta_k}/{\beta_{k+1}}$ for $k=1,2,3$) cases.} We use four subdomain TransNets for the  Multi-TransNet method with totally 400 hidden-layer neurons and the translating and scaling parameters $\left\{(\boldsymbol{x}^{(k)}_{c}, R_i)\right\}_{k=1}^4$, where $\boldsymbol{x}^{(k)}_{c}=(0.0, 0.0)$ for $k=1,2,3,4$, and $R_1=0.35$, $R_2=0.75$, $R_3=1.0$, $R_4=1.6$, to solve the problem \eqref{equ:multi-interf}. We uniformly sample the same number of points as those in \cite{li2023local}, i.e., $N_f = 3000$ in the interior, $N_g=500$ on the boundary, and $N_{\Gamma,i}=500, i=1,2,3$ on the interfaces, for training.  

\begin{figure}[!ht]
	\centering
	\includegraphics[width=.85\textwidth]{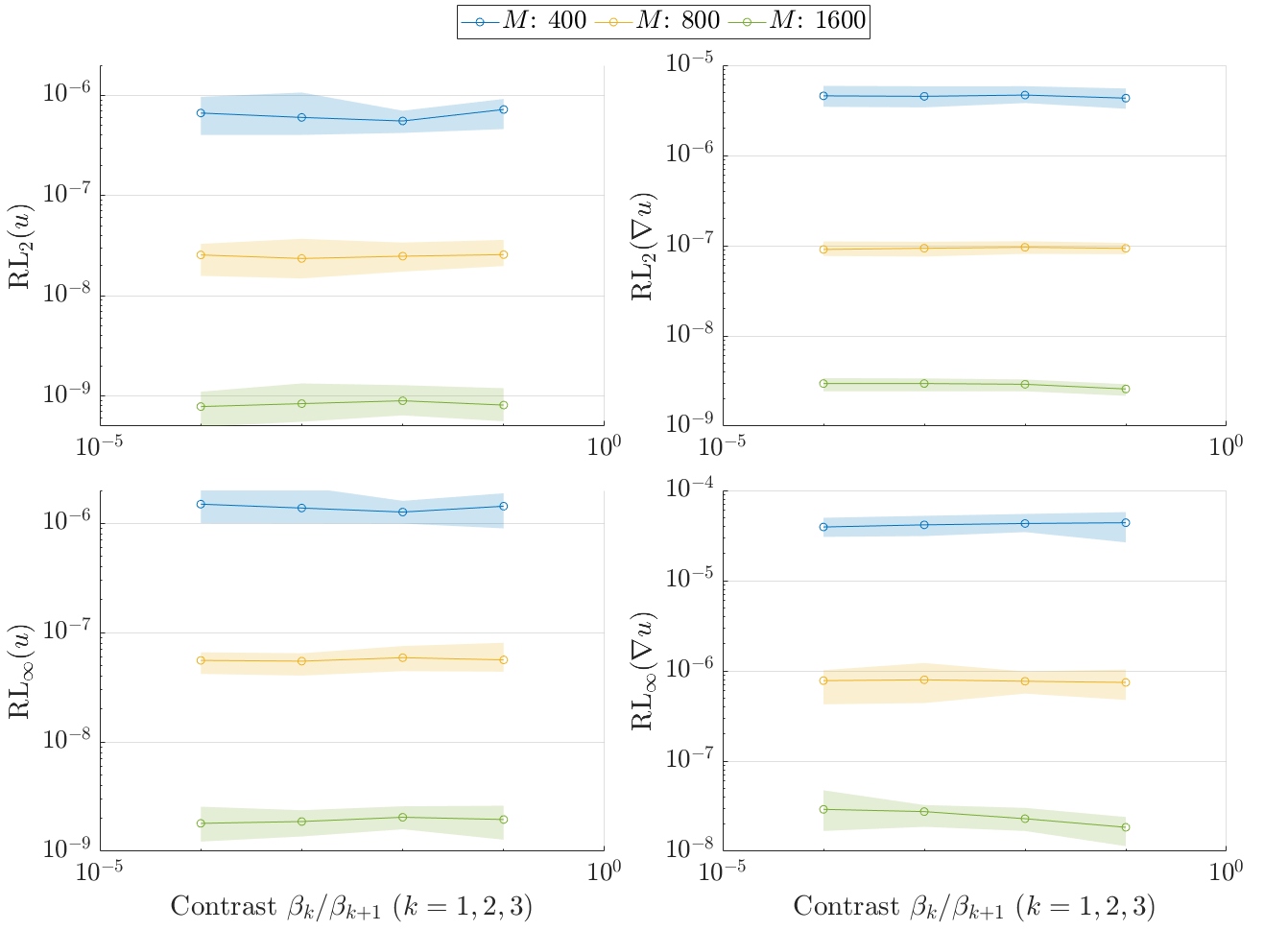}
	\vspace{-2mm}
	\caption{The relative L$_2$ and L$_\infty$ errors of the numerical solutions (left column) and their gradients (right column) produced by the Multi-TransNet with different total numbers of hidden-layer neurons ($M= 400, 800, 1600$) for the 2D diffusion interface problem \eqref{equ:multi-interf} with the piecewise constant diffusion coefficient under different contrasts {$\beta_k/\beta_{k+1} = 10^{-1}, 10^{-2}, 10^{-3}, 10^{-4}$, $k=1,2,3$}.}
	\label{fig:2D-multi-interf-accu-contrast-M}
\end{figure}

The case of diffusion coefficients $\beta_k=10^{k-1}$ ($k=2,3,4$) is first taken. The right of \autoref{fig:multi-interf-domain-exasol-abserr}  shows the pointwise   absolute  errors of the numerical solution produced by the Multi-TransNet method with totally $M=\sum_{k=1}^4 M_k=400$ hidden-layer neurons under the training loss-based optimization strategy for optimal  shape parameters $\{\gamma_k\}_{k=1}^4$. Meanwhile, an estimated value 4.1474e-2 is obtained for the empirical constant $C$ in the extended empirical formula \eqref{222}. Afterward, this empirical constant is employed to estimate  appropriate shape parameters for the Multi-TransNet method with other numbers of hidden-layer neurons, which will also be used to solve other cases of diffusion coefficient settings $\beta_k/\beta_{k+1} = 10^{-1}, 10^{-2}, 10^{-3}, 10^{-4}$ for  $k=1,2,3$. The corresponding numerical results are shown in \autoref{fig:2D-multi-interf-accu-contrast-M}. It is clearly observed that the relative L$_2$ and L$_\infty$ errors of the Multi-TransNet solutions and gradients rapidly decay as the total number of hidden-layer neurons gradually doubles from $M=400,800$ to $M=1600$, finally attaining a level of $O(10^{-9})$ for the solution $u$ and $O(10^{-8})$ for the gradient $\nabla u$. Moreover, all the results of the Multi-TransNet almost remain flat for the diffusion coefficients with low to high contrasts.

Finally, we compare the performance of the Multi-TransNet method with totally 800 hidden-layer neurons with the Local Randomized Neural Network (LRNN) \cite{li2023local} method with 1280 hidden-layer neurons for  this example with different diffusion coefficient cases $(\beta_1, \beta_2, \beta_3, \beta_4)$, and the results are reported in \autoref{tab:multi-interf}. It is evidently noticed that the Multi-TransNet remarkably outperforms the LRNN.

\begin{table}[!ht]
	\centering
	\caption{Comparison of the performance of the Multi-TransNet and the LRNN \cite{li2023local} for the 2D diffusion problem with multiple interfaces \eqref{equ:multi-interf} under different combinations of diffusion coefficients. Note that \#DOFs = $M$ for the Multi-TransNet, with the same number of constraints as those in LRNN, in this problem.}
	\vspace{0.2cm}
	\label{tab:multi-interf}
	{\footnotesize
		\begin{tabular}{ccrc}
				\toprule
				$(\beta_1, \beta_2, \beta_3, \beta_4)$ & Method & \#DOFs & RL$_2$ ($u$) \\
				\midrule
				\multirow{2}{*}{$(1, 10^2, 10^4, 10^6)$} 
				& LRNN          & 1,280 & 4.64e-5 \\
				\cmidrule{2-4}
				& Multi-TransNet & 800 & 2.41e-8 \\
				\cmidrule{1-4}
				\multirow{2}{*}{$(10^{-6}, 10^{-4}, 10^{-2}, 1)$} 
				& LRNN          & 1,280 & 4.44e-6 \\
				\cmidrule{2-4}
				& Multi-TransNet &  800 & 2.25e-8 \\
				\cmidrule{1-4}
                \multirow{2}{*}{$(1, 10^{-2}, 10^{-4}, 10^{-6})$} 
				& LRNN          & 1,280 & 7.68e-7 \\
				\cmidrule{2-4}
				& Multi-TransNet &  800 & 5.75e-8 \\
				\cmidrule{1-4}
				\multirow{2}{*}{$(10^{6}, 10^{4}, 10^{2}, 1)$} 
                & LRNN          & 1,280 & 7.16e-5 \\
				\cmidrule{2-4}
				& Multi-TransNet &  800 & 5.68e-8 \\
				\bottomrule
			\end{tabular}
		}
\end{table}

\subsubsection{A 3D elasticity interface problem with an ellipsoidal interface}\label{ssec:elastic-interf}
\begin{figure}[!ht]
	\centering
	\includegraphics[width=\linewidth]{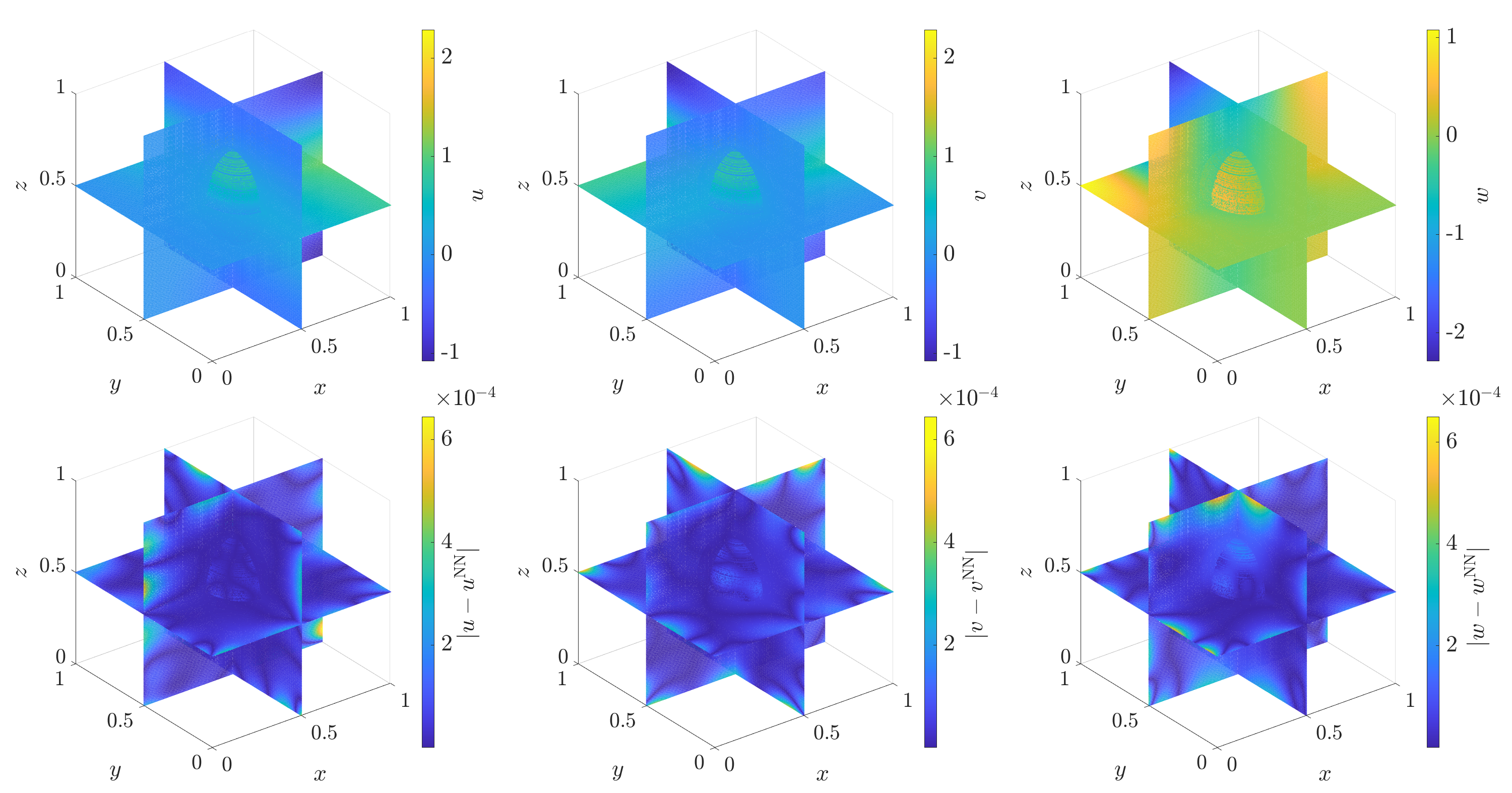}
	\vspace{-6mm}
	\caption{Top row: the exact solution of the 3D elasticity interface problem \eqref{equ:ex1} on the interface and three centered intersecting planes; Bottom row: the pointwise absolute errors of the numerical solution, when the Lam\'e parameters $(\lambda_1, \lambda_2)=(\mu_1, \mu_2)=(1,10)$, produced by the Multi-TransNet method with totally 500 hidden-layer neurons. From left to right: $u$, $v$ and $w$.}
	\label{fig:elastic-interf-exasol-abserr}
\end{figure}

In this example, we consider a three-dimensional elasticity interface problem with an ellipsoidal interface  \cite{zhang2022high}. The domain  $\Omega=[0,1]^3$ is partitioned into two subdomains $\Omega_1$ (inside) and $\Omega_2$ (outside) by the following ellipsoidal surface
\[
\Gamma:\ \frac{(x-0.5)^2}{a^2}+\frac{(y-0.5)^2}{b^2}+\frac{(z-0.5)^2}{c^2}=1,
\]
where $a=0.15, b=0.2$ and $c=0.25$. The problem is described by
\begin{equation}\label{equ:ex1}
    \left\{\begin{aligned}
        -\nabla \cdot \boldsymbol{\sigma}(\boldsymbol{u}) & =\boldsymbol{f}, & & (x,y,z) \in \Omega_1\cup \Omega_2, \\
        \left[\boldsymbol{u} \right] & =\boldsymbol{h}_1, & & (x,y,z) \in \Gamma, \\
        \left[\boldsymbol{\sigma}(\boldsymbol{u}) \boldsymbol{n}\right] & =\boldsymbol{h}_2, & & (x,y,z) \in \Gamma, \\
        \boldsymbol{u} & =\boldsymbol{g}, & & (x,y,z) \in \partial \Omega,
    \end{aligned}\right.
\end{equation}
where the displacement $\boldsymbol{u}=\left(u, v, w\right)^\top$, the stress tensor $\boldsymbol{\sigma}(\boldsymbol{u}) = \lambda\left(\nabla \cdot \boldsymbol{u}\right)\boldsymbol{I} + 2\mu \boldsymbol{\epsilon}\left(\boldsymbol{u}\right)$, the strain tensor $\boldsymbol{\epsilon}\left(\boldsymbol{u}\right) = \frac{1}{2}\left(\nabla \boldsymbol{u} + \nabla \boldsymbol{u}^\top\right)$, the Lam\'e parameters $\lambda(x,y,z) = \lambda_k$ and $\mu(x,y,z) = \mu_k$ for $(x,y,z)\in \Omega_k$, $k=1,2$. The corresponding exact solution is
chosen as
\begin{align*}
    & u(x,y,z)= \begin{cases}-\cos \left(x^2\right) \mathrm{e}^{-y^2} \sin (2 \pi z), & (x,y,z) \in \Omega_1 \\
        -\sin \left(x^2\right) \mathrm{e}^{y^2} \cos (2 \pi z), & (x,y,z) \in \Omega_2,\end{cases} \\
    & v(x,y,z)= \begin{cases}-\cos \left(y^2\right) \mathrm{e}^{-x^2} \sin (2 \pi z), & (x,y,z) \in \Omega_1 \\
        -\sin \left(y^2\right) \mathrm{e}^{x^2} \cos (2 \pi z), & (x,y,z) \in \Omega_2,\end{cases} \\
    & w(x,y,z)= \begin{cases}\cos \left(y^2\right) \mathrm{e}^{-z^2} \sin (2 \pi x), &\;\; (x,y,z) \in \Omega_1 \\
        \sin \left(y^2\right) \mathrm{e}^{z^2} \cos (2 \pi x), &\;\; (x,y,z) \in \Omega_2,\end{cases}
\end{align*}
which on the interface and three centered intersecting planes is plotted in the top row of \autoref{fig:elastic-interf-exasol-abserr}. {The Lam\'e parameters $\lambda_1=1$ and $\mu_1=1$ are fixed and other values will be set to reflect various contrast (i.e., $\lambda_1/\lambda_2$ and $\mu_1/\mu_2$) cases.}
Two subdomain TransNets are used for Multi-TransNet to solve the problem \eqref{equ:ex1}, and a combination of  translating and scaling parameters  is set to be $(\boldsymbol{x}^{(1)}_{c}, R_1)=((0.5, 0.5, 0.5), 0.35)$ for $\Omega_1$ and $(\boldsymbol{x}^{(2)}_{c}, R_2)=((0.5, 0.5, 0.5), 1.0)$ for $\Omega_2$. The numbers of sampling points for training are fixed to $N_f=29^3$ (in the interior), $N_g=31^2\times 6$ (on the boundary) and $N_\Gamma=60\times 30$ (on the interface), with an overall spacing of about 0.03. 

\begin{figure}[!ht]
	\centering
	\includegraphics[width=\textwidth]{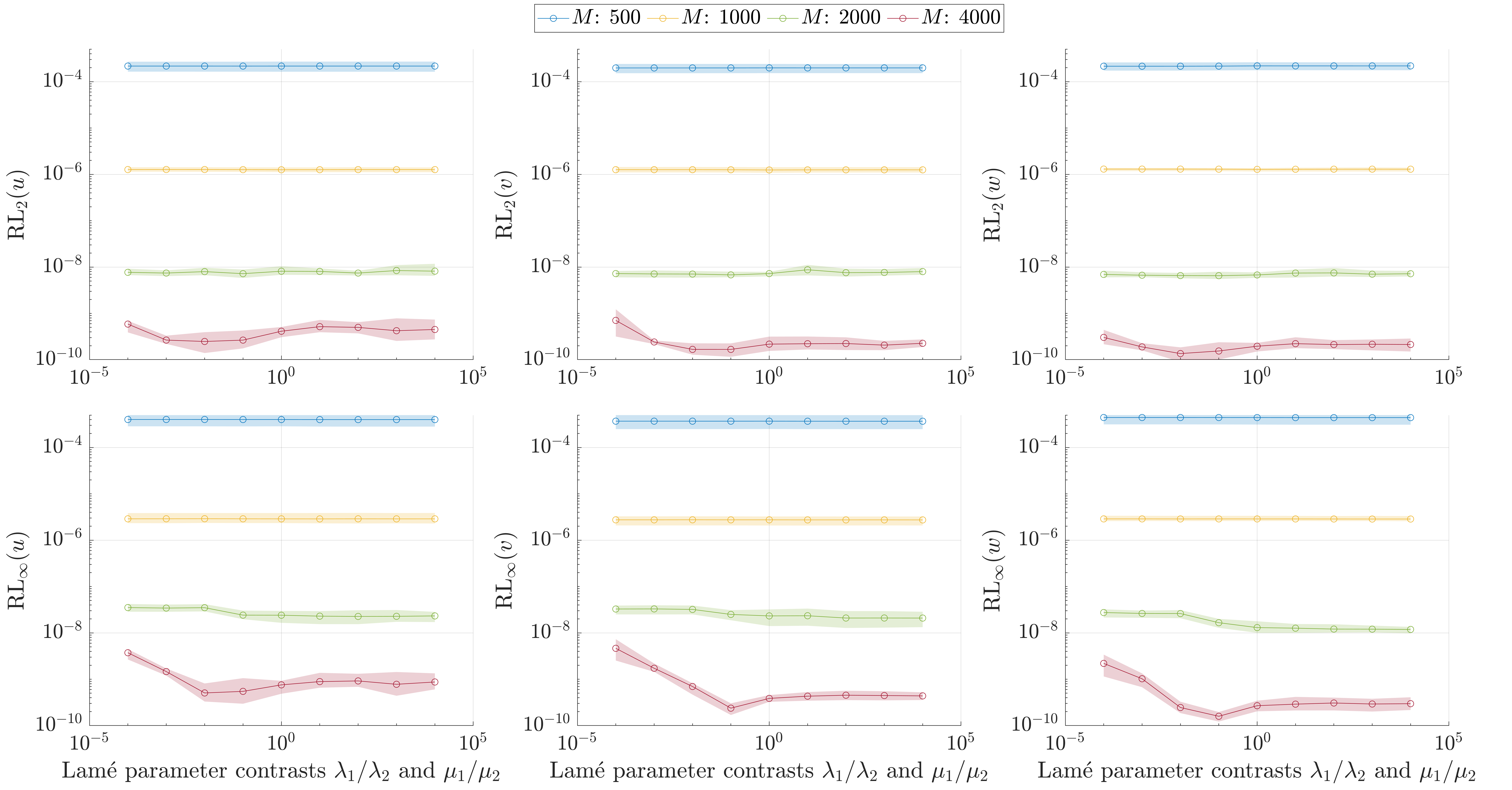}
	\vspace{-6mm}
	\caption{The relative L$_2$ (top row) and L$_\infty$ (bottom row) errors of the numerical solutions produced by the Multi-TransNet method with different total  numbers of hidden-layer neurons ($M = 500, 1000, 2000, 4000$) for the 3D elasticity interface problem \eqref{equ:ex1} with {different Lam\'e parameter contrasts}. From left to right:  $u$, $v$  and $w$.}
	\label{fig:elastic-accu-contrast-M}
\end{figure}

\begin{figure}[!ht]
	\centering
	\includegraphics[width=0.55\textwidth]{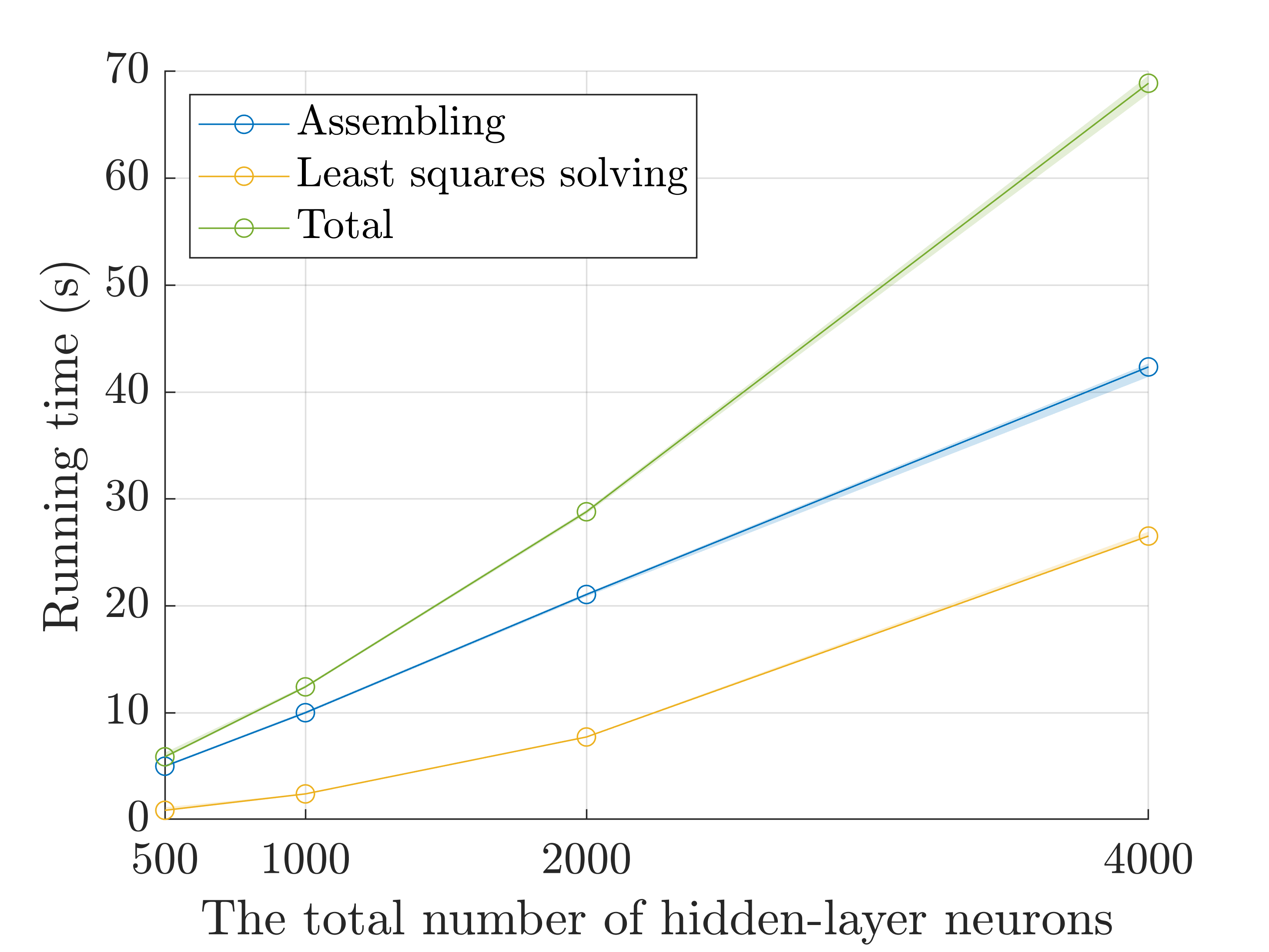}
	\vspace{-2mm}
	\caption{The average running times (measured in seconds) of the Multi-TransNet method with different total numbers of hidden-layer neurons ($M = 500, 1000, 2000, 4000$) for the 3D elasticity interface problem \eqref{equ:ex1}.}
	\label{fig:elastic-rt}
\end{figure}

We first tackle the low-contrast case of the Lam\'e parameters {$\lambda_2=\mu_2=10$} by using the Multi-TransNet with totally 500 hidden-layer neurons and the training loss-based optimization strategy for optimal shape parameters $(\gamma_1,\gamma_2)$. The pointwise absolute errors of the numerical solution are illustrated in the bottom row of \autoref{fig:elastic-interf-exasol-abserr}. We observe that the absolute errors on the interface and intersecting planes reach the level of approximately 6e-4. An estimated value 8.1606e-2 is obtained for the empirical constant $C$ in the extended  empirical formula \eqref{111} for the Multi-TransNet. 
Then the empirical formula-based prediction strategy is employed to determine appropriate shape parameters for the Multi-TransNet with other numbers of hidden-layer neurons, which will be used to solve different cases of the Lam\'e parameter contrasts $\lambda_1/\lambda_2$ and $\mu_1/\mu_2$.
\autoref{fig:elastic-accu-contrast-M} shows the numerical results for $\lambda_1/\lambda_2=\mu_1/\mu_2$ increasing from $10^{-4}$ to $10^{4}$ by a factor of 10 each time. From that, it is easy to see the relative L$_2$ and L$_\infty$ errors of the Multi-TransNet solutions quickly decrease when the total number of hidden-layer neurons gradually doubles from $M=500,1000,2000$ to $M=4000$, and finally reach a level of $O(10^{-10})$ for all variables $u$, $v$ and $w$. Furthermore, we observe that 
both relative  errors mostly  remain flat along the significant changes of the magnitude of the Lam\'e parameter contrasts, except that the relative L$_\infty$ errors exhibit some fluctuations when $M$ is taken as 4000. \autoref{fig:elastic-rt} reports the average running times of the Multi-TransNet with respect to different total numbers of hidden-layer neurons. We again find that the system assembling time doubles along with the doubling of the number of  hidden-layer neurons but the least squares solving time increases faster.
 
Finally, we compare the performance of  the Multi-TransNet method with  the high-order Non-symmetric Interface Penalty Finite Element Method (NIPFEM) \cite{zhang2022high}, which was recently proposed for solving the 3D elasticity interface problem with the high-contrast Lam\'e parameters $(\lambda_1, \lambda_2)=(\mu_1, \mu_2)=(1,1000)$.  The numerical results are reported in \autoref{tab:elastic-interf} and it is easy to find that   Multi-TransNet easily and significantly beats  the NIPFEM in terms of  the number of DOFs and accuracy.  

\begin{table}[!ht]
	\centering
	\caption{Comparison of the performance of the Multi-TransNet and the high-order NIPFEM \cite{zhang2022high} for the 3D elasticity interface problem \eqref{equ:ex1} with the high-contrast Lam\'e parameters $(\lambda_1, \lambda_2)=(\mu_1,  \mu_2)=(1,1000)$. Note that $\#$DOFs = 3$M$ for the  Multi-TransNet in this problem.}
	\label{tab:elastic-interf}
	\vspace{0.2cm}
	{\footnotesize
	\begin{tabular}{crrl}
		\toprule
		Method & \#Constraints & \#DOFs & RL$_2$ ($\boldsymbol{u}$) \\
		\midrule
		\multirow{3}{*}{NIPFEM} 
		& 23,550 & 23,550   & 2.964e-3 \\ 
		& 176,694 & 176,694  & 1.382e-4 \\ 
		& 1,369,446 & 1,369,446 & 6.959e-6 \\ 
		\cmidrule{1-4}
		\multirow{3}{*}{Multi-TransNet} & \multirow{3}{*}{101,262}
		& 1,500 & 2.115e-4 \\
		& & 3,000 & 1.282e-6 \\
		& & 6,000 & 7.073e-9 \\
		\bottomrule
	\end{tabular}
	}
\end{table}

\subsubsection{A 3D diffusion interface problem with a convoluted immersed interface}\label{ssec:3D-immer}
\begin{figure}[!t]
    \centering
    \includegraphics[width=0.48\linewidth]{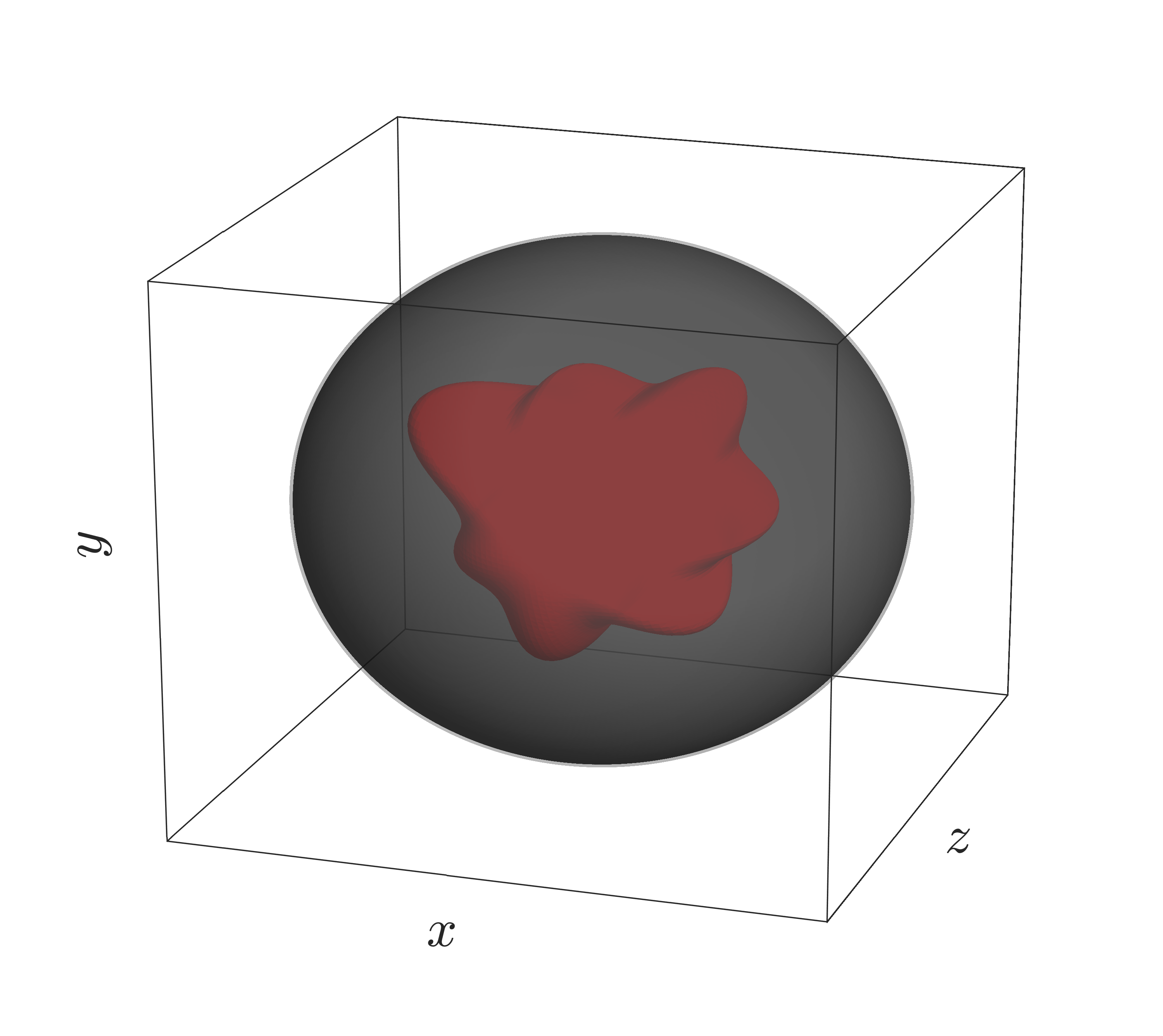}
      \includegraphics[width=0.48\linewidth]{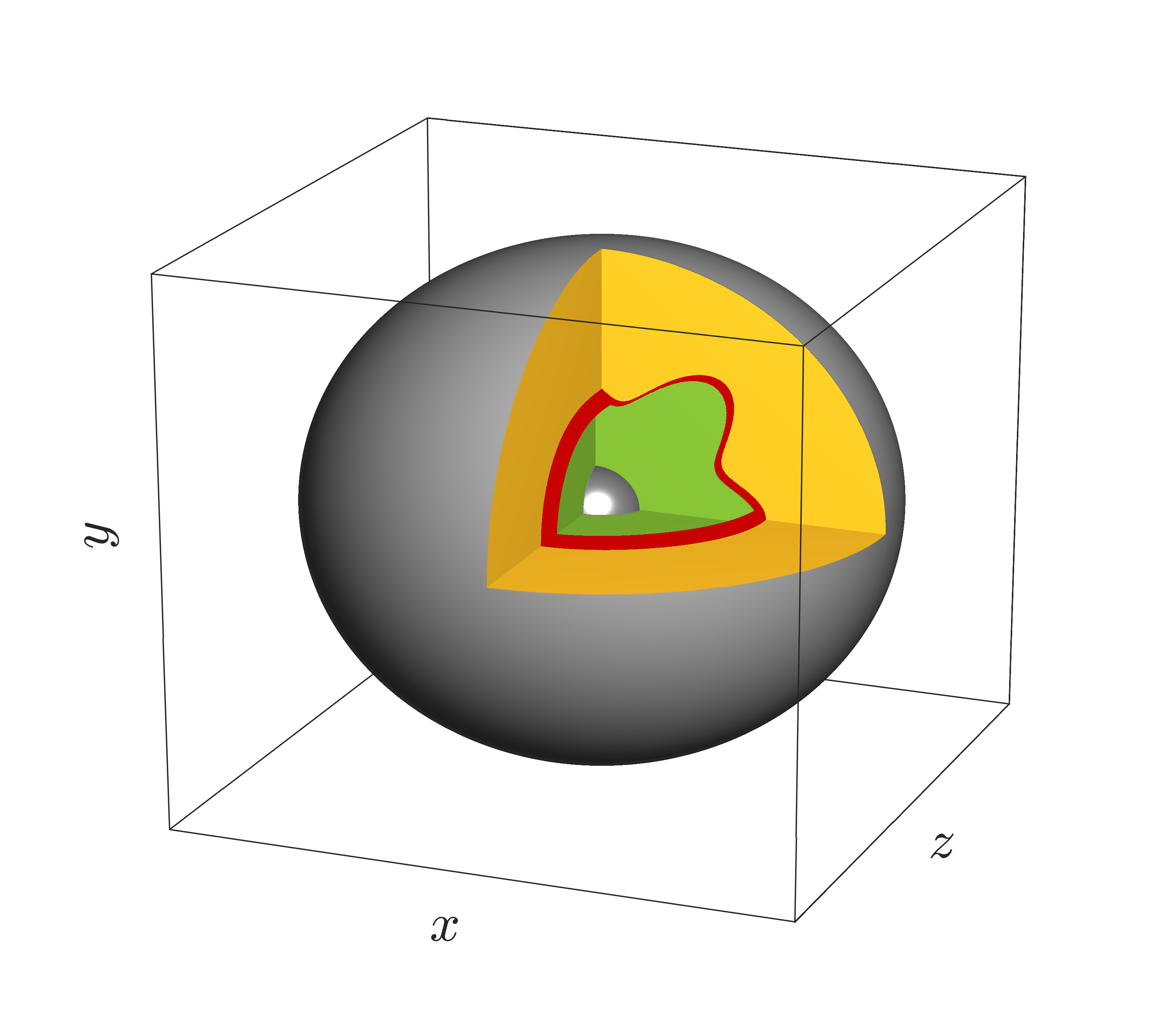}
    \vspace{-4mm}
    \caption{Profiles of the convoluted interface \eqref{equ:3D-immer-interf} immersed in a spherical shell-shaped domain for the 3D diffusion interface problem \eqref{equ:3D-immer}.}
    \label{fig:3D-immer-dom-interf}
\end{figure}

\begin{figure}[!t]
	\centering
	\includegraphics[width=\linewidth]{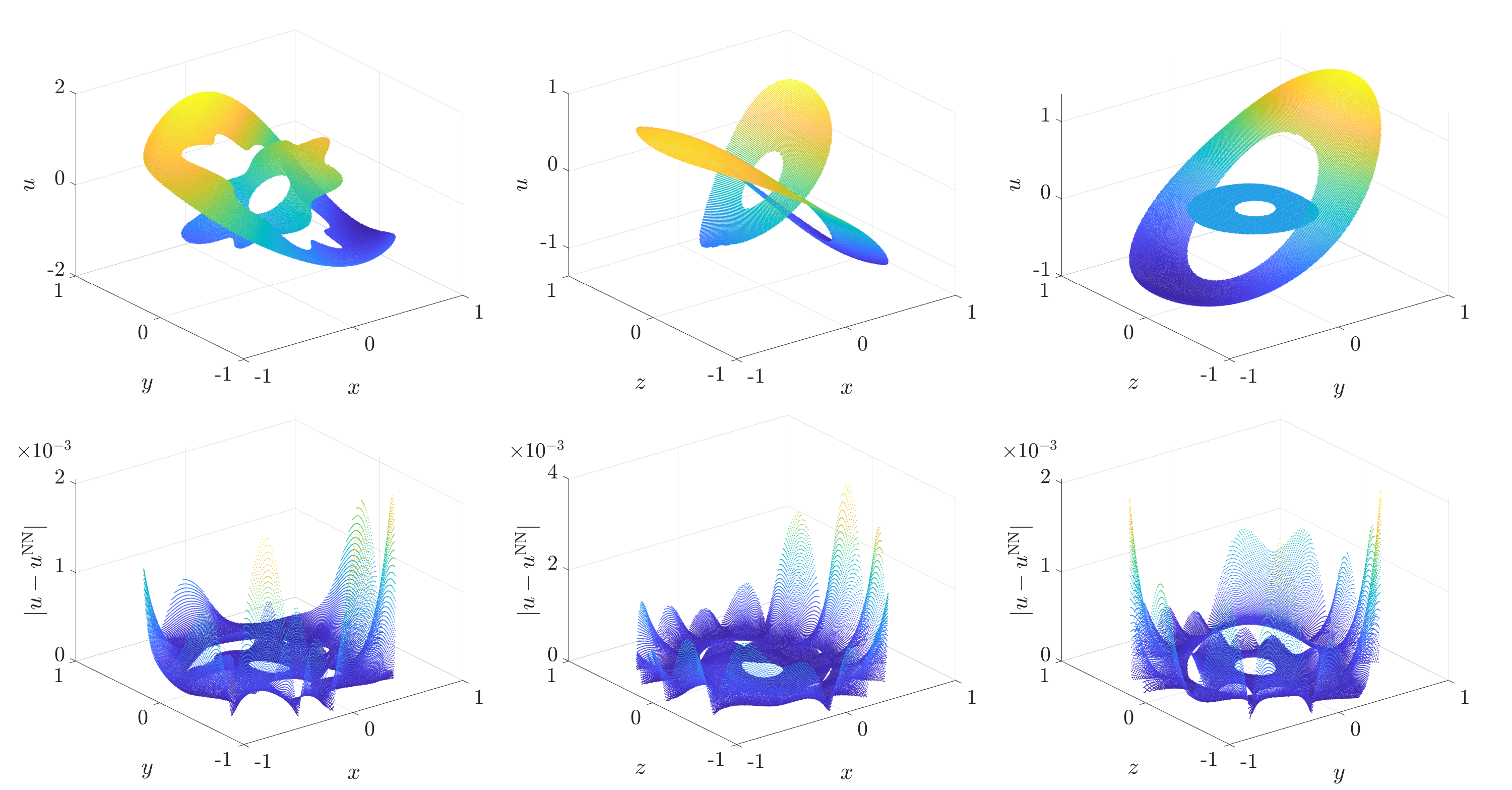}
	\vspace{-6mm}
	\caption{Top row: the exact solution of the 3D diffusion interface problem \eqref{equ:3D-immer} on the three coordinate planes; Bottom row: the pointwise absolute errors of the numerical solution on the three coordinate planes  (in the case of the spatially varying diffusion coefficient \eqref{equ:3D-immer-varying-coef}) produced by the Multi-TransNet with totally 250 hidden-layer neurons. From left to right: $xy$-plane, $xz$-plane and $yz$-plane.}
	\label{fig:3D-immer-varying-coef-exasol-abserr}
\end{figure}

In this example, we consider a 3D diffusion problem  with a convoluted interface \cite{bochkov2020solving}  immersed in a spherical shell-shaped domain $\Omega = \{ (x,y,z) \,|\,\ 0.151^2 \leqslant x^2 + y^2 + z^2 \leqslant 0.911^2 \}$. The interface $\Gamma$ is implicitly defined by 
\begin{equation}\label{equ:3D-immer-interf}
    \sqrt{x^2+y^2+z^2}-r_0\left(1+\left(\frac{x^2+y^2}{x^2+y^2+z^2}\right)^2 \sum_{k=1}^3 a_k \cos \left(n_k\left(\arctan \left(\frac{y}{x}\right)-\theta_k\right)\right)\right) = 0,
\end{equation}
with the parameters
\[
r_0=0.483, \quad
\begin{pmatrix}
a_1 \\
a_2 \\
a_3 \\
\end{pmatrix}=\begin{pmatrix}
0.1 \\
-0.1 \\
0.15
\end{pmatrix}, \quad
\begin{pmatrix}
n_1 \\
n_2 \\
n_3
\end{pmatrix}=\begin{pmatrix}
3 \\
4 \\
7
\end{pmatrix}, \text{ and }
\begin{pmatrix}
\theta_1 \\
\theta_2 \\
\theta_3
\end{pmatrix}=\begin{pmatrix}
0.5 \\
1.8 \\
0
\end{pmatrix},
\]
as plotted in \autoref{fig:3D-immer-dom-interf}. The problem is formulated as 
\begin{equation}\label{equ:3D-immer}
    \left\{\begin{aligned}
        -\nabla \cdot(\beta \nabla u) & = f, & & (x,y,z) \in \Omega_1 \cup \Omega_2, \\
        [u] & =h_1, & & (x,y,z) \in \Gamma, \\
        [\beta \nabla u \cdot \boldsymbol{n}] & =h_2, & & (x,y,z) \in \Gamma, \\
        u & =g, & & (x,y,z) \in \partial \Omega,
    \end{aligned}\right.
\end{equation}
where $\Omega_1$ denotes the interior subdomain and $\Omega_2$  the exterior subdomain.
The exact solution is chosen as
\begin{equation}\label{equ:3D-immer-exasol}
u(x, y, z) = 
\begin{cases}
\sin (2 x) \cos (2 y) \exp (z), & (x, y, z) \in \Omega_1, \\
\left(16\left(\frac{(y-x)}{3}\right)^5-20\left(\frac{(y-x)}{3}\right)^3+5\left(\frac{(y-x)}{3}\right)\right) \log (x+y+3) \cos (z), & (x, y, z) \in \Omega_2,
\end{cases}
\end{equation}
which on three coordinate planes ($xy$, $xz$ and $yz$) are shown in the top row of \autoref{fig:3D-immer-varying-coef-exasol-abserr}.

We employ two subdomain TransNets with the respective translating and scaling parameters $(\boldsymbol{x}^{(1)}_{c}, R_1)=((0.0, 0.0, 0.0), 0.75)$ and $(\boldsymbol{x}^{(2)}_{c}, R_2)=((0.0, 0.0, 0.0), 1.0)$ for Multi-TransNet to solve the problem \eqref{equ:3D-immer}. The numbers of sampling points for training are set to be around $N_f=12000$ (in the interior), $N_g=2000$ (on the boundaries, 200 on the interior one and 1800 on the exterior one) and $N_\Gamma=750$ (on the interface). We first carry out the convergence test for the case of a spatially varying diffusion coefficient, and then test for the case of a piecewise constant diffusion coefficient.

\begin{figure}[!t]
	\centering
	\begin{minipage}[t]{.48\textwidth}
		\centering
		\includegraphics[width=\linewidth]{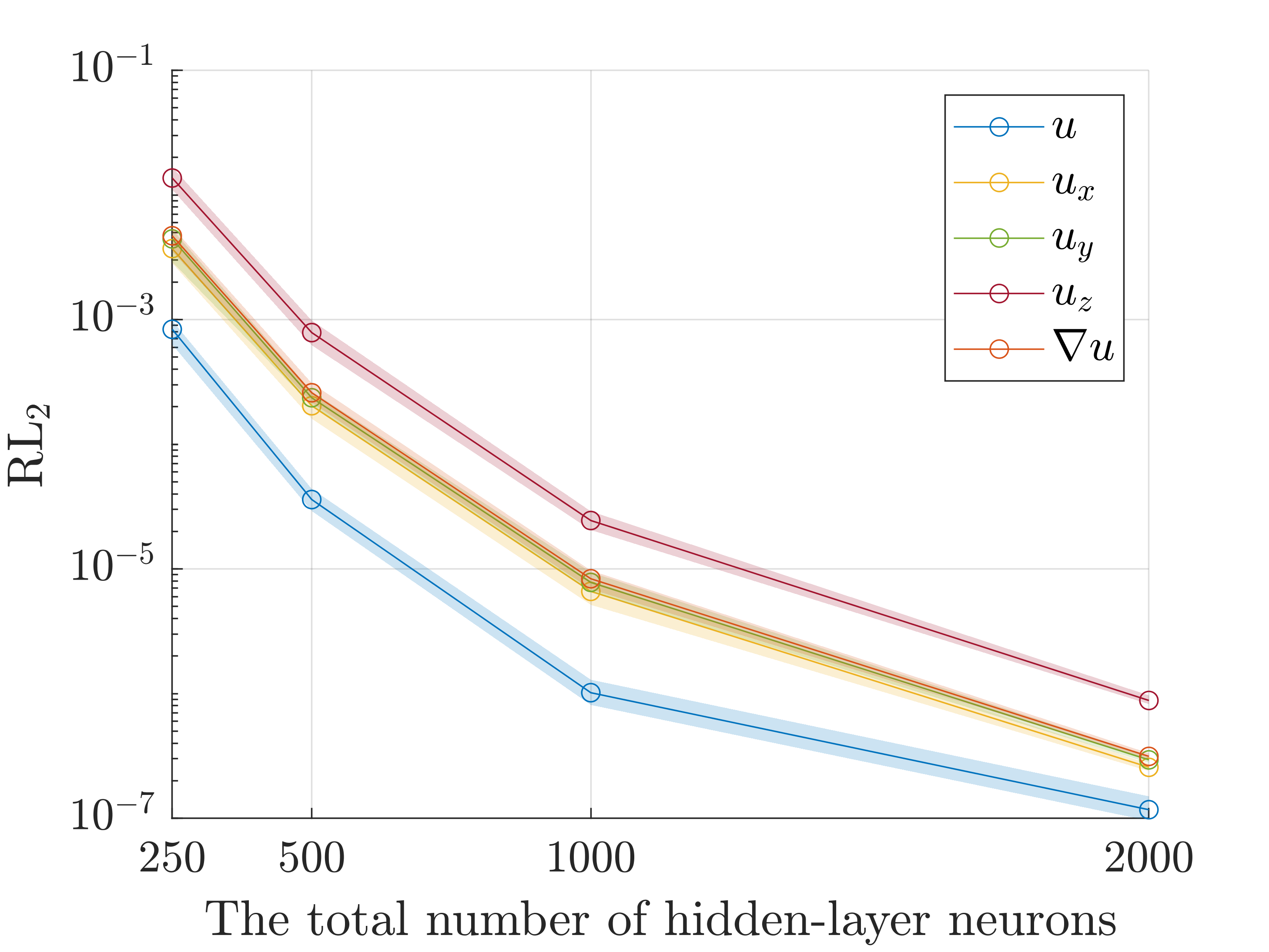}	
	\end{minipage}
	\begin{minipage}[t]{.48\textwidth}
		\centering
		\includegraphics[width=\linewidth]{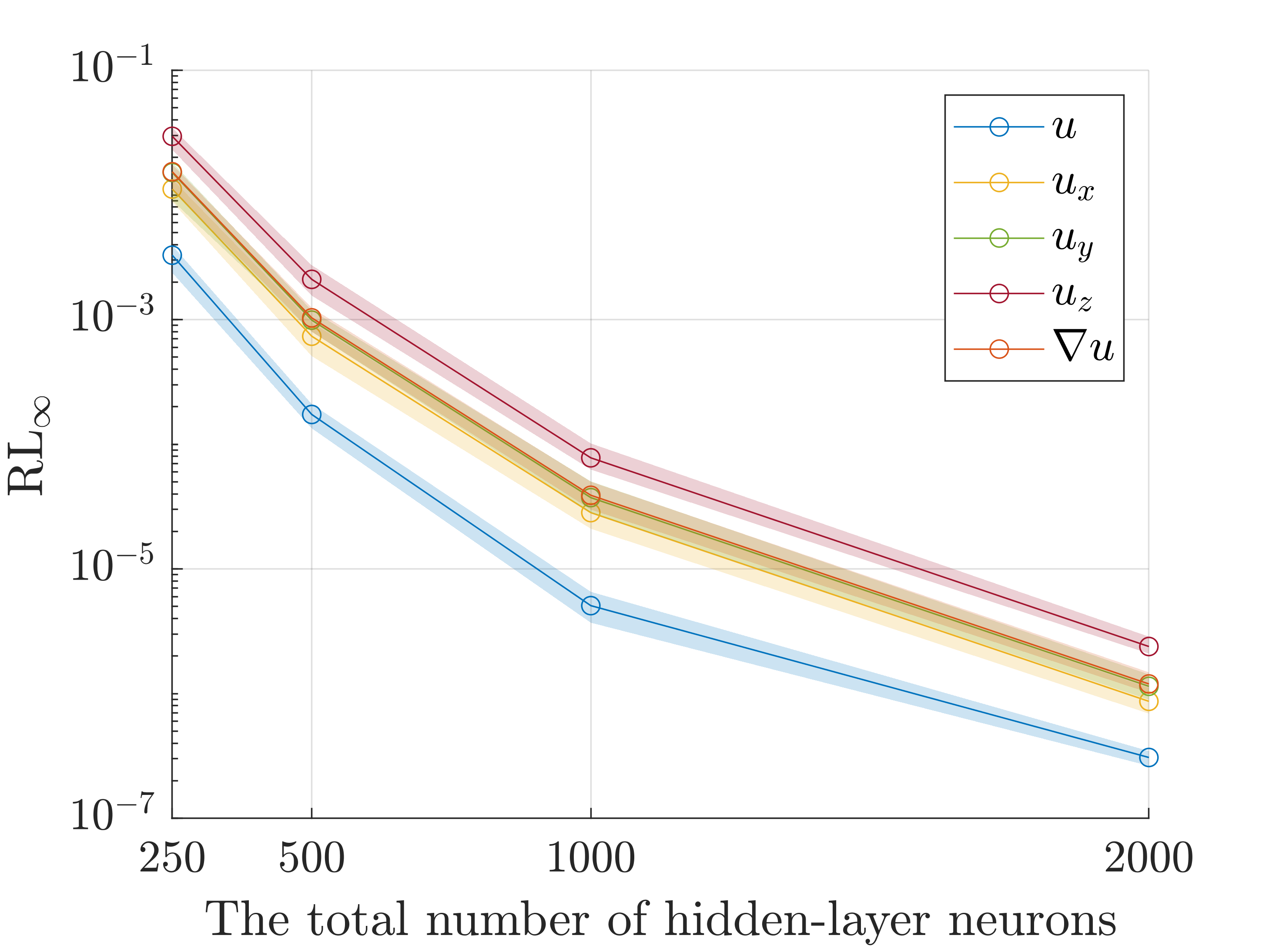}
	\end{minipage}
	\vspace{-2mm}
    \caption{The relative L$_2$ (left) and L$_\infty$ (right) errors of numerical solutions and their partial derivatives and gradients produced by the Multi-TransNet  for the 3D diffusion interface problem \eqref{equ:3D-immer} with the spatially varying diffusion coefficient \eqref{equ:3D-immer-varying-coef}.}\label{fig:3D-immer-varying-coef-err-time}
\end{figure}

\paragraph*{Case I: With a spatially varying diffusion coefficient}
The diffusion coefficient is defined by 
\begin{equation}\label{equ:3D-immer-varying-coef}
\beta(x, y, z)= 
\begin{cases}
10\left(1+\frac{1}{5} \cos (2 \pi(x+y)) \sin (2 \pi(x-y)) \cos (z)\right), & (x, y, z) \in \Omega_1, \\
1, & (x, y, z) \in \Omega_2.
\end{cases}
\end{equation}
We first use a total of 250 hidden-layer neurons for the Multi-TransNet along with the training loss-based optimization strategy for searching optimal shape parameters $(\gamma_1,\gamma_2)$ and the  pointwise absolute errors  of the numerical solution on three coordinate planes  are shown  in the bottom row of \autoref{fig:3D-immer-varying-coef-exasol-abserr}. The whole golden-section searching process takes about 1.11 seconds. With the estimated value 4.4180e-2 for the empirical constant $C$ in the empirical formula \eqref{111},  the empirical
formula-based prediction strategy is then employed to determine appropriate shape parameters for the Multi-TransNet with totally 500, 1000 and 2000 hidden-layer neurons to solve the problem \eqref{equ:3D-immer}. The convergence results on the relative errors of the Multi-TransNet solutions and their partial derivatives and gradients for the solution process are shown in \autoref{fig:3D-immer-varying-coef-err-time}. We observe that all relative L$_2$ and L$_\infty$ errors rapidly and steadily decay along the doubling of the total number of hidden-layer neurons. 

Furthermore, we compare the accuracy and running times of the Multi-TransNet method with the cusp-capturing PINN method \cite{tseng2023cusp}, which was recently proposed for solving the elliptic interface problem.  The numerical results are reported in \autoref{tab:3D-immer-varying-coef}, from which it is observed that the running times of the Multi-TransNet method are significantly less than those of the cusp-capturing PINN when achieving nearly the same level of relative errors. This is  because the Multi-TransNet method can efficiently solve the parameters of the output layer using the well-established least-squares techniques.

\begin{table}[!ht]
	\centering
	\caption{Comparison of the accuracy and running times of the Multi-TransNet method and the cusp-capturing PINN method \cite{tseng2023cusp} for the 3D diffusion interface problem \eqref{equ:3D-immer} with the spatially varying diffusion coefficient \eqref{equ:3D-immer-varying-coef}. Note that $M$ and \#DOFs represent the number of the hidden-layer neurons and the number of trainable parameters, respectively.}
	\label{tab:3D-immer-varying-coef}
	\vspace{0.2cm}
	{\footnotesize
	\begin{tabular}{ccllr}
		\toprule
		\multirow{2}{*}{\begin{tabular}[c]{@{}l@{}}Method\end{tabular}} & 
        \multirow{2}{*}{\begin{tabular}[c]{@{}l@{}}($M$, \#DOFs)\end{tabular}} & 
        \multirow{2}{*}{\begin{tabular}[c]{@{}l@{}}RL$_2$($u$)\end{tabular}} & 
        \multirow{2}{*}{\begin{tabular}[c]{@{}l@{}}RL$_\infty$($u$)\end{tabular}} &
        \multirow{2}{*}{\begin{tabular}[c]{@{}l@{}}Running\\time (s)\end{tabular}} \\
        &&&& \\
		\midrule
		\multirow{3}{*}{cusp-capturing PINN} 
		& (25, 150)  & 1.10e-3 & 1.86e-3 & 55.06 \\ 
		& (50, 300)  & 2.01e-5 & 5.11e-5 & 75.12 \\ 
		& (100, 600) & 2.10e-6 & 5.66e-6 & 91.23 \\ 
		\cmidrule{1-5}                        
		\multirow{3}{*}{Multi-TransNet}
		& (500, 500)  & 3.60e-5 & 1.73e-4 & 0.33 \\ 
		& (1000, 1000)& 1.02e-6 & 5.08e-6 & 0.76 \\ 
		& (2000, 2000)& 1.17e-7 & 3.07e-7 & 1.70 \\ 
		\bottomrule
	\end{tabular}
	}
\end{table}

\begin{figure}[!htb]
	\centering
	\includegraphics[width=.85\textwidth]{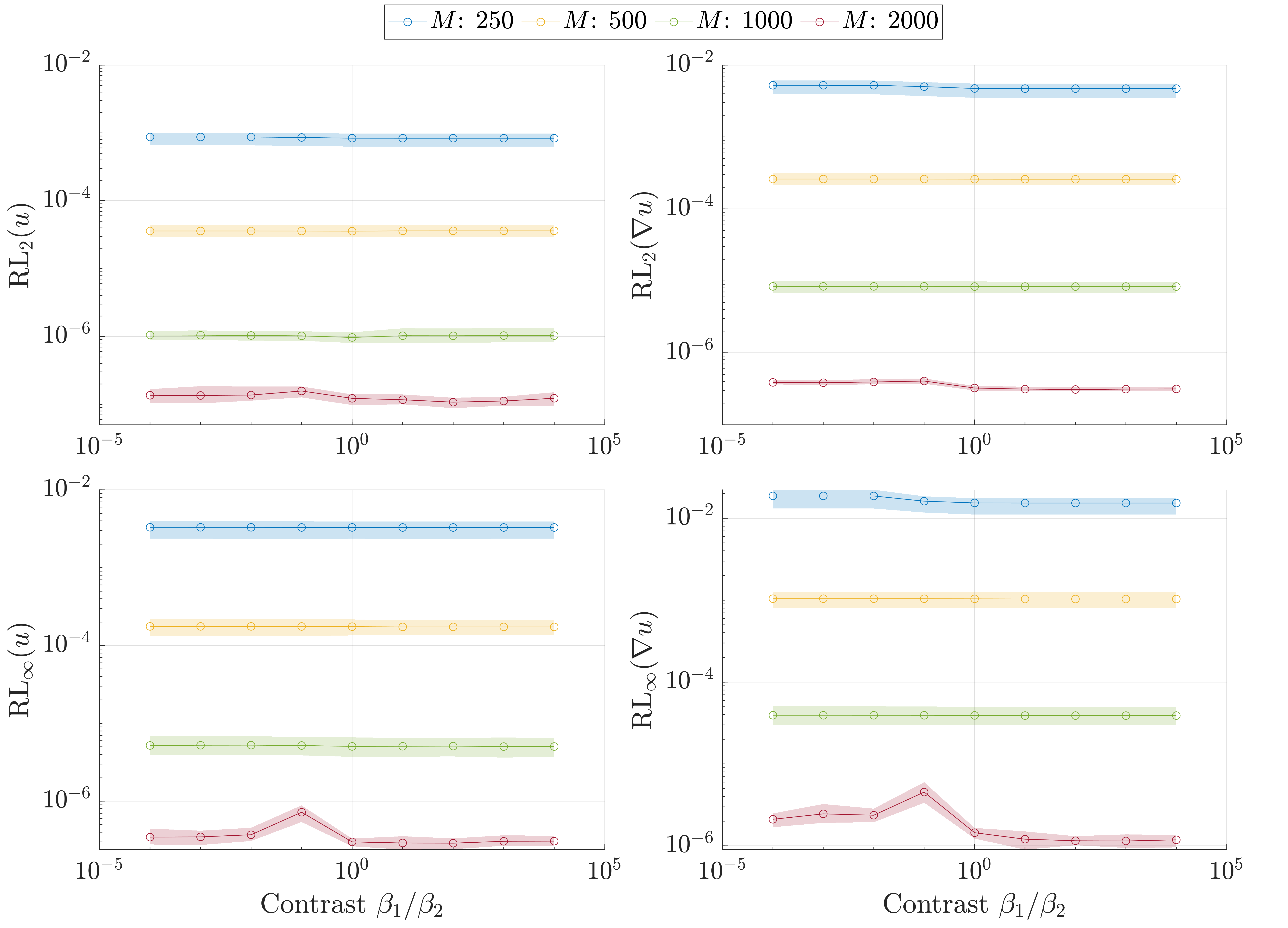}
	\vspace{-2mm}
	\caption{The relative L$_2$ and L$_\infty$ errors of the numerical solutions (left column) and their gradients (right column) produced by the Multi-TransNet with different total numbers of hidden-layer neurons ($M= 250, 500, 1000, 2000$) for the 3D diffusion interface problem \eqref{equ:3D-immer} with the piecewise constant diffusion coefficient under different contrasts.}
	\label{fig:3D-immers-interf-pwconst-coef-accu-contrast-M}
\end{figure}

\paragraph*{Case II: With a piecewise constant diffusion coefficient}
The diffusion coefficient is defined to $\beta(x,y,z) = \beta_k$ for $(x,y,z)\in\Omega_k$ ($k=1,2$). We fix $\beta_2=1$, and adopt all of the parameters of the Multi-TransNet from \emph{Case I} to solve \emph{Case II}, where the diffusion coefficient contrast $\beta_1/\beta_2$ increase from $10^{-4}$ to $10^{4}$ by a factor of 10 each time. The relative L$_2$ and L$_\infty$ errors of the Multi-TransNet solution and its gradient are illustrated in \autoref{fig:3D-immers-interf-pwconst-coef-accu-contrast-M}. It is clearly observed that for a fixed total number of hidden-layer neurons, all the results of the Multi-TransNet basically remain flat for different contrasts, except for a bump arising in the relative L$_\infty$ errors of the solution and its gradient when $M=2000$. With the doubling of the total number of hidden-layer neurons, the relative errors of both the Multi-TransNet solution and its gradient rapidly decrease.

\section{Conclusions}\label{sec:conclusions}
In this paper, we propose a novel Multi-TransNet method, which integrates multiple distinct TransNets using the nonoverlapping domain decomposition approach, to solve various types of elliptic interface problems. Specifically, we design an efficient empirical formula-based prediction strategy for automatically selecting proper shape parameters for the subdomain TransNets, greatly reducing the tuning cost. Apart from that, the globally uniform distribution of the hidden-layer neurons across the entire domain is developed, with adaptive assignment of the number of hidden-layer neurons for each subdomain TransNet. During training, the weighting parameters are adaptively determined based on a normalization strategy, which balances the terms of the loss function and enhances the robustness. Extensive ablation studies and comparisons with some recently proposed neural network methods and traditional numerical techniques for solving classic 2D and 3D elliptic interface problems energetically demonstrate that the proposed Multi-TransNet method can offer striking performance in terms of accuracy, efficiency and robustness. 
To the best of our knowledge, this is the first time the TransNet technique has been extended to solve elliptic interface problems. Some potential future work includes: (1) improving the assembling and solving efficiency of the resulting least squares problem; (2) developing more effective approaches for generating the hidden-layer neurons based on specific target domains; (3) extending the present method to dynamic interface problems.

\section*{Acknowledgement}
L. Zhu's work was partially supported by National Natural Science Foundation of
China under grant number 11871091.

\bibliographystyle{elsarticle-num} 
\bibliography{elsarticle-bibliography.bib}

\end{document}